%% file: main.tex
\documentclass{art}
\usepackage{mypreamble}

\begin{document}
\begin{frontmatter}

  \title{Bicompleteness Theorems for Team Logics with the Dual Negation}
  \runtitle{Bicompleteness Theorems for Team Logics}

  \author{\fnms{Aleksi}
    \snm{Anttila}
    \ead[label=e1]{aleksi.anttila@unipd.it}
  }
  \address{Department of Philosophy, Sociology, Education and Applied Psychology\\
    University of Padua\\
    Palazzo del Capitanio, Piazza Capitaniato 3 \\
    35139 Padua\\
    ITALY\\
    \printead{e1}\\
  }%
  \runauthor{A.~Anttila}

\begin{abstract}
The dual or game-theoretical negation $\lnot$ of independence-friendly logic (IF) and dependence logic (D) exhibits an extreme degree of semantic indeterminacy in that for any pair of sentences $\phi$ and $\psi$ of IF/D, if $\phi$ and $\psi$ are incompatible in the sense that they share no models, there is a sentence $\theta$ of IF/D such that $\phi\equiv \theta$ and $\psi\equiv \lnot \theta$ (as shown originally by Burgess in the equivalent context of the prenex fragment of Henkin quantifier logic). Together with its converse, a result of this type can be seen as an expressive completeness theorem with respect to properties of pairs: a pair (A,B) of classes of models definable in IF/D is disjoint just in case A is the class of models of some IF/D sentence $\phi$, and B is the class of models of its negation $\lnot \phi$---we formulate a notion of expressive completeness for pairs (\emph{bicompleteness}) to make this precise. We prove a number of bicompleteness theorems with respect to different classes of pairs for propositional and modal team logics with dual-like or bilateral negations, including Aloni's bilateral state-based modal logic, Hawke and Steinert-Threlkeld's semantic expressivist logic for epistemic modals, and the dual-negation variant of propositional dependence logic.
\end{abstract}



\end{frontmatter}

\input{sections/introduction}
\input{sections/preliminaries}

\input{sections/inconsistencies}

\input{sections/bsml_result}
\input{sections/propositional}
\input{sections/conclusion}

\input{sections/bibliography}

 \begin{acks}

The research for this paper was supported in part by grants 336283 and 368671 of the Academy of Finland, grant 101020762 of the European Research Council (ERC), and Research Funds of the University of Helsinki. %Part of the research was conducted while the author was affiliated with the Institute for Logic, Language and Computation, University of Amsterdam.

A preliminary version of this paper appeared in my PhD dissertation \cite{anttila2025}; I would like to thank my supervisors, Maria Aloni and Fan Yang, for extensive feedback. I am also grateful to Ivano Ciardelli, Marco Degano, Tomasz Klochowicz, Søren Knudstorp, Juha Kontinen, and Woxuan Zhou for helpful discussions on the contents of this paper.

%The author received support from the Academy of Finland, decision numbers 336283 and 368671, and from the European Research Council (ERC), grant 101020762.
 
 % The research for this paper was supported in part by 

 % of the second and third author was supported by grant 336283 of the Academy of Finland and Research Funds of the University of Helsinki.
 \end{acks}
%removed for blind review

\end{document}

%% file: sections/introduction.tex
\section{Introduction} \label{neg:section:introduction}

\subsection{The Dual Negation, Burgess' Theorem, and Bicompleteness} \label{neg:section:dual_negation}

Henkin quantifier logic (H) extends classical first-order logic (FO) with \emph{Henkin quantifiers} \cite{henkin} (also known as \emph{branching} or \emph{partially ordered quantifiers}) such as the following:
\begin{equation} \label{neg:equation:Henkin}
    \begin{pmatrix}
\forall x & \exists y \\
\forall u & \exists v 
\end{pmatrix} \phi(x,y,u,v).
\end{equation}
The intuitive meaning of (\ref{neg:equation:Henkin}) is that the value of $y$ depends only on that of $x$, and the value of $v$ only on that of $u$---(\ref{neg:equation:Henkin}) is equivalent to the existential second-order sentence $\exists f \exists g\forall x \forall u\phi(x,f(x),u,g(u))$ (cf. the FO-sentence $\forall x\exists y \forall u \exists v \phi(x,y,u,v)$). Hintikka and Sandu's \emph{independence-friendly logic} (IF) \cite{hintikka1989,hintikka1996}, in a similar manner, extends FO with \emph{slashed quantifiers} such as $(\exists y/\{x\})$, with the intended meaning of $(\exists y/\{x\})$ being that the value of $y$ must be chosen independently of the value of $x$ (an IF-sentence equivalent to (\ref{neg:equation:Henkin}): $\forall x\exists y \forall u (\exists v/\{x\} )\phi(x,y,u,v)$). Finally, \emph{dependence logic} (D) \cite{vaananen2007}, Väänänen's refinement of IF, replaces the slashed quantifiers with \emph{dependence atoms}: the atom $\dep{x_1,\ldots,x_n}{y}$ asserts that the value of $y$ is functionally determined by those of $x_1,\ldots ,x_n$ (a D-sentence equivalent to (\ref{neg:equation:Henkin}): $\forall x \exists y \forall u\exists v (\phi(x,y,u,v) \land \dep{u}{v})$).

%\emph{Team semantics}---originally introduced by Hodges \cite{hodges1997,hodges1997b} to provide a compositional semantics for Hintikka and Sandu's \cite{hintikka1989,hintikka1996} \emph{independence-friendly logic} (IF logic) and developed further by Väänänen \cite{vaananen2007} in his work on \emph{dependence logic}---is a generalization of the standard semantics of first-order logic in which formulas are interpreted with respect to \emph{sets} of variable assignments called \emph{teams} rather than single assignments. The shift to teams allows for a perspicuous way of representing quantifier dependencies 

%one to express relationships between assignment values not expressible in first-order logic (FO)---for instance, dependence logic extends first-order logic with \emph{dependence atoms}: the dependence atom $=\hspace{-0.1cm}(x,y)$ is true in a team if the values of $y$ in the team are functionally dependent on the values of $x$ in the team.

The prenex fragment of H (denoted $\mathrm{H_p}$; the formulas of $\mathrm{H_p}$ are of the prenex form Q$\phi$, where the prefix Q is a set of quantifiers with a partial order and the matrix $\phi$ is a quantifier-free FO-formula) is expressively equivalent to existential second-order logic $\Sigma^1_1$ \cite{walkoe,enderton}, as are the full logics IF \cite{hintikka1996} and D \cite{vaananen2007}.
%\footnote{To be more precise, $\Sigma^1_1$, IF, D, the prenex fragment of H are expressively equivalent in the sense that the classes of models definable using sentences of one of these logics are the same as the classes definable in any other of the logics. There is another sense in which IF and D are strictly less expressive than $\Sigma^1_1$: only 
%only a specific subclass of $\Sigma^1_1$-definable 
%non-empty \emph{downward} \cite{kontinen2009}.}
These logics are therefore not closed under classical negation $\bneg$ (where $\bneg\phi$ is true just in case $\phi$ is false).\footnote{This is because there are $\Sigma^1_1$-sentences $\phi$ such that $\bneg\phi$ is not equivalent to any $\Sigma^1_1$-sentence. This, in turn, follows, for instance, from the compactness of $\Sigma^1_1$ together with the fact that one can define a $\Sigma^1_1$-sentence which is true in a model just in case the model is infinite.} The original semantics for IF were game-theoretical, with the meaning of $(\exists y/\{x\})$ captured with a game rule to the effect that the player choosing the value of $y$ must do so without knowing the value of $x$---the game is one of \emph{imperfect information} (whereas in the game for FO the players have perfect information). This naturally led to the adoption (in both IF and D) of what is known as the \emph{dual} or \emph{game-theoretical} negation $\lnot$. The rule associated with this negation in the semantic game is simply the standard negation rule (also used in the game-theoretical semantics for FO) whereby the players switch their verifier/falsifier roles\footnote{\label{neg:footnote:dual_neg}The truth conditions of the dual negation can equivalently be obtained by stipulating that negated atomic formulas have their standard truth conditions, and that for complex formulas, the following dual equivalences must hold: $\phi\equiv \lnot \lnot \phi$; $\lnot(\phi\land \psi)\equiv \lnot \phi \vee \lnot \psi$; $\lnot(\phi\vee \psi)\equiv \lnot \phi \land \lnot \psi$; $\lnot \forall x\phi\equiv \exists x\lnot \phi$; $\lnot \exists x\phi\equiv \forall x\lnot \phi$; $\lnot (\exists y/\{x\})\phi\equiv (\forall y/\{x\} ) \lnot \phi$ and $\lnot (\forall y/\{x\})\phi\equiv (\exists y/\{x\})\lnot \phi$ (in IF); and $\lnot \dep{x_1,\ldots,x_n}{y}\equiv \bot$ (in D). (Note that there are multiple distinct versions of IF in the literature. I use the version in \cite{mann}.)}, and indeed, IF and D are conservative extensions of FO, with the dual negation extending the classical negation of FO\footnote{\label{neg:footnote:classical_negation}That is, for $\alpha$ an FO-formula, the dual-negated $\lnot\alpha $ is equivalent to $\alpha$ negated with the classical negation of FO. Note that this does not imply that IF/D are closed under (the full second-order) classical negation $\sim$.}. % in the FO-fragments of these logics, the dual negation is equivalent to the classical negation of FO.%---hence the name `dual negation'.

In the context of the full logics, however, the dual negation exhibits non-classical behaviour. Most germane to our purposes is the fact that \emph{preservation of equivalence under replacement} fails in negated contexts---for instance, $\phi \equiv \psi$ need not imply $\lnot \phi \equiv \lnot\psi$. (A concrete example: in dependence logic, we have $\lnot \dep{x}{y}\equiv \bot$, but $\lnot \lnot \dep{x}{y}\equiv \dep{x}{y} \nequiv \bot \equiv \lnot \lnot \bot$.) Put another way, each sentence $\phi$ defines a class of models $\left\Vert\phi\right\Vert=\{\MMM\;|\;\MMM\vDash\phi\}$ which we may think of as the meaning of $\phi$, and whereas negating a sentence in FO corresponds to the semantic operation of complementation on these classes of models, the dual negation fails to correspond to any operation on such classes: the class of models $\left\Vert\phi\right\Vert$ of $\phi$ does not determine the class of models $\left\Vert\lnot \phi\right\Vert$ of its dual negation. This does not mean that $\left\Vert\lnot \phi\right\Vert$ is fully unconstrained: one can show that $\phi$ and $\lnot \phi$ must be incompatible in that they share no models ($\left\Vert\phi\right\Vert\cap \left\Vert\lnot \phi\right\Vert=\emptyset$), with the effect that despite some standard inferential principles involving negation (such as negation introduction $\Gamma,\phi\vDash \bot \implies \Gamma \vDash \lnot \phi$) being invalidated due to the lack of determination, others (such as \emph{ex falso} $\phi,\lnot \phi\vDash \psi$) remain valid.

Burgess \cite{burgess2003}, however, showed that this failure of determination is extreme in the sense that incompatibility as formulated above is the \emph{only} constraint that $\left\Vert\phi\right\Vert$ places on $\left\Vert\lnot \phi\right\Vert$. That is, he showed (in the context of $\mathrm{H_p}$\footnote{Burgess proved his result for sentence of $\mathrm{H_p}$ and their \emph{contraries}, where the contrary of $Q\phi$ is obtained by swapping all quantifiers in $Q$ with their duals ($\forall$ with $\exists$ and vice versa), and (classically) negating the FO-matrix $\phi$. This is the same notion as the dual negation in that the dual negation of an IF/D-sentence $\phi$ can be computed in essentially the same manner---see note \hyperref[neg:footnote:dual_neg]{2}. Note that it is \emph{not} the case that given a sentence $\phi$ of $\mathrm{H_p}$ and an equivalent sentence $\psi$ of IF/D, the contrary of $\phi$ and the dual negation of $\psi$ are also equivalent---this is due to failure of replacement of equivalents under dual negation/contraries. It is the case, however, that given a sentence $\phi$ of $\mathrm{H_p}$ and an equivalent sentence $\psi$ of IF/D, there is some sentence $\chi$ of IF/D that is equivalent to both $\phi$ and $\psi$ and whose dual negation is equivalent to the contrary of $\phi$---this follows by the Burgess theorems for $\mathrm{H_p}$/IF/D.}) that for any sentences $\phi$ and $\psi$, if $\phi$ and $\psi$ share no models ($\left\Vert\phi\right\Vert\cap \left\Vert\psi\right\Vert=\emptyset$), then there is a sentence $\theta$ such that $\phi\equiv \theta$ and $\psi \equiv \lnot \theta$ ($\left\Vert\phi\right\Vert=\left\Vert\theta\right\Vert$ and $\left\Vert\psi\right\Vert=\left\Vert\lnot \theta\right\Vert$). So knowing only $\left\Vert\phi\right\Vert$ (without knowing the sentence $\phi$) tells us nothing whatsoever about $\left\Vert\lnot \phi\right\Vert$, save for the fact that the two classes are disjoint and the fact that $\left\Vert\lnot \phi\right\Vert$ is expressible in $\mathrm{H_p}$: for a given $\left\Vert\phi\right\Vert$, take any class $\mathcal{X}$ expressible in $\mathrm{H_p}$ and disjoint with $\left\Vert\phi\right\Vert$---by Burgess' theorem, there is then a $\theta$ with $\left\Vert\phi\right\Vert=\left\Vert\theta\right\Vert$ and $\mathcal{X}=\left\Vert\lnot \theta\right\Vert$, so for all we know, $\left\Vert\lnot \phi\right\Vert$ might be $\mathcal{X}$. See Figure \ref{neg:fig:burgess}. Dechesne \cite{dechesne} later reformulated and proved the result for IF. Kontinen and Väänänen \cite{kontinen2011} (working in D and IF) generalized the theorem to arbitrary formulas (that may contain free variables). Mann (in \cite{mann}) proved an analogue of the result for sentences for the \emph{perfect-recall} fragment of IF.\footnote{See also \cite{bellier2023}, which introduces an extension of IF featuring a negation with a natural game-theoretic interpretation which does not exhibit failure of determination.} %\todoa{Add comment about result only holding for models of cardinality $\geq 2$?}

{ 
  \begin{figure}[t]  
\centering
  \subfigure[FO]{\hspace{0.1cm}
\begin{tikzpicture}[>=latex,scale=.85]
\label{neg:fig:burgess_a}

 % Frame
 \draw (-2,1.5) rectangle (2, -1.5);
  \draw (0,1.5) -- (0, -1.5);

 % Indices
\draw (-1,0) node[scale=0.9] (yy) {$\left\Vert\phi\right\Vert$};
\draw (1,0) node[scale=0.9] (yy) {$\left\Vert\lnot \phi\right\Vert$};

\end{tikzpicture}\hspace{0.1cm} }
\hspace{0.2cm}
 \subfigure[\footnotesize $\mathrm{H_p}$]{ 
 \hspace{0.1cm}
 
 \begin{tikzpicture}[>=latex,scale=.85]
 \label{neg:fig:non-zero}

 % Possibilities
\draw (-1,0.5) node[scale=0.9] (yy) {$\left\Vert\phi\right\Vert$};
\draw (-1,0) node[scale=0.9] (yy) {$\left\Vert\theta_1\right\Vert$};
\draw (-1,-0.5) node[scale=0.9] (yy) {$\left\Vert\theta_2\right\Vert$};
\draw (-1,0.5) node[scale=0.9] (yy) {$\left\Vert\phi\right\Vert$};
\draw (0.75,0.5) node[scale=0.9] (yy) {$\left\Vert\lnot\theta_1\right\Vert$};
\draw (0.75,-0.5) node[scale=0.9] (yy) {$\left\Vert\lnot\theta_2\right\Vert$};
 \draw[opaque,rounded corners] (-1.5,1) rectangle (-0.5, -1);
 \draw[opaque,rounded corners] (0,1) rectangle (1.5, 0);
 \draw[opaque,rounded corners] (0.2,0.2) rectangle (1.3, -1);

 % Frame
  \draw (-2,1.5) rectangle (2, -1.5);

\end{tikzpicture}
}
\caption{Failure of determination: in FO, $\left\Vert\lnot\phi\right\Vert$ is the complement of $\left\Vert\phi\right\Vert$. In $\mathrm{H_p}$, given only $\left\Vert\phi\right\Vert$, $\left\Vert\lnot \phi\right\Vert$ can be any class disjoint with $\left\Vert\phi\right\Vert$ that is expressible in $\mathrm{H_p}$.} \label{neg:fig:burgess}
\end{figure}
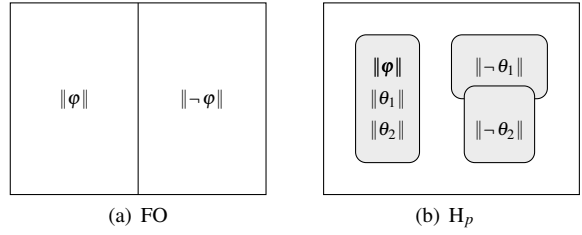  
}

Burgess observed that his theorem can be seen as a result concerning the expressive power of $\mathrm{H_p}$:

\begin{quote}
The Enderton-Walkoe theorem says that for any PC [pseudo-elementary (i.e., $\Sigma^1_1$-definable) class of models], call it $K$, there is a Henkin sentence $\theta$ such that $K=\left\Vert \theta \right\Vert$. The corollary just proved allows this theorem to be strengthened to say that for any two disjoint PCs, call them $K_0$ and $K_1$, there is a Henkin sentence $\theta$ such that $K_0=\left\Vert \theta \right\Vert$ and $K_1=\left\Vert \lnot\theta \right\Vert$. \cite[Notation amended.]{burgess2003}
\end{quote}

\noindent
Dechesne \cite{dechesne} noticed that, together with its converse (any $\mathrm{H_p}$/IF/D-sentence and its dual negation share no models), Burgess' theorem further yields an expressive completeness result for pairs: the class of pairs $(\left \Vert \phi\right \Vert,\left \Vert \lnot\phi\right \Vert )$ of classes of models definable by a $\mathrm{H_p}$/IF/D-sentence $\phi$ and its dual negation is precisely the class of pairs $(\mathcal{X},\mathcal{Y})$ such that $\mathcal{X}$ and $\mathcal{Y}$ are disjoint and both are $\Sigma^1_1$-definable. We will call this type of expressive completeness result with respect to properties/classes of pairs a \emph{bicompleteness} theorem: on the level of sentences, $\mathrm{H_p}$/IF/D is bicomplete for disjoint pairs of $\Sigma^1_1$-definable classes. We also say that, relative to its expressive power, each of these logics is bicomplete for disjoint pairs.% on the level of sentences.%; by Kontinen and Väänänen's result each of IF and D is also bicomplete for 

\subsection{Burgess and Bicompleteness Theorems for Modal and Propositional Team Logics} \label{neg:section:team_logics}

In this paper, we prove analogues of Burgess' theorem as well as bicompleteness theorems for a number of modal and propositional \emph{team logics}.

Team logics are logics such as dependence logic that are primarily intended to be interpreted using \emph{team semantics}. In team semantics (introduced in an early form by Hodges \cite{hodges1997,hodges1997b} to provide a compositional semantics for IF, and later developed into its contemporary form by Väänänen and Hodges \cite{vaananen2007,vaananenhodges}; also independently developed in the guise of \emph{inquisitive semantics} chiefly by Ciardelli, Groenendijk, and Roelofsen \cite{CiardelliRoelofsen2011,inqsembook,ciardellibook}), formulas are evaluated with respect to \emph{sets} of evaluation points called \emph{teams} rather than single points, as in standard Tarskian semantics. For instance, in first-order team semantics, teams are sets of variable assignments; in modal team semantics \cite{vaananen2008}, teams are sets of possible worlds, etc.\footnote{\label{neg:footnote:team_logic}The term `team logic' is from Väänänen \cite{vaananen2007}, who originally used it to refer to the extension of D with the \emph{Boolean} or \emph{contradictory negation} $\bneg$ (where $\bneg \phi$ is true in a team iff $\phi$ is false), and it is now frequently used to refer to any logic employing team semantics that incorporates $\bneg$. The term is used both in this sense and in the more general sense (whereby `team logics' are logics intended to be interpreted using team semantics) in the literature. 

We noted earlier that D and IF are not closed under $\bneg$, and referred to $\bneg$ as the `classical negation'. This name is appropriate when discussing the relationship between D, IF, and $\Sigma^1_1$ because on the level on which D and IF correspond to $\Sigma^1_1$ (essentially the level of teams), $\bneg$ corresponds to the standard, classical negation of second-order logic. In the context of team semantics, the name `classical negation' would potentially be confusing because %, while on the level of teams $\bneg $ has the classical negation truth conditions, 
there is another level (essentially the level of assignments or worlds---the elements of teams) on which $\bneg$ behaves nonclassically (on this level it is the dual negation rather than $\bneg$ which corresponds, in the FO-fragments of D and IF, to the classical negation of FO---see note \hyperref[neg:footnote:classical_negation]{3} and Fact \ref{neg:fact:ML_team_classical_correspondence}.)}

%Dependence logic and its descendants are usually interpreted using \emph{team semantics} (first introduced by Hodges \cite{hodges1997,hodges1997b} to provide a compositional semantics for IF)---we will accordingly refer to these logics as `team logics'\footnote{\label{neg:footnote:team_logic}Väänänen \cite{vaananen2007} originally used `team logic' to refer to the extension D with $\bneg$ (where $\bneg \phi$ is true in a team iff $\phi$ is false), and it is now frequently used to refer to any logic employing team semantics that incorporates $\bneg$. The term is used both in this sense and in the more general sense in the literature. 

%In the team semantics literature, $\bneg$ is usually called the \emph{Boolean} or \emph{contradictory negation}. The term `classical negation' is potentially confusing because, while on the level of teams $\bneg $ has the classical negation truth conditions, there is another level on which it behaves nonclassically, and the negation in the classical fragment of a team logic is not $\bneg$. See Fact \ref{neg:fact:ML_team_classical_correspondence} and Section \ref{neg:section:classical_logics}.
%}. In team semantics, formulas are evaluated with respect to \emph{sets} of evaluation points called \emph{teams} rather than single points, as in standard Tarskian semantics (so, for instance, in first-order team semantics, teams are sets of variable assignments; in modal team semantics \cite{vaananen2008}, teams are sets of possible worlds, etc.).

Contemporary team logics in the lineage of dependence logic typically do not incorporate the dual negation. However, interest in this negation has recently been reinvigorated due to team logics in the philosophical logic and formal semantics literature which employ \emph{bilateral} notions of negation similar to the dual negation. We will consider two such logics in this paper: Aloni's \emph{bilateral state-based modal logic} ($\BSML$) \cite{aloni2022}, and Hawke and Steinert-Threlkeld's semantic expressivist logic ($\mathrm{HS}$) for epistemic modals \cite{hawke}. These are both propositional modal team logics developed to account for linguistic phenomena such as \emph{free choice inferences} and \emph{epistemic contradictions}. The semantics of both of these logics is formulated bilaterally---that is, using both a positive primitive semantic relation $\vDash$ (interpreted as assertibility), as well as a negative primitive relation $\Dashv$ (interpreted as rejectability). The bilateral semantics is used to define a \emph{bilateral negation} $\lnot$: for a team $s$, $s\vDash \lnot \phi$ iff $s\Dashv \phi$, and $s\Dashv \lnot \phi$ iff $s\vDash \phi$. The bilateral negations of these logics are similar to the dual negation in that they exhibit failure of determination; dual equivalences similar to those of IF/D hold in $\BSML$; and, as we will show, an analogue of Burgess' theorem holds in both $\BSML$ and $\mathrm{HS}$.

In the first part of this paper, we show the analogue and associated bicompleteness theorem for $\BSML$ (as well as for $\BSMLI$, the extension of $\BSML$ with the \emph{global} or \emph{inquisitive disjunction} $\vvee$). Of the aforementioned theorems, the most directly comparable to our analogues is Kontinen and Väänänen's result for D/IF. The main substantial difference between our analogues and this result is that the notion of incompatibility must be adjusted, or, in other words, relative to expressive power, $\BSML$/$\BSMLI$ and D/IF are bicomplete with respect to different properties of pairs:  whereas Kontinen and Väänänen showed that, relative to expressive power, D/IF are bicomplete for what we call \emph{$\bot$-incompatible} pairs, $\BSML$/$\BSMLI$ are bicomplete for \emph{ground-incompatible} pairs. These incompatibility notions/pair properties interact with \emph{team-semantic closure properties} in interesting ways. Formulas of D/IF are \emph{downward closed}---the truth of a formula in a team implies its truth in all subteams---whereas formulas of $\BSML$/$\BSMLI$ need not be. Ground-incompatibility is in general stronger than $\bot$-incompatibility, but the two notions are equivalent in a downward-closed setting (this has the effect that, relative to expressive power, D/IF are in fact bicomplete for both $\bot$-incompatible pairs as well as for ground-incompatible pairs; $\BSML$/$\BSMLI$ are only bicomplete for the latter.) The incompatibility notions are also conceptually suggestive; we will briefly comment on possible intuitive interpretations of the notions as well as their connection to \emph{epistemic contradictions}.

Whereas the logics in the first part of the paper employ modal team semantics---team semantics on Kripke models---in the second part of the paper, we consider logics with dual or bilateral negations interpreted using propositional team semantics. The simpler propositional setting allows us to readily prove multiple Burgess/bicompleteness theorems with respect to many different properties of pairs, giving us a better handle on the notion of bicompleteness as well as on the diverse behavior this type of negation can be induced to exhibit. We prove bicompleteness theorems for Hawke and Steinert-Threlkeld's logic $\HS$ \cite{hawke} (which we treat here as a propositional logic despite its featuring a modality since it is interpreted on propositional teams); for the propositional fragments of $\BSML$, and $\BSMLI$; %(both introduced in \cite{yang2017}); 
and for \emph{propositional dependence logic} \cite{yangvaananen2016} with the dual negation, the propositional version of D. %These theorems require the introduction of yet more pair properties/incompatibility notions. 
%We show that ground-incompatibility yields a bicompleteness theorem for the propositional fragment of $\BSMLI$, whereas the propositional fragment of $\BSML$ (both introduced in \cite{yang2017}) requires an even stronger notion. We show that both this stronger notion as well as a further strengthened version of it yield a bicompleteness theorem for \emph{propositional dependence logic} \cite{yangvaananen2016} with a dual-like negation, the propositional version of D.\
Whereas the proofs in the first part of the paper are analogous to Burgess', the results in the second part are, for the most part, proved by modifying a standard unilateral expressive completeness theorem proof for the logic in question. We also list some trivial bicompleteness theorems for propositional logics which do not exhibit failure of determination (for instance, classical propositional logic interpreted on teams), and provide an example of maximal failure of determination by introducing a logic which is bicomplete for all propositional pair properties.

%for the modal logics (as well as that for the propositional fragment of $\BSMLI$) are analogous to Burgess' (these results are essentially corollaries of interpolation for classical modal/propositional logic; Burgess', similarly, is a corollary of interpolation for FO), we employ different methods to prove the other propositional theorems: the $\HS$-result follows easily from the basic properties of the negation in $\HS$, while the results for the propositional fragment of $\BSML$ and propositional dependence logic are shown to follow from the expressive completeness theorems for these logics. 

%While we distinguish the different pair properties/incompatibility notions for technical reasons, many of the notions we consider are also conceptually suggestive. In the third part of the paper, we sketch some possible intuitive interpretations of some of the notions.

We prove the $\BSML$$\BSMLI$-theorems in Section \ref{neg:section:BSML} and the propositional theorems in Section \ref{neg:section:propositional}. 
%We consider possible interpretations of the pair properties in Section \ref{neg:section:interpretations}. 
We conclude, in Section \ref{neg:section:conclusion}, with some discussion of the results. In the remainder of this introduction, I comment briefly on how Burgess' theorem has been interpreted.

A preliminary version of this paper appeared in my PhD dissertation \cite{anttila2025}.

\subsection{Remarks on Burgess' Theorem} \label{neg:section:remarks_on_Burgess_theorem} Burgess himself intended his theorem to serve, in part, as a point against IF and Hintikka's philosophical ambitions. He writes:
\begin{quote}
In recent years Hintikka and co-workers have revived a variant version of the
logic of Henkin sentences under the label “independence-friendly” logic, have restated many theorems about existential second-order sentences for this “new” logic,
and have made very large claims about the philosophical importance of the theorems
thus restated. In discussion, pro and con, of such philosophical claims it has not been sufficiently emphasized that contrariety [dual negation], the only kind of “negation” available, fails to correspond to any operation on classes of models. For this reason it seemed worthwhile to set down, in the form of the corollary above, a clear statement of just how total the failure is. \cite{burgess2003}
\end{quote}
Accordingly, the result is occasionally cited as attesting that the behaviour of the negation is anomalous or problematic (as in \cite{Humberstone2018}). A common gloss has it that Burgess establishes or shows that the dual negation is ``not a semantic operation''.

I do not intend to weigh in on this debate, or to argue for or against the dual negation here. What I do hope to do is to point out why the common characterization of the theorem is potentially misleading in multiple ways.\footnote{It is also worth pointing out, in connection with the passage from Burgess concerning Hintikka and IF, that Hintikka did also consider an extension of IF with the Boolean/contradictory negation $\bneg$ (\emph{extended independence-friendly logic}), and that he ultimately viewed each negation as indispensable \cite{hintikka1996,hintikkasandu1997,hintikka2002}. Consider, for instance, the following passage \cite[p. 154]{hintikka1996}:\\

\indent\begin{minipage}{.8\textwidth}
    %\begin{quote}
$\ldots$\emph{in any sufficiently rich language, there will be two different notions of negation present.} Or if you prefer a different formulation, our ordinary concept of negation is intrinsically ambiguous. The reason is that one of the central things we certainly want to express in our language is the contradictory negation. But $\ldots$ a contradictory negation is not self-sufficient. In order to have actual rules for dealing with negation, one must also have the dual negation present, however implicitly.
%\end{quote}
\end{minipage}}

My first qualm with the gloss is that it is not Burgess' theorem that establishes the ``non-semantic'' nature of the dual negation. The failure of the negation to correspond to any operation on classes of models is an easily observable, simple fact (recall our dependence logic example: $\lnot \dep{x}{y}\equiv \bot$, but $\lnot \lnot \dep{x}{y} \nequiv \lnot \lnot \bot$). The theorem is a further fact, characterized by Burgess in the passage above as concerning the degree of this failure: the (implied) problem with the negation is its failure to correspond to any operation on classes of models; the theorem demonstrates the extreme degree of this failure. The failure itself is independent of the theorem. The results in this paper help to reinforce this point: we prove variants of Burgess' theorem with respect to multiple distinct pair properties, which we may think of as points on a scale for measuring the degree of failure of correspondence, in that a Burgess theorem with respect to a weaker property corresponds to a higher degree of failure (see Section \ref{neg:section:conclusion}). We have, then, failure of determination in multiple different logics, regardless of whether the direct analogue of Burgess' original theorem (employing $\bot$-incompatibility) holds for that logic (and also regardless of whether some variant of the theorem can be established), with the Burgess theorems providing insight into the degree and nature of the failure in each case.

Second, as essentially already pointed out by Hodges \cite{hodges1997}, and later by many others \cite{vaananen2007,dechesne,kontinen2011,gradel2012}, %(and as is already implicit in our discussion of bicompleteness)
 there is a sense in which the dual negation is a semantic operation, which, while perhaps trivial, should not be ignored. Indeed, Hodges' \cite{hodges1997} goal in developing team semantics was to demonstrate that logics of imperfect information such as IF could be given a semantics in which all connectives are semantic operations in the sense that replacing a subformula occurrence with one with the same meaning does not change the meaning of a formula.
%---a compositional semantics with the \emph{full abstraction} property\footnote{\label{neg:footnote:abstraction}In Hodges' \cite{hodges1997} formulation, the semantics of a language have the full abstraction property if\\
%\indent\begin{minipage}{.8\textwidth}
%Two formulas $\phi(\vec{x})$ and $\psi(\vec{x})$, say of signature $\Sigma$, have the same meaning if and only if for every signature $\Sigma'$ containing $\Sigma$, every sentence $\chi$ of signature $\Sigma'$, and every structure $A$ of signature $\Sigma'$, the truth or otherwise of $\chi$ in $A$ is not affected if we replace an occurrence of $\phi$ in $\chi$ by an occurrence of $\psi$.
%\end{minipage}}, as Hodges puts it.) 
Instead of taking the meaning of $\phi$ to be $\left\Vert \phi\right\Vert$, let it be represented by the pair $(\left\Vert \phi\right\Vert, \left\Vert \lnot\phi\right\Vert)$. Then, given that double negation elimination is valid in $\mathrm{H}_p$/IF/D, the meaning of $\lnot \phi$ is determined by that of $\phi$ since to obtain the former from the latter, one need only flip the elements of the pair: $(\left\Vert \lnot \phi\right\Vert, \left\Vert \phi\right\Vert)$.

One could of course perform a similar trick with any other operator or connective no matter its semantics, but there are potentially good reasons for insisting that, applied to the negation, this move is not entirely \emph{ad hoc}. For instance, on the view known both as \emph{rejectivism} and as \emph{bilateralism} (see, for instance, \cite{price1983,price1990,smiley1996,rumfitt2000}), the speech act of rejection is not reducible to that of assertion---these notions should be treated as being on par, with neither being conceptually reducible to the other. If one endorses this view, it is natural to conceive of the pair $(\left\Vert \phi\right\Vert^+, \left\Vert \phi\right\Vert^-)$ as the full meaning of $\phi$, where $\left\Vert \phi\right\Vert^+$ denotes the models in which in $\phi$ is assertible, and $\left\Vert \phi\right\Vert^-$ those in which it is rejectable. And if, as in $\BSML$, one additionally associates rejectability with a negation operator (cf. varieties of rejectivism with a specialized rejection operator such as \cite{smiley1996,rumfitt2000}), one may equate $(\left\Vert \phi\right\Vert^+, \left\Vert \phi\right\Vert^-)$ with $(\left\Vert \phi\right\Vert, \left\Vert \lnot\phi\right\Vert)$.\footnote{\label{neg:footnote:BSML_bilateralism}The bilateral semantics of $\BSML$ and Hawke and Steinert-Threlkeld's logic are not motivated by rejectivism \emph{per se}; the primary intended function of bilateralism in these logics is, rather, to allow these logics to correctly capture and predict linguistic data. Aloni \cite{aloni2022} makes a point in defence of failure of replacement by drawing on these empirical considerations:\\

\indent\begin{minipage}{.8\textwidth}
As we will see, however, this non-classical behavior is precisely what we need to explain the effects of pragmatic enrichment in negative contexts. $[\ldots]$ when replacement under $\lnot$ holds, the rejection conditions are derivable from the support conditions, and so without failure of replacement bilateralism would not give different predictions from unilateral systems and so it would be empirically unjustified.
\end{minipage}} Whether it is possible to formulate a tenable form of rejectivism featuring the dual negation in this way is beyond the scope of this paper; at any rate, given the clear similarity between Hodges' notion of meaning and that of the rejectivists, it seemed worthwhile to make this connection explicit (rejectivism has not thus far, to my knowledge, been discussed in the literature on the dual negation). It should also be noted that, as alluded to in endnote \hyperref[neg:footnote:BSML_bilateralism]{8}, bilateral notions of meaning of this kind have been proposed for reasons other than rejectivism. %\todoa{Should I put references here? Is it correct to say that bilateralism in truthmaker semantics is not motivated by rejectivism? I assumed it was correct to say this for BSML; is that the case?}

%% file: sections/preliminaries.tex
\section{Burgess Theorems for $\BSMLBF$ and $\BSMLIBF$} \label{neg:section:BSML}
In this section, we list the required preliminaries concerning the logics $\BSML$ and $\BSMLI$ and team semantics (Section \ref{neg:section:preliminaries}); reformulate Burgess' theorem as an expressive completeness theorem for pairs satisfying a specific notion of incompatibility (Section \ref{neg:section:inconsistencies}); distinguish between multiple notions of incompatibility and examine how they are related (Section \ref{neg:section:inconsistencies}); and show Burgess theorems for $\BSML$ and $\BSMLI$ using some of these incompatibility notions (Section \ref{neg:section:BSML_theorems}).

\subsection{Preliminaries} \label{neg:section:preliminaries}
We first define the syntax and semantics of $\BSML$ and $\BSMLI$ (from \cite{aloni2022,aloni2023}); then recall the definitions of standard team-semantic closure properties and list some results relating these properties to our logics (from \cite{aloni2023,aknudstorp2024}); and, finally, recall basic notions and results from modal logic (see, e.g., \cite{blackburn2001,goranko2007}), together with team-based analogues (from \cite{hella2014,aloni2023}).

$\BSML$ is an extension of classical modal logic with the \emph{nonemptiness atom} $\NE$ (introduced in \cite{vaananen2014,yang2017}), which is true in a team just in case the team is nonempty. $\BSMLI$ is $\BSML$ extended with the \emph{global} or \emph{inquisitive disjunction} $\intd$, used in inquisitive logic \cite{CiardelliRoelofsen2011,inqsembook,ciardellibook} (see Section \ref{neg:section:classical_logics}) to model the meanings of questions.

\begin{definition}[Syntax of $\MLBF$, $\BSMLBF$, and $\BSMLIBF$] \label{neg:def:syntax}
Fix a (countably infinite) set $\mathsf{Prop}$ of propositional variables. The set of formulas of \emph{bilateral state-based modal logic} $\BSML$ is generated by:
        \begin{align*}
            \phi ::= p \sepp \bot \sepp \bnot \phi \sepp ( \phi \land \phi )  \sepp ( \phi \dis \phi ) \sepp \Di \phi   \sepp \NE
        \end{align*}
        where $p\in \mathsf{Prop}$.   
        \emph{Classical modal logic} $\ML$ is the $\NE$-free fragment of $\BSML$. $\BSMLI$ is the extension of $\BSML$ with the binary connective $\vvee$.
\end{definition}

We use the first Greek letters $\alpha,\beta$ to refer exclusively to formulas of \ML (also called {\em classical formulas}).
We write $\mathsf{P}(\phi)$ for the set of propositional variables in $\phi$. %We write $\phi (\psi/p)$ for the result of replacing all occurrences of $p$ in $\phi$ by $\psi$.

A Kripke model $M = (W,R,V)$ over $\mathsf{X}\subseteq \mathsf{Prop}$ is defined as usual, where in particular, $V: \mathsf{X} \to \wp(W)$. We call a subset $s\subseteq W$ of $W$ a (modal) \emph{team} on $M$. For any world $w$ in $M$, define, as usual, $R[w]:=\{v \in W\sepp wRv\}$. Similarly, for any team $s$ on $M$, define $R[s]:=\bigcup_{w\in s}R[w]$.

\begin{definition}[Semantics of $\MLBF$, $\BSMLBF$ and $\BSMLIBF$] \label{neg:def:semantics}
    For a model $M=(W,R,V)$ over $\mathsf{X}$, a
    team $s$ on $M$, and formula $\phi$ with $\mathsf{P}(\phi)\subseteq \mathsf{X}$, the notions of $\phi$ being \emph{supported by}/\emph{anti-supported by} $s$ in $M$, written $M,s \vDash \phi$/$M,s\Dashv \phi$ (or simply $s \vDash \phi$/$s\Dashv \phi$), are defined recursively as follows:
    \begin{center}
    \begin{longtable}{p{2cm} p{1cm} p{8.5cm}}
        $M,s \vDash   p$  & $:\iff$&  for all $ w \in s: w \in V(p)$ \\
        $M,s \Dashv  p  $& $:\iff$& for all $ w \in s: w \notin V(p)$\\
        &&\\
        $M,s \vDash   \bot$  & $:\iff$&  $s=\emptyset$ \\
        $M,s \Dashv  \bot  $& & always the case \\
        &&\\
                $M,s \vDash    \NE$ &$ :\iff$&$ s \neq \emptyset$ \\
        $M,s \Dashv    \NE $&$ :\iff$&$ s =\emptyset$\\
        &&\\
        $M,s \vDash    \bnot \phi$ & $:\iff$ & $M,s\Dashv \phi$ \\
        $M,s \Dashv    \bnot \phi $&$ :\iff $&$ M,s \vDash \phi $\\
        &&\\
        $M,s \vDash    \phi \land \psi  $&$ :\iff $&$ M,s \vDash  \phi\text{ and } M,s \vDash \psi$\\
        $M,s \Dashv    \phi \land \psi  $& $:\iff$& there exist $t,u$ s.t. $ s=t \cup u\text{ and } M,t \Dashv  \phi\text{ and } M,u \Dashv  \psi$ \\
        &&\\
        $M,s \vDash    \phi \dis \psi  $& $:\iff$& there exist $t,u $ s.t. $ s=t\cup u\text{ and } M,t \vDash  \phi\text{ and } M,u \vDash  \psi$ \\
        $M,s \Dashv    \phi \dis \psi  $&$ :\iff $&$ M,s \Dashv  \phi\text{ and } M,s \Dashv  \psi$\\
        &&\\
                $M,s \vDash    \phi \intd \psi  $&$ :\iff$&$ M,s \vDash  \phi \text{ or } M,s \vDash  \psi$\\
        $M,s \Dashv    \phi \intd \psi  $&$ :\iff$&$  M,s \Dashv  \phi\text{ and } M,s \Dashv  \psi$\\
        &&\\
        $M,s \vDash    \Di \phi  $& $:\iff$& for all $ w \in s$ there exists $  t \subseteq R[w]$ s.t. $ t\neq \emptyset \text{ and }M,t \vDash \phi$ \\
        $M,s \Dashv    \Di \phi  $&$ :\iff$& for all $ w \in s : M,R[w]\Dashv \phi$
    \end{longtable}
    \end{center}
\end{definition}   
\setcounter{table}{0}
\vspace{-0.5cm}
We also refer to support by $s$ as \emph{truth in $s$} for convenience, although we caution that most interpretations of team logics rely on distinguishing the way in which a team satisfies a formula (often called `support') from the way in which an element of a team satisfies one (often called `truth'; see Fact \ref{neg:fact:ML_team_classical_correspondence}). We write $\phi\vDash\psi$ and say $\phi$ \emph{entails} $\psi$ if for all $M$ and all $s$ on $M$, $M,s\vDash \phi$ implies $M,s\vDash \psi$. If both $\phi\vDash\psi$ and $\psi\vDash\phi$, then we write $\phi\equiv\psi$ and say that $\phi$ and $\psi$ are \emph{equivalent}. If both $\phi \equiv \psi$ and $\lnot \phi \equiv \lnot \psi$, then $\phi$ and $\psi$ are said to be \emph{bi-equivalent} or \emph{strongly equivalent}, written $\phi \equiv^\pm \psi$.

%The box modality $\Bo$ is defined as the dual of the diamond: $\Box \phi:=\lnot \Di\lnot \phi$; the resulting support/antisupport clauses are:
%\begin{center}
%    \begin{tabular}{p{2cm} p{1cm} p{8.5cm}}
%        $M,s \vDash    \Bo \phi  $&$ \iff$&for all $ w \in s:  M,R[w]\vDash \phi$\\
 %       $M,s \Dashv    \Bo \phi  $& $\iff$&for all $ w \in s$ there exists $ t \subseteq R[w]$ s.t. $ t\neq \emptyset \text{ and }M,t \Dashv \phi$ \\
%        &&\\
%    \end{tabular}
%\end{center}

    We refer to the atom $\bot$ as the \emph{weak contradiction}. We also define the \emph{strong contradiction} $\Bot:=\bot \land \NE$, and $\top:=\lnot \bot$. The weak contradiction is true only in the empty team; the strong contradiction in no team; $\top$ in all teams. Note that \emph{ex falso} with respect to $\Bot$ holds for all formulas $\phi$: $\Bot\vDash \phi$; with respect to $\bot$ it holds only for formulas $\psi$ with the empty team property (see below): $\bot\vDash \psi$ if $M,\emptyset\vDash \psi$ for all $M$. We let $\bigdis \emptyset := \bot$; $\bigintd \emptyset := \Bot$; and $\Box \phi:=\lnot \Di\lnot \phi$. Observe that $\lnot (\phi \vee \psi)\equiv \lnot \phi \land \lnot \psi\equiv \lnot (\phi \vvee \psi)$; $\lnot \lnot \phi\equiv \phi$; and $\lnot \Diamond \phi \equiv \Box \lnot \phi$. It follows from these equivalences that each formula is equivalent to one in negation normal form.

%   We use these contradictions and tautologies to interpret the empty disjunctions and conjunction:

    \begin{definition}[Closure properties] \label{neg:def:closure_properties}
We say that a formula $\phi$ 
\begin{itemize}
    \item[--] is \emph{downward closed}, provided   $[M,s \vDash \phi \text{ and } t \subseteq s] \implies M,t \vDash \phi$;
    \item[--] is \emph{convex}, provided   $[M,s \vDash \phi ;$ $M,t \vDash \phi; \text{ and } t \subseteq u\subseteq  s] \implies M,u \vDash \phi$;
     \item[--] is \emph{union closed}, provided $[M,s \vDash \phi \text{ for all } s\in S\neq \emptyset] \implies M, \bigcup S  \vDash \phi$;
    \item[--] has the \emph{empty team property}, provided $M, \emptyset \vDash \phi \text{ for all }M$;
     \item[--] is \emph{flat}, provided $M, s\vDash\phi  \iff M,\{w\}\vDash \phi\text{ for all }w\in s$.
\end{itemize}
\end{definition}
  It is easy to check that a formula is flat if and only if it is downwards closed and union closed, and has the empty team property. One can show by induction that, in the context of the connectives we have introduced, all $\vvee$-free formulas (in particular, all formulas of \BSML) are union closed and convex, and that all $\NE$-free formulas are downward closed and have the empty team property (so that all formulas of \ML are flat). Formulas of $\BSML$ or $\BSMLI$ clearly need not be downward closed or have the empty team property (consider $p\land \NE$); formulas of $\BSMLI$ need not be union closed or convex (consider $(a\land ((p \land \NE) \vee \top)) \vvee (b\land ((q \land \NE) \vee \top))$. Another easy induction shows:

\begin{fact}[Conservativity of $\MLBF$ over standard semantics] \label{neg:fact:ML_team_classical_correspondence}
For any $\alpha\in \ML$:
\begin{align*}
    &M,s \vDash \alpha &&\iff &&M,\{w\}\vDash\alpha\text{ for all }w\in s&&\iff &&M,w \vDash \alpha \text{ for all }w \in s,
\end{align*}
\end{fact}
\noindent
where $\vDash$ on the right is the standard single-world truth relation for $\ML$. This implies that for $\alpha,\beta \in \ML$: $\alpha\vDash \beta$ in the team-semantics sense iff $\alpha \vDash \beta$ in the usual single-world-semantics sense; we may therefore use the symbol `$\vDash$' to refer interchangeably to either type of entailment when discussing $\ML$-formulas.

We write $\chi[\xi]$ to refer to a specific occurrence of the subformula $\xi$ in $\chi$, and $\chi[\psi/\xi]$ for the result of replacing this occurrence of $\xi$ in $\chi$ with $\phi$. Our Burgess theorems will build on the fact that preservation of equivalence under replacement fails in negative contexts:  if $\chi[p]$ is within the scope of an odd number of $\lnot$ in $\chi$, $\phi\equiv \psi$ need not imply $\chi[\phi/p]\equiv\chi[\psi/p]$. For instance, $ \bot\equiv \lnot \NE$ but $ \lnot \bot\equiv \top \notequiv \NE \equiv \lnot \lnot\NE$, and $\lnot (\phi \lor \psi)\equiv \lnot (\phi \intd \psi)$ but $\lnot\lnot (\phi \lor \psi) \equiv\phi \lor \psi \notequiv  \phi \intd \psi \equiv\lnot \lnot (\phi \intd \psi)$. Bi-equivalence, on the other hand, is preserved under replacement: if $\phi \equiv^\pm \psi$, then $ \chi[\phi/p] \equiv^\pm \chi[\psi/p]$. (This implies that if we take meanings to be represented by pairs $\left \Vert \phi\right\Vert^\pm=(\left\Vert\phi\right\Vert,\left\Vert\lnot\phi\right\Vert )$---see below for the definition of this notation---the semantics is compositional in Hodges' sense. See Section \ref{neg:section:remarks_on_Burgess_theorem}.) And, crucially, replacement does hold in non-negated and other positive contexts:
\begin{fact}[Replacement in positive contexts] \label{neg:fact:replacement}
        If $\chi[p]$ is in the scope of an even number of $\lnot$ in $\chi$, then $\phi\equiv \psi$ implies $\chi[\phi/p]\equiv\chi[\psi/p]$.
\end{fact}

%\todoa{Introduce dual negation, bilateral negation, negation is non-semantic (for singletons), negation is non-semantic for pairs}

We officially reserve the term \emph{dual negation} to refer to any negation operator that engenders failure of preservation of equivalence under replacement in negated contexts (as noted above, in $\BSML$ and $\BSMLI$ this failure occurs only in \emph{negative} negated contexts---contexts in the scope of odd number of negations). However, for simplicity, we also occasionally use this term to refer to any negation operator with bilateral semantics (such as the negation of $\ML$ as defined above).

%We now recall some standard notions and results for classical modal logic, and see how they can be extended to the team-based setting. For further details on the classical results, see, e.g. \cite{blackburn2001,goranko2007}; for more on team-based notions and $\BSML$, see \cite{hella2014,aloni2023}.

A \emph{pointed model} over $\mathsf{X}$ is a pair $(M,w)$ where $M$ is a model over $\mathsf{X}$ and $w\in W$. We define \emph{$k$-bisimilarity} (for $k\in \mathbb{N}$) between pointed models with respect to a set of propositional variables $\mathsf{X}$ in the usual way:
    \begin{itemize}
        \item[--] $M,w \bisim^{\mathsf{X}}_0 M,w :\iff $ for all $p \in \mathsf{X}$ we have $M,w \vDash p$ iff $M',w' \vDash p$.
        \item[--] $M,w \bisim^{\mathsf{X}}_{k+1} M',w' :\iff M,w \bisim_0 M',w'$ and
        \begin{itemize}
            \item {[forth]} for all $v \in R[w]$ there is a $v' \in R' [w']$ such that $M,v \bisim_k M',v'$;
            \item {[back]} for all $v' \in R'[w']$ there is a $v \in R [w]$ such that $M,v \bisim_k M',v'$,
        \end{itemize}
    \end{itemize}
and we write $M,w \bisim^{\mathsf{X}}_k M',w'$ (or simply $w\bisim^k w'$) if $(M,w)$ and $(M',w')$ are $k$-bisimilar with respect to $\mathsf{X}$. The \emph{modal depth} $md(\phi)$ of a formula $\phi$ is defined as usual; note that we put $md(\NE):=0$ and $md(\phi \intd \psi):=max\{md(\phi),md(\psi)\}$. We say that $(M,w)$ and $(M,w')$ are \emph{${\mathsf{X}},k$-equivalent}, written $M,w \equiv^{\mathsf{X}}_k M',w'$ (or simply $w\equiv_k w'$), if for all $\alpha({\mathsf{X}}) \in \ML$ with $md(\alpha)\leq k$: $M,w\vDash \alpha \iff  M',w'\vDash \alpha$. We define the \emph{$k$-th Hintikka formula} or \emph{characteristic formula} $\chi^{\mathsf{X},k}_{M,w} \in \ML$ (or simply $\chi^k_{w} $) of $(M,w)$ as follows:
    \begin{align*}
        &\chi^{\mathsf{X},0}_{M,w}   &&:= &&\bigwedge \{p \sepp p \in \mathsf{X} , w \in V(p)\} \land \bigwedge \{\lnot p \sepp p \in \mathsf{X} , w \notin V(p)\}\\
        &\chi^{{\mathsf{X}},k+1}_{M,w} &&:= &&\chi^{k}_{w}
            \land \underset{v \in R[w]}{\bigwedge}\Di \chi^{k}_{v}
            \land \Bo \underset{v \in R[w]}{\bigdis}\chi^{k}_{v}
    \end{align*}
    
    We then have that:

\begin{theorem} (See, for instance, \cite{blackburn2001,goranko2007}.
)\label{neg:theorem:hintikka_bisimulation} 
    \begin{align*}
        &w \equiv_k w'  &&\iff &&w \bisim_k w' &&\iff   &&w'\vDash \chi^{k}_{w} && \iff && \chi^k_w \equiv \chi^k_{w'}
    \end{align*}
\end{theorem}

 Similarly, a \emph{pointed (team) model} over $\mathsf{X}$ is a pair $(M,s)$ where $M$ is a model over $\mathsf{X}$ and $s$ is a team on $M$. For a given $k\in \mathbb{N}$, $(M,s)$ and $(M,s')$ are \emph{$\mathsf{X},k$-equivalent} in a logic $L$, written $M,s \equiv^\mathsf{X}_k M',s'$, if for all $\phi(\mathsf{X})\in L$ with $md(\phi)\leq k: M,s\vDash \phi \iff M',s'\vDash \phi$. It can be shown that $s \equiv^\mathsf{X}_k s'$ in $\BSML$ iff $s \equiv^\mathsf{X}_k s'$ in $\BSMLI$, so we use $s \equiv^\mathsf{X}_k s'$ to refer to equivalence in either of these two logics. We say that $(M,s)$ and $(M',s')$ are \emph{(team) $k$-bisimilar} with respect to $\mathsf{X}$, and write $M,s \bisim^\mathsf{X}_k M',s'$ (or simply $s\bisim_k s'$) if for each $w \in s$ there is some $w' \in s'$ such that $w \bisim_k w'$; and for each $w' \in s'$ there is some $w \in s$ such that $w \bisim_k w'$.  We define the \emph{(weak) $k$-th Hintikka formula} $\chi^{\mathsf{X},s}_{M,s} $ (or simply $\chi^k_{s} $) of $(M,s)$ by $\chi^{\mathsf{X},k}_{M,s}:=\bigvee_{w\in s}\chi^k_{w}$. We have:
\begin{theorem}\cite{hella2014,aloni2023} \label{neg:theorem:hintikka_bisimulation_team}
    \begin{align*}
        &[s \equiv_k s'  &&\iff &&s \bisim_k s'] &&\text{and}  &&[s'\vDash \chi^{k}_{s} && \iff && s' \bisim_k t \text{ for some $t\subseteq s$}]
    \end{align*}
\end{theorem} 
 
 %We define the \emph{$k$-th Hintikka formula} $\theta^{\mathsf{X},s}_{M,s} $ (or simply $\theta^k_{s} $) by $\theta^{\mathsf{X},k}_{M,s}:=\bigvee_{w\in s}(\chi^k_{w}\land \NE)$. Then:
%\begin{theorem} \cite{aloni2023} \label{neg:theorem:hintikka_bisimulation_team}
%    \begin{align*}
%        &s \equiv_k s'  &&\iff &&s \bisim_k s' &&\iff   &&s'\vDash \theta^{k}_{s} && \iff && \theta^k_s \equiv \theta^k_{s'}
%    \end{align*}
%\end{theorem}

A \emph{(world) property} (over $\mathsf{X}$) is a class of pointed models over $\mathsf{X}$. Each formula $\alpha\in \ML$ \emph{expresses} a property $\llbracket \alpha \rrbracket_{\mathsf{X}}$ (or simply $\llbracket \alpha \rrbracket$) over any $\mathsf{X}\supseteq \mathsf{P}(\alpha)$, where
$$\llbracket \alpha \rrbracket_{\mathsf{X}}:=\{(M,w)\text{ over $\mathsf{X}$}\sepp M,w\vDash \alpha\}.$$
A \emph{(team) property} (over $\mathsf{X}$) is a class of pointed team models over $\mathsf{X}$. Each formula $\phi$ \emph{expresses} a property $\left\Vert \phi \right\Vert_{\mathsf{X}}$ (or simply $\left\Vert \phi \right\Vert$) over any $\mathsf{X}\supseteq \mathsf{P}(\phi)$, where
$$\left\Vert\phi\right\Vert_{\mathsf{X}}:=\{(M,s)\text{ over }\mathsf{X}\sepp M,s \vDash \phi\}.$$
The \emph{ground team} of a team property $\PPP$, denoted $\bigcupdot \PPP$, is the world property
$$\bigcupdot \PPP:=\{(M,w)\sepp \exists (M,s)\in \PPP: w\in s  \},$$
and the \emph{ground team} of a formula $\phi$ (over $\mathsf{X}$), denoted $|\phi|_\mathsf{X}$, is defined by $|\phi|_\mathsf{X}:=\bigcupdot \left\Vert\phi\right\Vert_{\mathsf{X}}$, so that
$$|\phi|_\mathsf{X}=\{(M,w)\sepp \exists (M,s)\in \left\Vert\phi\right\Vert_{\mathsf{X}}: w\in s \}.$$
(Note that $\bigcupdot \PPP$ is essentially a team modulo bisimulation\footnote{$\bigcupdot \PPP$ is bisimilar to the pointed model $(\biguplus \{M' \mid (M',w')\in \bigcupdot\PPP\},\{w'\mid (M',w')\in \bigcupdot\PPP  \})$, where $\biguplus$ denotes disjoint union, in the sense that for each $(M,w)\in \bigcupdot\PPP$ there is a $w'\in \{w'\mid (M',w')\in \bigcupdot\PPP  \}$ such that $(M,w)$ is bisimilar to $(\biguplus \{M' \mid (M',w')\in \bigcupdot\PPP\},w') $, and vice versa.}, whence the name.)\footnote{\label{neg:footnote:ground_team}While the other notions in Section \ref{neg:section:preliminaries} are from the literature, the term `ground team' is introduced in this paper. The notion itself is, however, essentially identical or very similar to a number of other notions in the literature. In this endnote we list some of these notions to stave off confusion, to help the reader make connections, and also to highlight the fact that ground teams (and hence potentially also the incompatibility notions making use of ground teams---see Sections \ref{neg:section:inconsistencies} and \ref{neg:section:propositional}) can (at least in some settings) be given meaningful, interesting interpretations. %(We will discuss some possible interpretations of ground teams and the incompatibility notions in more detail in Section \ref{neg:section:interpretations}.)

Hodges \cite{hodges1997} defines (a first-order analogue) of the following notion: the \emph{flattening operator} $\downarrow$ is such that $\left\Vert\downarrow\phi\right\Vert=\{(M,t)\mid \forall w\in t: \exists (M,s)\in \left\Vert \phi\right\Vert: w\in s\}$. In other words, the flattening $\downarrow \phi$ of $\phi$ expresses the property consisting of teams consisting of elements of $|\phi|$, the ground team of $\phi$. This is essentially the power set of the ground team---it is the flat team property corresponding to the world property $|\phi|$. Hodges' flattening is similar to but distinct from Väänänen's \cite{vaananen2007} \emph{flattening} $\phi^f$ of $\phi$ (see Sections \ref{neg:section:PLNE} and \ref{neg:section:propositional_dependence_logic})---$\phi^f$ is, roughly, $\phi$ with each nonclassical subformula $\psi$ replaced with a classical subformula $\alpha$ such that $\psi\vDash \alpha$ and $|\psi|=|\alpha|$. %\todoa{It is not the case that $\left\Vert\phi^f\right\Vert=\left\Vert\downarrow\phi\right\Vert$. Take $\phi:=Ra\vee\forall x\con{x}$. Then $\left\Vert\downarrow\phi\right\Vert=\left\Vert Ra\right\Vert$ and $\left\Vert\phi^f\right\Vert=\left\Vert \top \right\Vert$.}
%with all nonclassical   The inquisitive analogue of Väänänen's flattening of $\phi$ is the \emph{classical variant} of $\phi$.
%The flattening is then essentially the power set of the truth set, and in Hodges' downward closed setting, the power set of the ground team.
%We make use of this notion in the definition of flat-incompatibility; see Section \ref{neg:section}. 
%Hodges' flattening is not to be confused with Väänänen's \cite{vaananen2007} (syntactic) \emph{flattening} (see Sections \ref{neg:section:PLNE} and \ref{neg:section:propositional_dependence_logic}), which is a classical formula rather than a property. The inquisitive analogue of Väänänen's flattening of $\phi$ is the \emph{classical variant} of $\phi$.

In inquisitive semantics/logic (see \cite{CiardelliRoelofsen2011,ciardelli2014,inqsembook,ciardellibook}, as well as Section \ref{neg:section:classical_logics}; we consider here only the propositional setting for simplicity), the union $\bigcup \PPP$ of a property $\PPP$ is referred to and represents the \emph{informative content} of the (inquisitive) proposition represented by $\PPP$ (in the propositional setting---see Section \ref{neg:section:propositional}---we define $\bigcupdot \PPP:=\bigcup \PPP$, so $\bigcup\PPP$ is the same notion as the ground team of $\PPP$). Similarly, the informative content of $\phi$ is $\bigcup \left\Vert\phi\right\Vert$. There is an operator on propositions $!$ such that $!\PPP=\wp(\bigcup \PPP)$ and a corresponding connective $!$ such that $\left\Vert!\phi\right\Vert=\wp(\bigcup\left\Vert\phi\right\Vert)$, and hence $!\phi$ is equivalent to the Hodges' flattening $\downarrow \phi$ of $\phi$. Applied to $\PPP$/$\phi$, the operator/connective yields a proposition/formula that has the same informative content as $\PPP$/$\phi$, but with no \emph{inquisitive content}. The set of worlds $w$ such that $\{w\}\vDash\phi$ is referred to as the \emph{truth-set} of $\phi$ (as in $\BSML$, the fundamental semantic notion, defined with respect to teams, is \emph{support}; truth is defined as support with respect to singleton teams). Due to the downward closure of inquisitive logic, the truth-set of $\phi$ is equal to its ground team/informative content. The truth-set/informative content of $\phi$ can also represent the information that $\phi$ presupposes, if $\phi$ expresses a question. Inquisitive logic also features a notion analogous to Väänänen's flattening of $\phi$; this is called the \emph{classical variant} of $\phi$.} Observe that for classical $\alpha$, $\llbracket \alpha\rrbracket =|\alpha|$. 
%Given a model $M$, a formula $\phi$ and $\alpha\in\ML$, we also define the \emph{local world property} $\llbracket \alpha \rrbracket_M :=\{ w\mid M,w\vDash \alpha\}$, the \emph{local team property} $\left\Vert \phi \right\Vert_M :=\{ s\mid M,s\vDash \phi\}$ and the \emph{local ground team} $|\phi|_M:=\bigcupdot \left\Vert\phi\right\Vert_M$.

We say that a logic $L$ is \emph{expressively complete} for a class of (team) properties $\mathbb{P}$, written $\left\Vert L\right\Vert=\mathbb{P}$, if for each finite $\mathsf{X}$, the class $\mathbb{P}_\mathsf{X}$ of properties over $\mathsf{X}$ in $\mathbb{P}$ is precisely the class of properties over $\mathsf{X}$ expressible by formulas of $L$, that is, if
$$\left\Vert L \right\Vert_\mathsf{X}=\mathbb{P}_\mathsf{X}, \text{ where }\left\Vert L \right\Vert_\mathsf{X}:=\{\left\Vert \phi \right\Vert_\mathsf{X}\sepp \phi \in L\}.$$
The definition of closure properties (\ref{neg:def:closure_properties}) is extended to team properties in the obvious way; for instance, a property $\PPP$ is \emph{union closed} just in case $ [(M,s)\in \PPP $ for all  $s\in S\neq \emptyset]$ implies $ (M, \bigcup S  )\in \PPP$. We additionally say that $\PPP$ is \emph{invariant under $\mathsf{X},k$-bisimulation ($k\in \N$)} if $[(M,s)\in \PPP $ and $M,s\bisim^\mathsf{X}_k M',s']$ implies  $(M',s')\in \PPP $, and that $\PPP$ over $\mathsf{X}$ is \emph{invariant under bounded bisimulation} if $\PPP$  is invariant under $\mathsf{X},k$-bisimulation for some $k\in \N$.

\begin{theorem}[Expressive completeness of $\BSMLBF$ and $\BSMLIBF$] \cite{aknudstorp2024,aloni2023} \label{neg:theorem:BSML_expressive_completeness}
\begin{enumerate}
    \item[(i)] \makeatletter\def\@currentlabel{(i)}\makeatother\label{neg:theorem:BSML_expressive_completeness_i} $\BSML$ is expressively complete for
    $\{\PPP \mid \PPP $ is convex, union closed, and invariant under bounded bisimulation$\}$.
    \item[(ii)] \makeatletter\def\@currentlabel{(ii)}\makeatother\label{neg:theorem:BSML_expressive_completeness_ii} $\BSMLI$ is expressively complete for $\{\PPP \mid \PPP $ is invariant under bounded bisimulation$\}$.
\end{enumerate}
\end{theorem}

%It follows from this theorem that each formula (in any of our logics) is equivalent to a $\BSMLI$-formula in disjunctive normal form: $\phi\equiv \bigintd_{(M,s)\in \left\Vert\phi\right\Vert}\theta^{\mathsf{P}(\phi),md(\phi)}_s$.

%refers to the analogue of truth-set of $\phi$ in his setting as the \emph{flattening} of $\phi$ (Hodges' setting is downward closed, so this is again equal to the ground team of $\phi$). 

%% file: sections/inconsistencies.tex
\subsection{Bicompleteness and Incompatibility} \label{neg:section:inconsistencies}
In this section, we reformulate Burgess' theorem as an expressive completeness result for pairs satisfying a specific notion of incompatibility. We then distinguish between multiple different notions of incompatibility, examine how they are related, and show that some of them are not sufficiently strong to yield Burgess theorems for $\BSML$ and $\BSMLI$ due to the failure of downward closure and the empty team property in these logics.

\begin{figure}[t]
\begin{tikzpicture}
\node[] (flat) [] {F-I};
\node[] (intcon) [right=of flat] {D-I} ;
\node[] (intcon2) [right=of intcon] {E-D-I} ;
%\node[] (flatcon1) [right=of flat] {F-C1} ;
%\node[] (flatcon) [right=of flatcon1] {F-C2} ;
%\node[] (intcon) [right=of world] {Intuitionistically I} ;
\node[] (gc) [right=of intcon2] {G-C} ;
\node[] (gcmodbot) [right=of gc] {G-C mod $\Bot$} ;
\node[] (world) [below=of gc] {W-I} ;
\node[] (gi) [right=of gcmodbot] {G-I} ;
\node[] (bot) [right=of gi] {$\bot$-I} ;
\node[] (empty) [above=of bot] {$\NE$-I};
\node[] (Bot) [below=of bot] {$\Bot$-I} ;
\node[] (team) [left=of Bot] {T-I} ;

\path[->] (flat) edge[color=black] (intcon);
\path[->] (flat) edge[color=black,bend left=17] (gc);
\path[->] (intcon2) edge[color=black] (world);
\path[->] (intcon) edge[color=black] (intcon2);
\path[->] (intcon2) edge[color=blue,dashed] node[below] {\tiny DC} (gc);
%\path[->] (flatcon1) edge[color=black, bend left=5] (world);
%\path[->] (flat) edge[color=black] (flatcon1);
%\path[->] (flatcon1) edge[color=black,bend right=10]  (flatcon);
%\path[->] (flatcon) edge[color=blue,dashed,bend right=10] node[above] {\tiny DC} (flatcon1);
%\path[->] (flatcon) edge[color=black] (intcon);
%\path[->] (world) edge[color=black] (intcon);
%\path[->] (flatcon) edge[color=black] (gc);
\path[->] (gc) edge[color=black,bend right=10] (gcmodbot);
\path[->] (gcmodbot) edge[color=red,bend right=10,dashed] node[above] {\tiny ET} (gc);
\path[->] (gcmodbot) edge[color=black] (gi);
\path[->] (gi) edge[color=black,bend right=10] (bot);
\path[->] (bot) edge[color=red,bend right=10,dashed] node[right] {\tiny ET} (empty);
\path[->] (empty) edge[color=teal,dashed] node[above left] {\tiny Conv.} (gi);
\path[->] (empty) edge[color=black, bend right=10] (bot);
\path[->] (Bot) edge[color=black] (bot);
\path[->] (flat) edge[color=black,bend left=10] (empty);
\path[->] (intcon2) edge[color=black,bend right=22] (bot);
\path[->] (bot) edge[color=blue,bend right=10,dashed] node[above] {\tiny DC} (gi);
\path[->] (team) edge[color=black] (Bot);
%\path[->] (bot) edge[color=black,bend right=10] (intcon);
\path[->] (world) edge[color=violet,dashed,bend left=10] node[below left] {\tiny UC \& ET} (flat);
\path[<->] (gc) edge[color=blue,dashed] node[left] {\tiny DC} (world);
%\path[->] (world) edge[color=blue,dashed] node[below] {\tiny DC} (gi);
%\path[->] (team) edge[color=red] node[below] {\tiny ET} (gc);
%\path[->] (Bot) edge[color=red,bend right=10] node[above] {\tiny ET} (team);
\end{tikzpicture}
    \caption{Some implications among the incompatibility notions. Black arrows indicate implications which always hold. Dashed arrows in other colors indicate implications conditional on closure properties (e.g., DC = downward closure, ET = empty team property, etc.). Implications which trivialize the notions have been omitted; for instance, in a setting with the empty team property, no pair can be team-incompatible, so team-incompatibility with the empty team property implies all other notions. Observe also that if a pair $\phi_0,\phi_1$ is $\NE$-incompatible, then both formulas have the empty team property. Similarly flat-incompatibility forces the flatness of both $\phi_0$ and $\phi_1$, while D-I forces the downward closure and the empty team property of $\phi_1$. Each of team-incompatibility and $\Bot$-incompatibility is inconsistent with both formulas having the empty team property.}
       \label{neg:figure:inconsistency_relations} 
    	\end{figure}
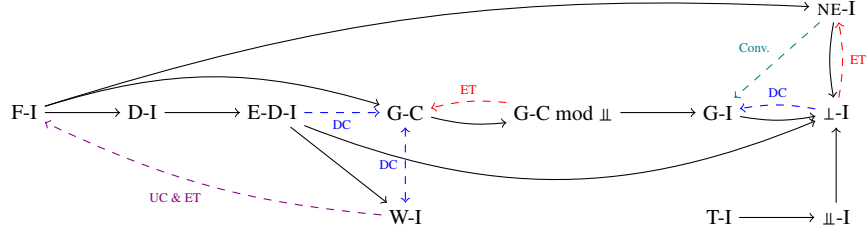

Burgess \cite{burgess2003} formulates his result as follows: if two sentences $\phi$ and $\psi$ are incompatible in that they share no models, then there is a sentence $\theta$ such that $\phi\equiv \theta$ and $\psi\equiv \lnot\theta$. Kontinen and Väänänen \cite{kontinen2011} (working in D and IF, whose formulas are downward closed and have the empty team property), make use of an equivalent notion of incompatibility, generalized to make it applicable to arbitrary formulas (that may contain free variables): the formulas $\phi$ and $\psi$ are incompatible if, for any given model and (first-order) team on that model, if $\phi$ and $\psi$ are both true in the team, the team must be empty.\footnote{For sentences, Kontinen and Väänänen's notion is equivalent to Burgess'. This follows from %(i) the empty team property; and (ii) 
the fact that a sentence, in first-order team semantics, is defined to be true in a model just in case it is true in all first-order teams on the model, and this is equivalent to it being true in some nonempty team on the model.} Let us say that two formulas which are incompatible in this way are \emph{$\bot$-incompatible}. In our modal setting:
\begin{definition}[$\bot$-incompatibility] \label{neg:def:bot_incompatibility}
$\phi_0$ and $\phi_1$ are \emph{$\bot$-incompatible ($\bot$-I)} if $[s\vDash \phi_1$ and $s\vDash \phi_1]$ implies $s=\emptyset$ (equivalently, if $\phi_0,\phi_1\vDash \bot$).
\end{definition}
Table \ref{neg:table:inconsistencies} lists all the incompatibility notions we consider in this paper, formulated for the simpler propositional setting (see Section \ref{neg:section:propositional}). Figure \ref{neg:figure:inconsistency_relations} displays some of the implications between the incompatibility notions.

    \begin{table}[t] 
    \begin{center}
        \begin{tabular}{ l | l | l } 
            Incompatibility notion & Definition(s) & \texttt{Bicomplete} logics\\
            \hline\hline
                \hyperref[neg:def:bot_incompatibility]{$\bot$-incompatible ($\bot$-I)} & $\bullet$ $[s\vDash \phi_0$ and $s\vDash \phi_1]\implies s=\emptyset$ & \hyperref[neg:theorem:negation_IF_D]{D}, \hyperref[neg:theorem:negation_IF_D]{IF} \\
                        & $\bullet$ $\phi_0,\phi_1\vDash \bot$ & \\
            & $\bullet$ $\NE$-incompatible or $\Bot$-incompatible & \\
            \hline
            \hyperref[neg:def:ground_incompatibility]{Ground-incompatible (G-I)} & $\bullet$ $[s\vDash \phi_0$ and $t\vDash \phi_1]\implies s\cap t=\emptyset$ & \hyperref[neg:theorem:negation_IF_D]{D}, \hyperref[neg:theorem:negation_IF_D]{IF}, \hyperref[neg:coro:bsml_negation_completeness]{$\BSML$},\\
            &  $\bullet$ $|\phi_0|\cap |\phi_1|=\emptyset$ &  \hyperref[neg:coro:bsml_negation_completeness]{$\BSMLI$}, \hyperref[neg:coro:PLVVEENE_negation_completeness]{$\PL(\NE,\vvee)$} \\
            \hline
        \hyperref[neg:def:Bot_NE_incompatibility]{$\Bot$-incompatible ($\Bot$-I)} & $\bullet$ $\phi_0$ and $\phi_1$ never jointly true  &  \\
             & $\bullet$ $\phi_0,\phi_1\vDash \Bot$ & \\
            & $\bullet$ $\left\Vert \phi_0 \right\Vert \cap \left\Vert \phi_1 \right\Vert  =\emptyset = \left\Vert \Bot \right\Vert$\\
            \hline
            \hyperref[neg:def:Bot_NE_incompatibility]{$\NE$-incompatible ($\NE$-I)} &$\bullet$ $[s\vDash \phi_0$ and $s\vDash \phi_1]\iff s=\emptyset$&  \hyperref[neg:theorem:negation_IF_D]{D}, \hyperref[neg:theorem:negation_IF_D]{IF}\\
            & $\bullet$ $\phi_0,\phi_1\vDash \bot$ and $\phi_0,\phi_1\nvDash \Bot$\\
            & $\bullet$ $\left\Vert \phi_0 \right\Vert \cap \left\Vert \phi_1 \right\Vert  =\{\emptyset\} = \left\Vert \bot \right\Vert$\\
            \hline
    \hyperref[neg:def:world_team_incompatibility]{World-incompatible (W-I)} & $\bullet$ $\{w\}\vDash \phi_0\iff \{w\}\nvDash \phi_1$& \hyperref[neg:coro:dep_negation_completeness]{$\PL(\con{\cdot})$}, \hyperref[neg:coro:PL_negation_completeness]{$\PL$} \\
            \hline
        \hyperref[neg:def:world_team_incompatibility]{Team-incompatible (T-I)} & $\bullet$ $s\vDash \phi_i\iff s\nvDash \phi_{1-i}$  & \hyperref[neg:coro:PL_negation_completeness]{$\PL(\sim)$ (w.r.t. $\sim$)}\\
            & $\bullet$ $\left\Vert\phi_i\right\Vert=\left\Vert \top\right\Vert\setminus \left\Vert \phi_{1-i}\right\Vert$ & \\
            \hline
                        \hyperref[neg:def:flat_incompatibility]{Flat-incompatible (F-I)} & $\bullet$  world-incompatible and $\phi_0,\phi_1$ flat& \hyperref[neg:coro:PL_negation_completeness]{$\PL$ }\\
            & $\bullet$ $\left\Vert\phi_i\right\Vert=\wp(|\top|\setminus |\phi_{1-i}|) $& \\
             & $\bullet$ $\left \Vert \phi_i\right \Vert=\{s \mid t\vDash \phi_{1-i} \implies s \cap t=\emptyset\}$ &\\
            & $\bullet$ $\left \Vert \phi_i\right\Vert=\bigcup\{\PPP\subseteq \left\Vert \top\right \Vert\mid \PPP,\left\Vert \phi_{1-i}\right\Vert$ G-I$\}$ &\\
            \hline 
            \hyperref[neg:def:down_set_incompatibility]{$\phi_1$ down-set incompatibility} &  $\bullet$ $s\vDash \phi_1\iff $ & \hyperref[neg:coro:PL_negation_completeness]{$\mathrm{InqB}$ (w.r.t. $\lnot_i$)}, \hyperref[neg:coro:PL_negation_completeness]{$\PL$} \\
             \hyperref[neg:def:down_set_incompatibility]{(D-I) of $\phi_0$}   & \quad$[[t\vDash \phi_{0} \text{ and }t\subseteq s]\implies t=\emptyset]$ &\\
             \hline
            \hyperref[neg:def:two_side_down_set_incompatibility]{Down-set incompatibilities} &  $\bullet$ $\phi_1$ D-I of $\phi_0$ or $\phi_0$ D-I of $\phi_1$ & \hyperref[neg:coro:HS_negation_completeness]{$\HS$}, \hyperref[neg:coro:PL_negation_completeness]{$\PL$} \\
             \hyperref[neg:def:two_side_down_set_incompatibility]{(on either side) (E-D-I)} &  &\\
             \hline
             \hyperref[neg:def:ground-complementary]{Ground-complementary (G-C)} & $\bullet$ $|\phi_i|=|\top|\setminus|\phi_{1-i}|$ & \hyperref[neg:coro:dep_negation_completeness]{$\PL(\con{\cdot})$}, \hyperref[neg:coro:PL_negation_completeness]{$\PL$}\\
            \hline
            \hyperref[neg:def:ground-complementary]{Ground-complementary mod $\Bot$} & $\bullet$ $|\phi_i|=|\top|\setminus|\phi_{1-i}|$ or $\phi_0 \equiv \Bot$ or $\phi_1 \equiv \Bot$ & \hyperref[neg:coro:PLNE_negation_completeness]{$\PL(\NE)$}, \hyperref[neg:coro:dep_negation_completeness]{$\PL(\con{\cdot})$}, \hyperref[neg:coro:PL_negation_completeness]{$\PL$}\\
            \hline
            %\todoa{Intuitionistically}  & $\{w\}\vDash \phi_i\implies \{w\}\nvDash \phi_{1-i}$& $\mathrm{InqB}$ \\
            %\todoa{incompatible?} & &\\
            %\hline
           % $\phi_1$ flat contradictory (2) & $\bullet$ $\left\Vert\phi_1\right\Vert=\wp(|\top|\setminus |\phi_{0}|) $& $\mathrm{InqB}$ (w.r.t. $\lnot_i$), $\PL$ \\
           %  of $\phi_0$ (F-C2)  & &\\
           %  \hline
            %$\phi_1$ flat contradictory (1) & $\bullet$  world-incompatible and $\phi_1$ flat & $\mathrm{InqB}$ (w.r.t. $\lnot_i$), $\PL$ \\
            % of $\phi_0$ (F-C1)  & &\\
            %\hline 
            All pairs &   & \hyperref[neg:coro:NEstar_negation_completeness]{$\PL(\NE^*,\vvee)$} \\
            \hline
        \end{tabular}
    \end{center}
    \caption{Incompatibility notions/pair properties, formulated for the propositional setting. Each bulleted ($\bullet$) item is a definition for the relevant notion; all definitions for a given notion are equivalent. %\todoa{missing at least some inqb completeness results. remove IF?}
    }
\label{neg:table:inconsistencies} 
    \end{table}

    We noted in Section \ref{neg:section:introduction} that a sentence (of D or IF) and its dual negation share no models. The analogous result also holds for arbitrary formulas: any formula $\phi$ and its dual negation $\lnot \phi$ are $\bot$-incompatible: $\phi,\lnot\phi\vDash \bot$. Observe that this fact immediately yields the converse of Kontinen and Väänänen's Burgess theorem: for any formulas $\phi$ and $\psi$, if there is a formula $\theta$ such that $\phi\equiv \theta$ and $\psi\equiv \lnot \theta$, then since $\theta$ and $\lnot \theta$ are $\bot$-incompatible, $\phi$ and $\psi$ must also be $\bot$-incompatible.

    We may view a Burgess theorem, together with its converse, as a type of expressive completeness result. In our modal setting, we let $\left\Vert \phi\right\Vert^{\pm,\neg}_\mathsf{X}$ (or simply $\left\Vert \phi\right\Vert^{\pm}_\mathsf{X}$) denote the pair $(\left\Vert \phi\right \Vert_\mathsf{X},\left\Vert\lnot \phi\right \Vert_\mathsf{X})$. A \emph{pair property} $\mathscr{P}$ over $\mathsf{X}$ is a class of pairs of team properties $(\PPP,\QQQ)$, where $\PPP$ and $\QQQ$ are over $\mathsf{X}$. We say that a logic $\LOGIC$ is \emph{bicomplete} for a pair property $\mathscr{P}$ (with respect to $\lnot$), written $\left\Vert L\right\Vert^{\pm,\lnot}=\mathscr{P}$ (or simply $\left\Vert L\right\Vert^{\pm}=\mathscr{P}$), if for all finite $\mathsf{X}\subseteq\mathsf{Prop}$,
    
    \begin{minipage}{\textwidth}
    \begin{align*}
    \noalign{\centering $\left\Vert \LOGIC\right\Vert^{\pm,\lnot}_\mathsf{X}=\mathscr{P}_\mathsf{X},$ where\\
        $\left\Vert \LOGIC\right\Vert^{\pm,\lnot}_\mathsf{X}:=\{\left\Vert \phi\right\Vert^{\pm,\lnot}_\mathsf{X}\mid \phi\in \LOGIC\}\text{ and }\mathscr{P}_\mathsf{X}:=\{(\PPP,\QQQ)\in \mathscr{P}\mid \PPP,\QQQ\text{ over }\mathsf{X}\}$,}
    \end{align*}
    \end{minipage}
and we say it is \emph{bicomplete modulo expressive power} (or simply {\texttt{bicomplete}}, written in monospace font) for $\mathscr{P}$ (with respect to $\lnot)$ if it is bicomplete for $\mathscr{P}_{\LOGIC}:=\mathscr{P}\cap \{(\left\Vert \phi\right\Vert,\left\Vert \psi\right\Vert)\mid \phi,\psi\in L\}$ (with respect to $\lnot$). A Burgess theorem for an incompatibility notion (pair property) $\mathscr{P}$ yields an inclusion (abusing our notation) $\mathscr{P}_L \subseteq \left\Vert L\right\Vert^\pm$, and its converse yields the converse inclusion $\left\Vert L\right\Vert^\pm \subseteq \mathscr{P}_L$; together, then, they constitute a \texttt{bicompleteness} theorem. Say that $L_1$ and $L_2$ are \emph{(expressively) bi-equivalent} or \emph{strongly equivalent} (with respect to $\lnot_1$ and $\lnot_2$) if $\left\Vert L_1\right\Vert^{\pm,\lnot_1}_\mathsf{X}=\left\Vert L_2\right\Vert^{\pm,\lnot_2}_\mathsf{X}$ for all finite $\mathsf{X}\subseteq\mathsf{Prop}$.

We extend the definition of $\bot$-incompatibility to pairs of properties in the obvious way: $\PPP$ and $\QQQ$ are \emph{$\bot$-incompatible} if $(M,s)\in \PPP \cap \QQQ$ implies $s=\emptyset$. The first-order analogues of these notions (formulated for all formulas rather than only sentences) give us the following reformulation of Kontinen and Väänänen's result:

\begin{theorem}[Bicompleteness of D and IF] \label{neg:theorem:negation_IF_D}\cite{kontinen2011}
    Each of D and IF is \texttt{bicomplete} for $\bot$-incompatible pairs.
\end{theorem}

Let us now consider the logics $\BSML$ and $\BSMLI$. While it does hold, in these logics, that a formula and its dual negation are $\bot$-incompatible, $\bot$-incompatibility does not yield a Burgess theorem or bicompleteness theorem for either of these logics. To show why, we make use of the following proposition:
\begin{proposition} \label{neg:prop:incompatibility}
    If $s\vDash \phi$ and $t\vDash \lnot \phi$, then $s\cap t=\emptyset$.
\end{proposition}
\begin{proof}
    By induction on the complexity of $\phi$, which we may assume to be in negation normal form. We only show some cases; for the rest, see \cite[Proposition 3.3.9]{anttila2021}.

    If $\phi=\psi\vee \chi$, we have $s=s_1\cup s_2$ where $s_1\vDash \psi$ and $s_2\vDash \chi$, and $t\Dashv \psi \vee\chi$, whence $t\Dashv \psi$ and $t\Dashv \chi$, and therefore $t\vDash \lnot \psi$ and $t\vDash \lnot \chi$. By the induction hypothesis, $s_1 \cap t = s_1\cap t=\emptyset$, so that also $s\cap t=(s_1\cup s_2)\cap t= (s_1\cap t)\cup (s_2\cap t)=\emptyset $.

    If $\phi=\Diamond \psi$, assume for contradiction that $s\cap t\neq \emptyset$. Let $w\in s\cap t $. By $s\vDash \Diamond \phi$, there is a nonempty $t\subseteq R[w]$ such that $t\vDash \psi$. By $t\vDash \lnot \Diamond \psi$, we have $t\Dashv \Diamond \psi$ whence $R[w]\Dashv \psi$ so that also $R[w]\vDash \lnot \psi$. But then by the induction hypothesis, $t=t\cap R[w]=\emptyset$, contradicting the fact that $t$ is nonempty.
\end{proof}

Now let $\phi:=p$ and $\psi:=((p\land \NE)\lor (\lnot p \land \NE))$. Then $\phi,\psi\vDash \Bot\vDash \bot$, so $\phi$ and $\psi$ are $\bot$-incompatible. A Burgess theorem employing $\bot$-incompatibility would then produce a $\theta$ such that $\phi\equiv \theta$ and $\psi \equiv \lnot \theta$. But consider the teams $\{w_p\}$ and $\{w_p,w_{\lnot p}\}$ (where $w_p\vDash p$, etc.). We have $\{w_p\} \vDash\phi$, so $\{w_p\}\vDash \theta$; and $\{w_p,w_{\lnot p}\} \vDash\psi$, so $\{w_p,w_{\lnot p}\}\vDash \lnot \theta$. Therefore, by Proposition \ref{neg:prop:incompatibility}, $\{w_p\}\cap \{w_p,w_{\lnot p}\}=\{w_p\}=\emptyset$, a contradiction.

Drawing on Proposition \ref{neg:prop:incompatibility}, we can define an incompatibility notion more appropriate for $\BSML$ and $\BSMLI$:
\begin{definition}[Ground-incompatibility] \label{neg:def:ground_incompatibility}
$\phi_0$ and $\phi_1$ are \emph{ground-incompatible (G-I)} if [$s\vDash \phi_0$ and $t\vDash\phi_1$] implies $s\cap t=\emptyset$ (equivalently, if $|\phi_0 |_\mathsf{X}\cap |\phi_1|_\mathsf{X}=\emptyset =|\bot|_\mathsf{X}$ for any $\mathsf{X}\supseteq \mathsf{P}(\phi_0) \cup \mathsf{P}(\phi_1)$).
\end{definition}

 We will show in Section \ref{neg:section:BSML_theorems} that ground-incompatibility does yield Burgess theorems for $\BSML$ and $\BSMLI$. It is easy to see that ground-incompatibility is strictly stronger than $\bot$-incompatibility:

\begin{fact} \label{neg:fact:ground_incompatibility_bot_incompatibility} If $\phi$ and $\psi$ are ground-incompatible, then they are $\bot$-incompatible, and the converse implication does not hold in general.
\end{fact}
\begin{proof}
Given ground-incompatible $\phi$ and $\psi$, if there is some $s$ such that $s\vDash \phi$ and $s\vDash \psi$, this implies $s=\emptyset$ whence $\phi,\psi\vDash \bot$. Otherwise, if there is no such $s$, then $\phi,\psi\vDash \Bot\vDash \bot$. Either way, $\phi,\psi\vDash \bot$.

    We have already seen that the converse does not hold (in the example directly below Proposition \ref{neg:prop:incompatibility}).
\end{proof}

The analogous fact clearly also holds in the first-order setting. This means that the Burgess theorems for D and IF also hold with respect to ground-incompatibility. In general, if a Burgess theorem holds with respect to a given incompatibility notion---a given pair property---it also holds with respect to all stronger notions: if $\mathscr{P}_1\subseteq \mathscr{P}_2$, then $\mathscr{P}_{2\LOGIC}\subseteq \left\Vert \LOGIC \right\Vert^\pm$ implies $\mathscr{P}_{1\LOGIC}\subseteq \left\Vert \LOGIC \right\Vert^\pm$. In a downward-closed setting (such as that of D and IF), ground-incompatibility is in fact equivalent to $\bot$-incompatibility:

\begin{fact} \label{neg:fact:downward-closed_incompatibility} For downward-closed $\phi$ and $\psi$: $\phi$ and $\psi$ are ground-incompatible iff they are $\bot$-incompatible. 
\end{fact}
\begin{proof}
    Given downward-closed and $\bot$-incompatible $\phi$ and $\psi$, if $s\vDash \phi$ and $t\vDash \psi$, then by downward closure, $s\cap t\vDash \phi \land \psi$, whence $s\cap t =\emptyset$.
\end{proof}

This, in turn, means that the converses of the Burgess theorems for D and IF also hold with respect to ground-incompatibility, and hence that D and IF are also \texttt{bicomplete} for ground-incompatible pairs. In general we clearly have that if the converse of a Burgess theorem holds with respect to a given notion, it also holds with respect to all weaker notions: if $\mathscr{P}_1\subseteq \mathscr{P}_2$, then $ \left\Vert \LOGIC \right\Vert^\pm \subseteq \mathscr{P}_{1\LOGIC}$ implies $\left\Vert \LOGIC \right\Vert^\pm\subseteq \mathscr{P}_{2\LOGIC}$.

So $\bot$-incompatibility and ground-incompatibility are equivalent in a downward-closed setting, but come apart when downward closure fails. Are there any natural interpretations of $\bot$-/ground-incompatibility that would also yield some intuitive understanding of this fact? We discuss this below. Before we do so, we further dissect $\bot$-incompatibility into two more notions.

%In Section \ref{neg:section:interpretations}, we provide rough sketches of some possible sets of interpretations. It will be helpful, for the discussion in Section \ref{neg:section:interpretations}, to further dissect $\bot$-incompatibility into two more natural notions of incompatibility; this we do now.

\begin{definition}[$\Bot$-incompatibility and $\NE$-incompatibility] \label{neg:def:Bot_NE_incompatibility}\
\begin{enumerate}
    \item[(i)] $\phi_0$ and $\phi_1$ are \emph{$\Bot$-incompatible ($\Bot$-I)} if there is no team (including the empty team) that makes both $\phi_0$ and $\phi_1$ true (equivalently, if $\phi_0,\phi_1\vDash \Bot $; equivalently, if $\left\Vert\phi_0 \right\Vert_\mathsf{X}\cap \left\Vert\phi_1\right\Vert_\mathsf{X}=\emptyset=\left\Vert\Bot\right\Vert_\mathsf{X}$ for any $\mathsf{X}\supseteq \mathsf{P}(\phi_0)\cup \mathsf{P}(\phi_1)$).
    \item[(ii)] \emph{Incompatible in non-empty teams} or \emph{$\NE$-incompatible ($\NE$-I)} if $[[s\vDash \phi_0$ and $s\vDash \phi_1]$ if and only if $s=\emptyset]$ (equivalently, if $\phi_0,\phi_1\vDash \bot$ and $\phi_0,\phi_1\nvDash \Bot$; equivalently, if $\left\Vert\phi_0 \right\Vert_\mathsf{X}\cap \left\Vert\phi_1\right\Vert_\mathsf{X}=\left\Vert\bot\right\Vert_\mathsf{X}$ for any $\mathsf{X}\supseteq \mathsf{P}(\phi_0) \cup \mathsf{P}(\phi_1)$).
    %; equivalently, if $\left\Vert\phi_0 \right\Vert_\mathsf{M}\cap \left\Vert\phi_1\right\Vert_\mathsf{M}=\{\emptyset\}$ for any relevant $M$).
\end{enumerate}
\end{definition}
Clearly $\phi$ and $\psi$ are $\bot$-incompatible iff they are $\Bot$-incompatible or they are $\NE$-incompatible, and they are $\NE$-incompatible iff they are $\bot$-incompatible and each of them has the empty team property. Therefore, in a setting with both downward closure and the empty team property such as that of D/IF, ground-incompatibility, $\bot$-incompatibility and $\NE$-incompatibility are all equivalent (and so D and IF are also \texttt{bicomplete} for $\NE$-incompatible pairs). We further have:

\begin{fact} \label{neg:fact:convex_contradictoriness}
    For convex $\phi$ and $\psi$, if $\phi$ and $\psi$ are $\NE$-incompatible, then they are ground-incompatible. $\NE$-incompatibility does not imply ground-incompatibility in general.
\end{fact}
\begin{proof} If $\phi$ and $\psi$ are $\NE$-incompatible, then each has the empty team property. The empty team property together with convexity implies downward closure, so the first claim now follows by Fact \ref{neg:fact:downward-closed_incompatibility}. For the second claim, consider the formulas $\lnot p \vvee ((p\land \NE)\vee (\lnot p\land \NE))$ and $p$.
\end{proof}

In particular, $\NE$-incompatibility implies ground-incompatibility in $\BSML$. The Burgess theorem with respect to ground-incompatibility in Section \ref{neg:section:BSML_theorems} will, therefore, also imply a Burgess theorem with respect to $\NE$-incompatibility for $\BSML$. The converse of this theorem clearly does not hold, since, for instance, for $\phi:= \NE$, $\phi$ and $\lnot \phi$ are not $\NE$-incompatible.

On to our proposed interpretation (for which I draw mainly from \cite{veltman,yalcin,yalcin2011,ciardelli2014,inqsembook}; see also \cite{anttila2025} for more extensive discussion of interpretations of the incompatibility notions). A common interpretation of team logics is one in which possible worlds/propositional valuations represent possible states of affairs, with the teams comprising them representing \emph{information states}: if one knows that the actual world is in a state $t$, one knows that whatever holds in all worlds in $t$ is true. We can express epistemic uncertainty under this interpretation as follows: define the \emph{epistemic might operator} $\filleddiamond$ by $\filleddiamond \phi:=(\phi \land \NE)\vee \top$. Then $t\vDash \phi $ iff $\exists s\subseteq t$ such that $s\neq \emptyset$ and $s\vDash \phi$---for instance, $t$ supports the assertion that it might be raining ($t\vDash \filleddiamond r$) just in case there is some world in $t$ in which it is raining (equivalently, some nonempty substate $s$ with $s\vDash r$). In contrast with formulas expressing epistemic information such as $\filleddiamond r$, we may think of classical formulas as expressing \emph{factual information}: information that places constraints (in the first place) on the space of possible worlds, rather than on the space of information states.

Assuming for simplicity a setting with only classical formulas and formulas such as $\filleddiamond r$ representing epistemic uncertainty, we may think of the ground team $|\phi|$ of $\phi$ as representing the factual information expressed by $\phi$: the ground-team of $\phi$ is a classical proposition which is true precisely in all the possible worlds which compose the information states in which $\phi$ in assertible, and which does not directly communicate any epistemic information. `It might be raining' expresses no factual information since $| \filleddiamond r|=|\top|$---knowing that it might be raining gives one no information as to whether it actually is raining. 

We may, accordingly, think of ground-incompatibility as contradiction in factual information: if $\phi$ and $\psi$ are ground-incompatible ($|\phi|\cap |\psi|=\emptyset$), the factual information expressed by $\phi$ rules out that expressed by $\psi$, and vice versa. On the other hand, a $\bot$-incompatibility that is not also a ground-incompatibility (which, assuming we are in a convex setting, must be a $\Bot$-incompatibility by Fact \ref{neg:fact:convex_contradictoriness}) such as the pair $r$, $\filleddiamond\lnot r$ (we have $r, \filleddiamond\lnot r\vDash \Bot \vDash \bot$ and $|r|\cap |\filleddiamond\lnot r|=|r|\neq \emptyset$) has to be a contradiction involving epistemic information, or an \emph{epistemic contradiction} (c.f. \cite{yalcin,yalcin2011}). Indeed, $r\land \filleddiamond \lnot r$ corresponds to the classic epistemic contradiction `it is raining and it might not be raining'. We then have the following intuitive explanation for Facts \ref{neg:fact:ground_incompatibility_bot_incompatibility} and \ref{neg:fact:downward-closed_incompatibility}: G-I---factual contradictoriness---and $\bot$-I---contradictoriness in some more general sense---come apart when expressive epistemic information is admitted, since in addition to conflicts in purely fact-based constraints we now also have conflicts involving epistemic constraints, which include the typical examples of epistemic contradictions. In a setting with only factual, downward-closed information the latter type of conflict does not arise, so G-I and $\bot$-I coincide.

%Returning to our interpretations of these notions, this means that in $\BSML$, pairs of formulas which are epistemically/pragmatically contradictory but not non-epistemically/pragmatically contradictory ($\bot$-incompatible but not ground-incompatible) must be $\Bot$-incompatible. Conversely, $\Bot$-incompatibility 

%% file: sections/bsml_result.tex
\subsection{The Theorems} \label{neg:section:BSML_theorems}
In this section, we prove the Burgess- and bicompleteness theorems for $\BSML$ and $\BSMLI$.

Burgess' original result is essentially a corollary of Craig's interpolation for classical first-order logic. Similarly, our results follow from interpolation for classical modal logic $\ML$ (that is, interpolation for the smallest normal modal logic K---as with the other logics in this paper, we use `$\ML$' to refer both to the language as well as the logic consisting of all $\ML$-formulas which are valid in all Kripke frames). This allows us to prove the following variant of interpolation for ground-incompatible formulas (c.f. the similar results for IF \cite[pp. 59--60]{hintikka1996} and D \cite[Theorem 6.7]{vaananen2007}; and note that restricted to $\ML$-formulas, the following is equivalent to interpolation for $\ML$/K).

\begin{theorem}[Separation theorem] \label{neg:theorem:separation}
    If $\phi$ and $\psi$ are ground-incompatible, then there is a classical formula $\gamma\in \ML$ such that $\mathsf{P}(\gamma)\subseteq \mathsf{P}(\phi)\cap \mathsf{P}(\psi)$, $\phi\vDash \gamma$, and $\psi\vDash \lnot \gamma$.
\end{theorem}
\begin{proof}
Take $\alpha:= \bigvee_{u\in \left\Vert\phi\right\Vert_{\mathsf{P}(\phi)}}\chi^{\mathsf{P}(\phi),md(\phi)}_u$. We show that $\phi\vDash \alpha$. If there are no $(M,t)$ such that $M,t\vDash \phi$, we have $\alpha=\bot$ and $\phi\equiv \Bot\vDash \bot$. Otherwise let $M,t\vDash \phi$, where $(M,t)$ is over $\mathsf{X}\supseteq \mathsf{P}(\phi)$. This implies $M',t\vDash \phi$, where $M'$ is the restriction of $M$ to $\mathsf{P}(\phi)$, i.e., $(M',t)\in \left\Vert\phi\right\Vert_{\mathsf{P}(\phi)}$. Clearly $M,t\bisim_{\mathsf{P}(\phi)}^{md(\phi)} M',t$ so $M,t\vDash \chi^{\mathsf{P}(\phi),md(\phi)}_{M',t}$ by Theorem \ref{neg:theorem:hintikka_bisimulation_team}. Then by the empty team property of $\ML$, $M,t\vDash \alpha$. Similarly, for $\beta:=\bigvee_{u\in \left\Vert\psi\right\Vert} \chi^{\mathsf{P}(\psi),md(\psi)}_u$, we have $\psi\vDash \beta$.

We now show that $\alpha$ and $\beta$ are ground-incompatible, so let $M,s\vDash \alpha$ and $M,t\vDash \beta$, where $M$ is over $\mathsf{X}\supseteq  \mathsf{P}(\phi)\cup \mathsf{P}(\psi)$, and assume for contradiction that $w\in s\cap t$. This implies that there are $(M_\phi,u_\phi) \in \left\Vert\phi\right\Vert_{\mathsf{P}(\phi)}$ and $w_\phi\in u_\phi$ such that $w\vDash \chi^{\mathsf{P}(\phi),md(\phi)}_{w_\phi}$; and $(M_\psi,u_\psi) \in \left\Vert\psi\right\Vert_{\mathsf{P}(\psi)}$ and $w_\psi\in u_\psi$ such that $w\vDash \chi^{\mathsf{P}(\psi),md(\psi)}_{w_\psi}$. By Theorem \ref{neg:theorem:hintikka_bisimulation}, $w\bisim^{\mathsf{P}(\phi)}_{md(\phi)} w_\phi$ and $w\bisim^{\mathsf{P}(\psi)}_{md(\psi)} w_\psi$. Let $(M^*_\phi,u^*_\phi)$ be a pointed model over $\mathsf{X}$ defined by taking $(M_\phi,u_\phi)$ and assigning arbitrary valuations for all $p\in \mathsf{X}\setminus\mathsf{P}(\phi)$, and similarly for $(M^*_\psi,u^*_\psi)$. Let $M'$ be the disjoint union $M\uplus M^*_\phi \uplus M^*_\psi$, and let $u_\phi':=(u^*_\phi \setminus \{w^*_\phi\})\cup \{w\}$ and $u_\psi':=(u^*_\psi \setminus \{w^*_\psi\})\cup \{w\}$. Clearly $M',u_\phi' \bisim^{\mathsf{P}(\phi)}_{md(\phi)} M_\phi,u_\phi$ and $M',u_\psi' \bisim^{\mathsf{P}(\psi)}_{md(\psi)} M_\psi, u_\psi$, whence $u_\phi'\vDash \phi$ and $u_\psi'\vDash \psi$ by Theorem \ref{neg:theorem:hintikka_bisimulation_team}. But we have $w\in u_\phi' \cap u'_\psi$, contradicting the assumption that $\phi$ and $\psi$ are ground-incompatible.

    We have shown $ |\alpha|_\mathsf{X}\cap |\beta|_\mathsf{X}=\emptyset$; therefore also $\llbracket\alpha\rrbracket_\mathsf{X}\cap \llbracket\beta\rrbracket_\mathsf{X}=\emptyset $, and so we have $\alpha \vDash \lnot \beta$. By interpolation for ML, there is an $\ML$-formula $\gamma$ such that $\alpha \vDash \gamma\vDash \lnot \beta$ and $\mathsf{P}(\gamma)=\mathsf{P}(\alpha)\cap\mathsf{P}(\beta)=\mathsf{P}(\phi)\cap\mathsf{P}(\psi)$. Then also $\phi \vDash \alpha \vDash \gamma$ and $\psi\vDash \beta\vDash \lnot \gamma$.
\end{proof}

The modal complications in the proof above might distract from the structure of the argument. In the analogous proof for the propositional setting (see Section \ref{neg:section:propositional} for details on this setting), to show $ |\alpha|_\mathsf{X}\cap |\beta|_\mathsf{X}=\emptyset$, one simply notes that $|\phi|_{\mathsf{P}(\phi)}=|\alpha|_{\mathsf{P}(\phi)}$ and $|\psi|_{\mathsf{P}(\psi)}=|\beta|_{\mathsf{P}(\psi)}$, which implies $|\phi|_{\mathsf{X}}=|\alpha|_{\mathsf{X}}$ and $|\psi|_{\mathsf{X}}=|\beta|_{\mathsf{X}}$, whence $ |\alpha|_\mathsf{X}\cap |\beta|_\mathsf{X}=\emptyset$. (In the modal case, we have $|\phi|_{\mathsf{P}(\phi)}\subseteq |\alpha|_{\mathsf{P}(\phi)}$, but $|\phi|_{\mathsf{P}(\phi)}\supseteq |\alpha|_{\mathsf{P}(\phi)}$ might only hold modulo bisimulation.)

We require one more simple lemma.

\begin{lemma} \label{neg:lemma:negation_empty_team}
    For any formula $\phi$, there is a formula $\phi'$ such that $\phi'\equiv \phi$, and $\lnot\phi'$ has the empty team property.
\end{lemma}
\begin{proof}
    Define $\phi'$ by putting $\phi$ in negation normal form and then replacing each occurrence of $\lnot \NE$ with $\bot$. An easy induction shows that $\phi'$ is as required.
\end{proof}

We can now prove the Burgess theorems for $\BSML$ and $\BSMLI$. We give one version of the proof which is completely analogous to Burgess' and which works for both logics, and then provide an alternative proof for $\BSMLI$. The alternative proof does not make essential use of modal operators, and it is easy to see that its propositional analogue establishes the analogous result for the propositional fragment of $\BSMLI$. The propositional analogue of this result does not hold for the propositional fragment of $\BSML$; see Section \ref{neg:section:PLNE}. We show both the Burgess theorems as well as their converses.

%We also show that such a formula need not be necessary for this type of result, by providing an alternative proof for $\BSMLI$ which does not make use of $\theta_0$ (a propositional analogue of the alternative proof also establishes that the analogous result holds for the propositional fragment of $\BSMLI$, whereas this is not the case for $\BSML$; see Section \ref{neg:}.)

\begin{theorem}[Burgess theorems for $\BSMLBF$ and $\BSMLIBF$] \label{neg:theorem:negation_bsml}
    In both $\BSML$ and $\BSMLI$, the following are equivalent:
    \begin{enumerate}
        \item[(i)] \makeatletter\def\@currentlabel{(i)}\makeatother\label{neg:theorem:negation_bsml_i} $\phi$ and $\psi$ are ground-incompatible.
        \item[(ii)] \makeatletter\def\@currentlabel{(ii)}\makeatother\label{neg:theorem:negation_bsml_ii} There is a formula $\theta$ such that $\phi\equiv \theta$ and $\psi \equiv \lnot \theta$ (and $\mathsf{P}(\theta)=\mathsf{P}(\phi)\cup \mathsf{P}(\psi)$).
    \end{enumerate}
\end{theorem}
\begin{proof}
\ref{neg:theorem:negation_bsml_ii} $\implies $ \ref{neg:theorem:negation_bsml_i} by Proposition \ref{neg:prop:incompatibility}. \ref{neg:theorem:negation_bsml_i} $\implies$ \ref{neg:theorem:negation_bsml_ii}: Let $\theta_0:=\Diamond (\Bot \vee \lnot \Bot) $ so that $\theta_0\equiv \Diamond \Bot \equiv \bot $ and $\lnot\theta_0\equiv \Box \lnot (\Bot \vee \lnot \Bot)\equiv \Box (\lnot \Bot \land \Bot)\equiv \Box \Bot \equiv \bot $. By Lemma \ref{neg:lemma:negation_empty_team}, let $\phi'$ be such that $\phi\equiv \phi'$ and $\lnot \phi'$ has the empty team property, and similarly for $\psi'$. Let $\phi_0:= \phi' \vee \theta_0$ and $\psi_0:= \psi '\vee \theta_0$ so that $\phi_0\equiv \phi$ and $\lnot\phi_0\equiv \bot$ (clearly $\lnot \phi_0\equiv \lnot\phi'\land \lnot\theta_0\vDash \bot$, and by the empty team property of $\lnot\phi'$, we also have $\bot \vDash \lnot\phi'\land \lnot\theta_0$), and similarly for $\psi_0$. By Theorem \ref{neg:theorem:separation}, let $\gamma \in \ML$ be such that $\phi_0\vDash \gamma$ and $\psi_0\vDash \lnot \gamma$. Finally, let $\theta:= (\phi_0\land  \gamma)\vee \lnot\psi_0$. Then:
\begin{align*}
    &\theta &&= && (\phi_0\land  \gamma)\vee \lnot\psi_0&&\equiv &&\phi_0\vee \bot &&\equiv&& \phi&&&&\\
    &\lnot\theta &&= && \lnot((\phi_0\land  \gamma)\vee \lnot\psi_0)&&\equiv &&(\lnot \phi_0\vee \lnot \gamma )\land \psi_0&&\equiv &&(\bot\vee \lnot \gamma )\land \psi_0 && \equiv &&\psi
\end{align*}
Alternatively, for $\phi,\psi\in\BSMLI$, by Theorem \ref{neg:theorem:separation} we can let $\gamma$ be such that $\phi\vDash \gamma$ and $\psi\vDash\lnot \gamma$. Let $\phi_\Bot:=\lnot \phi \vee \Bot$ and $\psi_\Bot:=\lnot\psi \vee \Bot$ so that $\phi_\Bot \equiv \Bot$ and $\lnot \phi_\Bot\equiv \phi$, and similarly for $\psi_\Bot$. Finally, let $\theta:= \lnot(\phi_\Bot\intd \lnot (\psi_\Bot\intd \gamma))$. Then:
 \begin{align*}
    &\theta&&= &&\lnot (\phi_\Bot \intd \lnot (\psi_\Bot\intd \gamma))&&\equiv &&\lnot \phi_\Bot \land (\psi_\Bot\intd \gamma)&&\equiv &&\phi \land \gamma&&\equiv &&\phi&\\
&\lnot\theta&&= &&\lnot \lnot (\phi_\Bot\intd \lnot (\psi_\Bot\intd \gamma))&&\equiv &&\phi_\Bot\intd (\lnot \psi_\Bot\land \lnot \gamma)&&\equiv &&\psi\land \lnot \gamma &&\equiv &&\psi&\qedhere
    \end{align*}
\end{proof}
Burgess' proof (as well as those of Dechesne, Kontinen and Väänänen, and Mann) likewise makes use of a formula $\theta_0$ such that $\theta_0\equiv \bot \equiv \lnot \theta_0$; let us briefly discuss the sentence used by Burgess (originally due to Väänänen) to allow for easy comparison with our $\theta_0$. Formulated in IF, this sentence is $\theta_0:=\forall x (\exists y/\forall x)(x=y)$. This is equivalent to the $\Sigma^1_1$-sentence $\exists f \forall x (x=f(x))$, and its dual negation is equivalent to the FO-sentence $\exists x\forall y (x \neq y)$. Clearly, then, $\theta_0$ is true in a model just in case the domain of the model has less than two elements; and $\lnot \theta_0$ is never true. Therefore, ignoring empty models and models of size one, we have $\theta_0\equiv \bot \equiv \lnot \theta_0$.

%let us briefly discuss Kontinen and Väänänen's simplified version of this formula here to allow for easy comparison (the original formula used by Burgess is also due to Väänänen). Formulated in D, this formula is $\theta_0:=\forall x \con{x}$. A unary dependence atom such as $\con{x}$ is called a \emph{constancy atom}---given a first-order model, $\con{x}$ is a true in a first-order team $T$ on the model (a set of partial assignments which have the same domain, and with codomain the domain of the model) just in case the value assigned to the variable $x$ is constant throughout the team (see also Section \ref{neg:section:propositional_dependence_logic}, where we discuss propositional dependence atoms). The sentence $\forall x \con{x}$ is true in a model with domain $A$ just in case the value of $x$ is constant in the team $\{(x,a)\mid a\in A\}$. This can only be the case if $A$ has less than two elements. On the other hand, we have $\lnot \con{x}\equiv \bot$, whence also $\lnot \forall x \con{x} \equiv \exists x \lnot \con{x} \equiv \exists x \bot\equiv \bot$ (see Section \ref{neg:section:propositional_dependence_logic} for discussion on why $\lnot \con{x}\equiv \bot$). Therefore, ignoring empty models and models of size one, we have $\theta_0\equiv \bot \equiv \lnot \theta_0$.\todoa{Unnecessary paragraph?}

 We extend the definition of ground-incompatibility to pairs of properties in the obvious way: $\PPP$ and $\QQQ$ are \emph{ground-incompatible} if $\bigcupdot P\cap \bigcupdot \QQQ=\emptyset$. By Theorems \ref{neg:theorem:BSML_expressive_completeness} and \ref{neg:theorem:negation_bsml},
\begin{corollary}[Bicompleteness of $\BSMLBF$ and $\BSMLIBF$]\ \label{neg:coro:bsml_negation_completeness}
    \begin{enumerate}
    \item[(i)] \makeatletter\def\@currentlabel{(i)}\makeatother\label{neg:coro:bsml_negation_completeness_i} $\BSML$ is bicomplete for
    $\{(\PPP,\QQQ) \mid \PPP,\QQQ $ are convex, union closed, and invariant under bounded bisimulation; $\PPP,\QQQ$ are ground-incompatible$\}$.
    \item[(ii)] \makeatletter\def\@currentlabel{(ii)}\makeatother\label{neg:coro:bsml_negation_completeness_ii} $\BSMLI$ is bicomplete for $\{(\PPP,\QQQ) \mid \PPP,\QQQ $ are invariant under bounded bisimulation; $\PPP,\QQQ$ are ground-incompatible$\}$.
    \item[(iii)] \makeatletter\def\@currentlabel{(iii)}\makeatother\label{neg:coro:bsml_negation_completeness_iii} Each of $\BSML$ and $\BSMLI$ is \texttt{bicomplete} for ground-incompatible pairs.
    \end{enumerate}
\end{corollary}

Note that a Burgess theorem can be viewed as a kind of strengthening of the relevant Separation theorem (\ref{neg:theorem:separation}): given ground-incompatible $\phi$ and $\psi$, there is a $\theta$ such that not only $\phi\vDash \theta$ and $\psi\vDash \lnot\theta$, but also $\theta\vDash \phi$ and $\lnot \theta\vDash \psi$.\footnote{I am grateful to Maria Aloni for pointing this out to me.} However, to get also these converse entailments, one must go from $\mathsf{P}(\theta)\subseteq \mathsf{P}(\phi)\cap \mathsf{P}(\psi)$ to $\mathsf{P}(\theta)= \mathsf{P}(\phi)\cup \mathsf{P}(\psi)$.

%% file: sections/propositional.tex
\section{Burgess Theorems for Propositional Team Logics} \label{neg:section:propositional}

In this section, we prove Burgess/bicompleteness theorems for Hawke and Steinert-Threlkeld's semantic expressivist logic for epistemic modals $\HS$ (Section \ref{neg:section:HS}); for $\PL(\NE)$ and $\PL(\NE, \intd)$, the propositional fragments of $\BSML$ and $\BSMLI$, respectively (Section \ref{neg:section:PLNE}); as well as for propositional dependence logic $\PL(\con{\cdot)}$ (Section \ref{neg:section:propositional_dependence_logic}). We also comment (Section \ref{neg:section:classical_logics}) on the bicompleteness of team logics with negations that do not exhibit failure of replacement: team-based classical propositional logic $\PL$, propositional inquisitive logic InqB, and the extension $\PL(\sim)$ of $
\PL$ with the Boolean negation $\sim$ (also known as \emph{propositional team logic}---see endnote \hyperref[neg:footnote:team_logic]{5}), and give an example (in Section \ref{neg:section:PLNE}) of a logic bicomplete for all pairs.

%Table \ref{neg:table:inconsistencies} lists the incompatibility notions we make use of (formulated for the propositional setting) as well as all negation-completeness$^*$ results we discuss (throughout the entire paper).

\subsection{Bicompleteness Without Failure of Replacement} \label{neg:section:classical_logics}

The notion of bicompleteness can also be applied to logics in which replacement of equivalents under negation does not fail. We briefly consider three such logics by way of illustration.

\begin{definition}[Syntax of $\PLBF$, $\PLBF(\sim)$, and $\INQBBF$] \label{neg:def:prop_syntax}
The set of formulas of \emph{classical propositional logic} $\PL$ is generated by:
        \begin{align*}
            \alpha ::= p \sepp \bot \sepp \bnot \alpha \sepp ( \alpha \land \alpha )  \sepp ( \alpha \dis \alpha ) 
        \end{align*}
        where $p\in \mathsf{Prop}$, and $\mathsf{Prop}$ is a countable set of propositional variables, as before.
        
        The set of formulas of \emph{propositional logic with the Boolean negation} $\PL(\sim)$ is generated by:
        \begin{align*}
            \phi ::= p \sepp \bot \sepp \bnot \alpha \sepp ( \phi \land \phi )  \sepp ( \phi \dis \phi )\mid \sim\phi 
        \end{align*}
    where $p\in \mathsf{Prop}$ and $\alpha \in \PL$.

   The set of formulas of \emph{propositional inquisitive logic} $\mathrm{InqB}$ is generated by:
        \begin{align*}
            \phi ::= p \sepp \bot \sepp ( \phi \land \phi ) \mid (\phi  \to\phi) \sepp ( \phi \vvee \phi )
        \end{align*}
    where $p\in \mathsf{Prop}$.
\end{definition}
 We use $\alpha,\beta$ to refer exclusively to formulas of $\PL$ as well as the $\vvee$-free fragment of $\mathrm{InqB}$ (\emph{classical formulas}). In $\mathrm{InqB}$, let $\lnot_{i} \phi:=\phi \to \bot$ and $\bigvvee\emptyset:=\bot$. As before, let $\top:=\lnot \bot$ or $\top:=\lnot_{i}\top$, and $\bigvee \emptyset:=\bot$.

%It should be noted that the usual formulations of all the propositional logics we consider do not include the dual negation; instead, one would typically incorporate a negation that may only occur in front of propositional variables, or one that may only occur in front of classical formulas. The resulting variants are expressively equivalent in the usual sense (of being able to express the same properties) to the present formulations, but since Burgess theorems do not hold for these variants, they are not expressively equivalent in the stronger sense of being able to express the same pairs of properties to the present ones.

A (propositional) \emph{team} is a set of valuations $s\subseteq 2^{\mathsf{X}}$, where $\mathsf{X}\subseteq \mathsf{Prop}$. We formulate the semantics of $\PL$ and $\PL(\sim)$ in terms of support/anti-support conditions as in $\BSML$. The conditions for the atoms and connectives of $\PL$ are the obvious propositional analogues of their modal support/anti-support conditions. The support condition for the Boolean negation is given by $$s\vDash \sim \phi :\iff s\nvDash \phi.$$ The Boolean negation is not permitted to occur in the scope of the dual negation, so it does not have anti-support conditions.\footnote{One could of course formulate anti-support conditions for $\sim$ (i.e., truth/support conditions for $\lnot \sim \phi$), but it is not clear what they should be. The anti-support conditions for the connectives in $\BSML$ are motivated by empirical considerations, whereas those in all other logics we consider ultimately derive from the game-theoretic semantics for IF and D. In these semantics, the dual negation corresponds to an in-game move in which the players switch their verifier/falsifier roles. On the other hand, in the extensions of IF and D with $\sim$, $\sim \phi$ is true just in case the verifier does not have a winning strategy for $\phi$. The Boolean negation is evaluated from outside the game, so there is no way to interpret $\lnot \sim \phi$ in these semantics. For more discussion, see \cite{hintikka1996, vaananen2007}. See also \cite{bellier2023}, where the extension of D with $\sim$ is provided with game-theoretic semantics in an indirect way, and Section \ref{neg:section:HS}, featuring the logic $\HS$, which has a bilateral negation and in which $\sim \phi$ is uniformly definable (using the $\HS$-definition of $\sim \phi$ from Section \ref{neg:section:HS}, we get $s\Dashv \sim \phi\iff \forall t\subseteq s: t\vDash \phi$)}. We only formulate support conditions for $\mathrm{InqB}$; the support conditions for $\to$, the \emph{intuitionistic implication}, are as follows: $$s\vDash \phi \to \psi :\iff \forall t\subseteq s: [t\vDash \phi \implies t\vDash \psi].$$
The resulting support conditions for $\lnot_i$, the \emph{intuitionistic negation}, are:
$$s\vDash \textstyle{\lnot_i}\phi \iff \forall t\subseteq s: [t\vDash \phi \implies t=\emptyset].$$

We define the propositional analogues of the team semantic closure properties, team properties, expressive completeness, and the notations $\llbracket \alpha\rrbracket_\mathsf{X}$, $\left\Vert \phi\right\Vert_\mathsf{X}$, $|\phi|_\mathsf{X}$, etc. in the obvious way. Note that now $\bigcupdot \PPP = \bigcup \PPP$. We let $\left\Vert \phi\right\Vert^{\pm,\lnot}_\mathsf{X}:=(\left\Vert\phi\right\Vert_\mathsf{X},\left\Vert\lnot \phi\right\Vert_\mathsf{X}) $, $\left\Vert \phi\right\Vert^{\pm,\sim}_\mathsf{X}:=(\left\Vert\phi\right\Vert_\mathsf{X},\left\Vert\sim \phi\right\Vert_\mathsf{X})$, and $\left\Vert \phi\right\Vert^{\pm,\lnot_i}_\mathsf{X}:=(\left\Vert\phi\right\Vert_\mathsf{X},\left\Vert\lnot_i \phi\right\Vert_\mathsf{X})$.

As with $\ML$, $\PL$-formulas are flat, and the analogue of Fact \ref{neg:fact:ML_team_classical_correspondence} holds for these formulas. $\PL$ is indeed, therefore, simply classical propositional logic with team semantics, and the negation $\lnot$ is the classical negation with respect to worlds. $\PL(\sim)$ additionally incorporates the Boolean negation, which is essentially the classical negation with respect to teams. $\mathrm{InqB}$ features a negation $\lnot_{i}$ defined in terms of the intuitionistic implication $\to$. Note that, unlike with all other negations considered in this paper, double negation elimination is not valid for $\lnot_{i}$ (as, for instance, $\lnot_{i} \lnot_{i} (p\vvee q)\equiv  p\vee q\nvDash p\vvee q $). The $\vvee$-free fragment of $\mathrm{InqB}$ is an alternative way of formulating team-based classical propositional logic; one can show that this fragment is flat and that an analogue of \ref{neg:fact:ML_team_classical_correspondence} holds. Formulas $\mathrm{InqB}$ are downward closed and have the empty team property.

The usual formulations of $\PL$ and $\PL(\sim)$ (as well as those of the logics in Sections \ref{neg:section:PLNE} and \ref{neg:section:propositional_dependence_logic}) in the literature do not include the dual negation; instead, these formulations feature what we can call the \emph{restricted classical negation} $\lnot_c$. This negation may only precede classical formulas, and its semantics are given by:
$$s\vDash \textstyle{\lnot_c} \alpha:\iff \forall w\in s: \{w\}\nvDash \alpha.$$ 
One can show that for all classical $\alpha$ (where $\alpha$ may also include $\lnot,\lnot_i,\lnot_c$): $\lnot \alpha\equiv \lnot_i\alpha\equiv \lnot_c\alpha$. These formulations are therefore expressively equivalent in the usual sense (of being able to express the same properties) to the present formulations; let us express this by writing, for instance, $\left\Vert \PL_\lnot\right\Vert=\left\Vert \PL_{\lnot_c}\right\Vert$. Since, in $\PL$ and $\PL(\sim)$, the dual negation only appears in front of classical formulas, the $\lnot_c$-formulations of these logics are also bi-equivalent to the present ones (with respect to both $\lnot,\lnot_c$ as well as $\sim,\sim$): $\left\Vert \PL_\lnot\right\Vert^{\pm,\lnot}=\left\Vert \PL_{\lnot_c}\right\Vert^{\pm,\lnot_c}$; $\left\Vert \PL(\sim)_\lnot\right\Vert^{\pm,\lnot}=\left\Vert \PL(\sim)_{\lnot_c}\right\Vert^{\pm,\lnot_c}$; and $\left\Vert \PL(\sim)_\lnot\right\Vert^{\pm,\sim}=\left\Vert \PL(\sim)_{\lnot_c}\right\Vert^{\pm,\sim}$. This second fact does not, however, hold for the logics we introduce in Sections \ref{neg:section:PLNE} and \ref{neg:section:propositional_dependence_logic}---in these logics, the dual negation may occur in front of non-classical formulas, so the dual-negation formulations allow one to capture more pairs than the $\lnot_c$-formulations. So for instance, for the propositional fragment of $\BSML$, which we denote by $\PL(\NE)$, while we do have $\left\Vert\PL(\NE)_\lnot\right\Vert=\left\Vert\PL(\NE)_{\lnot_c}\right\Vert$, it is not the case that $\left\Vert\PL(\NE)_\lnot\right\Vert^{\pm,\lnot}=\left\Vert\PL(\NE)_{\lnot_c}\right\Vert^{\pm,{\lnot_c}}$. We prove Burgess theorems for the dual-negation formulations; no such theorems hold for the $\lnot_c$-formulations (we will clarify below what we mean by `Burgess theorems' in this context).

Returning to the properties of the logics $\PL$, $\PL(\sim)$ and InqB, it can be shown:

\begin{theorem}[Expressive completeness of $\PLBF$, InqB, and $\PLBF\mathbf{(\sim)}$] \cite{CiardelliRoelofsen2011,yang2017} \label{neg:theorem:PL_expressive_completeness}
    \begin{enumerate}
    \item[(i)] \makeatletter\def\@currentlabel{(i)}\makeatother\label{neg:theorem:PL_expressive_completeness_i} $\PL$ is expressively complete for
    $\{\PPP \mid \PPP $ is flat$\}$.
        \item[(ii)] \makeatletter\def\@currentlabel{(ii)}\makeatother\label{neg:theorem:PL_expressive_completeness_ii} $\mathrm{InqB}$ is expressively complete for $\{\PPP \mid \PPP $ is downward closed and has the empty team property$\}$.
    \item[(iii)] \makeatletter\def\@currentlabel{(iii)}\makeatother\label{neg:theorem:PL_expressive_completeness_iii} $\PL(\sim)$ is expressively complete for the class of all team properties.
\end{enumerate}
\end{theorem}

We introduce the following natural incompatibility notions corresponding to the two types of classical negation $\lnot$ and $\sim$:
\begin{definition}[World-incompatibility and Team-incompatibility] \label{neg:def:world_team_incompatibility}\
\begin{enumerate}
    \item[(i)] $\phi_0$ and $\phi_1$ are \emph{world-incompatible (W-I)} if $\{w\}\vDash \phi_0\iff \{w\}\nvDash \phi_1$.
    \item[(ii)] $\phi_0$ and $\phi_1$ are \emph{team-incompatible (T-I)} if $s\vDash \phi_0 \iff s \nvDash \phi_1$ (equivalently, if $\left\Vert\phi_0\right\Vert_{\mathsf{X}}=\left\Vert \top\right\Vert_{\mathsf{X}}\setminus \left\Vert \phi_1\right\Vert_{\mathsf{X}}$ for all $\mathsf{X}\supseteq \mathsf{P}(\phi_0)\cup \mathsf{P}(\phi_1)$).
\end{enumerate}
\end{definition}

Observe that, for instance, $p$ and $\lnot p$ are world-incompatible but not team-incompatible (since $\emptyset\vDash p\land \lnot p$), whereas $p$ and $\sim p$ are both team-incompatible as well as world-incompatible. These notions yield bicompleteness theorems for $\PL$ and $\PL(\sim)$, respectively, but we also introduce a stronger notion for $\PL$; this will allow us to differentiate between the pairs expressible (modulo expressive power) in $\PL$ from those expressible in propositional dependence logic (Section \ref{neg:section:propositional_dependence_logic}). 

\begin{definition}[Flat-incompatibility] \label{neg:def:flat_incompatibility}
$\phi_0$ and $\phi_1$ are \emph{flat-incompatible (F-I)} if they are world-incompatible and flat; or, equivalently, if any of the following holds for all $\mathsf{X}\supseteq \mathsf{P}(\phi)\cup \mathsf{P}(\psi)$ and $i\in\{0,1\}$:
\begin{enumerate}
    \item[(i)] $\left\Vert\phi_i\right\Vert_{\mathsf{X}}=\wp(|\top|_{\mathsf{X}}\setminus |\phi_{1-i}|_\mathsf{X})$.
    \item[(ii)] $\left\Vert\phi_i\right\Vert_{\mathsf{X}}=\{s \subseteq 2^\mathsf{X}\mid t\vDash \phi_{1-i} \implies s \cap t=\emptyset\}$.
    \item[(iii)] $\left\Vert\phi_i\right\Vert_{\mathsf{X}}=\bigcup\{\PPP\subseteq 2^{2^\mathsf{X}}\mid \PPP$ and $\left\Vert \phi_{1-i}\right\Vert_\mathsf{X}$ are ground-incompatible$\}$.
\end{enumerate}
\end{definition}

Observe that $p$ and $\lnot p$ are flat-incompatible, whereas $p$ and $\sim p$ are not.

%In a downward-closed setting such as that of InqB, the intuitionistic negation always expresses a flat property: $\left\Vert \lnot_i\phi \right\Vert =\wp(|\top|\setminus|\phi|)$. Notice that the pair property 

As for the intuitionistic negation $\lnot_i$, the truth conditions yield the following incompatibility notion: 
\begin{definition}[(One-sided) Down-set Incompatibility] \label{neg:def:down_set_incompatibility}
$\phi_1$ is a \emph{down-set incompatibility (D-I)} of $\phi_0$ if [$s\vDash \phi_1$ iff [for all $t\subseteq s$, $t\vDash \phi_0$ implies $t=\emptyset$]].
\end{definition}
 Unlike the other notions we consider, this is not symmetric---note, for instance, that $\lnot p \land \lnot q$ (which is equivalent to $\lnot_i (p\vvee q)$) is a down-set incompatibility of $p\vvee q$, but not vice versa: we have that [$t\subseteq \{w_{p,\lnot q},w_{\lnot p,q}\}$ and $t\vDash \lnot p\land \lnot q$] implies $t=\emptyset$, but $\{w_{p,\lnot q},w_{\lnot p,q}\}\nvDash p\vvee q$. This reflects the fact that double negation elimination is not sound for $\lnot_i$. It is easy to check that D-I implies world-incompatibility, and that in a flat setting (such as that of $\PL$), flat-incompatibility, D-I, and world-incompatibility are equivalent; we will therefore have that $\PL$ is \texttt{bicomplete} for each of the corresponding pair properties. That D-I implies W-I tracks the fact that in InqB the intuitionistic negation $\lnot_i$ behaves classically with respect to worlds/singletons; for instance, double negation elimination holds with respect to singletons: $\{w\}\vDash\lnot_i\lnot_i\phi\iff \{w\}\vDash \phi$.%\footnote{As an alternative to D-I, one could also introduce, for $\lnot_i$ in InqB, a one-sided version of flat-incompatibility. This would make it clearer that in InqB, $\lnot_i \phi$ is always flat, and thus also that in InqB, $\lnot_i$ is very similar to the classical $\lnot$, save for the fact that double negation elimination does not hold. For instance, say $\phi_1$ is a \emph{flat contradictory} (in sense 1) (\emph{F-C1}) of $\phi_0$ if $\phi_0$ and $\phi_1$ are world-incompatible and $\phi_1$ is flat. Say that $\phi_1$ is a \emph{flat contradictory} (in sense 2) (\emph{F-C2}) of $\phi_0$ if $\left\Vert\phi_1\right\Vert_{\mathsf{X}}=\wp(|\top|_{\mathsf{X}}\setminus |\phi_{0}|_\mathsf{X})$ for all $\mathsf{X}\supseteq \mathsf{P}(\phi)\cup \mathsf{P}(\psi)$. In a downward-closed setting with the empty team property, D-I, F-C1 and F-C2 are equivalent. We therefore have that InqB and $\PL$ are \texttt{negation-complete} for both F-C1 as well as F-C2 pairs.}

Extending these incompatibility notions to pairs of properties in the obvious way, it is easy to see that the following bicompleteness results hold:

%It is easy to show that a ``Burgess theorem'' employing classical incompatibility (and formulated using $\sim$ rather than $\lnot$), as well as its converse, hold for $\PL(\sim)$; and similarly for $\PL$ and flat-incompatibility with DNE and for $\mathrm{InqB}$ and flat-incompatibility. We therefore have:
\begin{corollary}[Bicompleteness of $\PLBF$, InqB, and $\PLBF\mathbf{(\sim)}$]\ \label{neg:coro:PL_negation_completeness}
        \begin{enumerate}
    \item[(i)]\makeatletter\def\@currentlabel{(i)}\makeatother\label{neg:coro:PL_negation_completeness_i} $\PL$ is bicomplete for F-I/W-I/D-I pairs of flat properties and therefore \texttt{bicomplete} for F-I/W-I/D-I pairs (with respect to $\lnot$).
    \item[(ii)]\makeatletter\def\@currentlabel{(ii)}\makeatother\label{neg:coro:PL_negation_completeness_ii} $\mathrm{InqB}$ is bicomplete for D-I pairs of downward closed properties with the empty team property, and therefore \texttt{bicomplete} for D-I pairs (with respect to $\lnot_i$).
    \item[(iii)]\makeatletter\def\@currentlabel{(iii)}\makeatother\label{neg:coro:PL_negation_completeness_iii} $\PL(\sim)$ is both bicomplete and \texttt{bicomplete} for T-I pairs (with respect to $\sim$).
\end{enumerate}
\end{corollary}

Inclusions such as $\{(\PPP,\QQQ)\mid (\PPP,\QQQ)$ world-incompatible$\}_{\PL}\subseteq \left\Vert \PL\right\Vert^{\pm,\lnot}$, which form part of the \texttt{bicompleteness} facts above, are of the same form as Burgess theorems. However, world-incompatibility relativized to pairs of flat properties is a pair property such that the first element of the pairs $(\PPP,\QQQ)$ in the property determines the second element: $\QQQ=\wp(\left\Vert\top\right\Vert\setminus \bigcup\PPP)$. Given the converse inclusion, then, this inclusion is trivial: if $(\PPP,\QQQ)\in \{(\PPP,\QQQ)\mid (\PPP,\QQQ)$ world-incompatible$\}_{\PL}$, then $\QQQ=\wp(\left\Vert\top\right\Vert\setminus \bigcup\PPP)$ and $\PPP,\QQQ\in \left \Vert \PL\right\Vert$, whence $\PPP=\left\Vert\alpha\right\Vert$ and $\QQQ=\left\Vert\beta\right\Vert$ where $\alpha,\beta\in \PL$. By the converse inclusion, $\left\Vert\alpha\right\Vert^\pm\in  \{(\PPP,\QQQ)\mid (\PPP,\QQQ)$ world-incompatible$\}_{\PL} $, whence $\left\Vert \lnot \alpha\right\Vert=\wp(\left\Vert\top\right\Vert\setminus \bigcup\left\Vert\alpha\right\Vert)=\wp(\left\Vert\top\right\Vert\setminus \bigcup\PPP)=\QQQ=\left\Vert\beta\right\Vert$, so $(\PPP,\QQQ)=(\left\Vert\alpha\right\Vert, \left\Vert\beta\right\Vert)=(\left\Vert\alpha\right\Vert, \left\Vert\lnot \alpha\right\Vert)\in \left\Vert \PL\right\Vert^{\pm,\lnot}$. The same holds for all incompatibility notions considered in this section. Symmetric reasoning shows that the Burgess inclusion is similarly trivial if the second element always determines the first. In order to keep the term `Burgess theorem' meaningful, let us reserve it for a theorem which shows an inclusion of this sort for a property of pairs containing at least one pair whose first element does not determine the second, and at least one pair whose second element does not determine the first. Note that this allows for a Burgess theorem with respect to a property of pairs such that in each pair, either the first element determines the second or the second the first (as long as the first does not always determine the second, etc.). We will see such a theorem in Section \ref{neg:section:HS}.
%Similar results clearly also hold in the modal and first-order settings for team-based classical modal and first-order logic and their extensions with $\sim$.

%\aanch{We can make a similar distinction for the term `dual negation'. As mentioned in Section \ref{neg:section:preliminaries}, we use this term for all the notions of negation we consider that have bilateral semantics for the sake of convenience. However, in order to differentiate $\lnot$ as it appears in $\PL$ (and in the FO-fragments of D and IF) from dual negation proper  }

\subsection{Hawke and Steinert-Threlkeld's Logic} \label{neg:section:HS} \emph{Hawke and Steinert-Threlkeld's Logic} $\HS$ \cite{hawke} is essentially $\PL$ together with a operator $\diamondslash$ intended to represent epistemic modals such as the `might' in the sentence `It might be raining'. The anti-support clauses for the conjunction and disjunction are different from those defined in Section \ref{neg:section:classical_logics}; we will therefore use distinct symbols for the conjunction and disjunction of $\HS$ to avoid confusion.

%\todoa{Change conjunction disjunction symbols throughout.}

\begin{definition}[Syntax of $\HSBF$] \label{neg:def:hs_syntax}
The set of formulas of \emph{Hawke and Steinert-Threlkeld's Logic} $\HS$ is generated by:
        \begin{align*}
            \phi ::= p \sepp \bot \sepp \bnot \phi \sepp ( \phi \wedgedot \phi )  \sepp ( \phi \veedot \phi ) \sepp \diamondslash \phi
        \end{align*}
        where $p\in \mathsf{Prop}$.
\end{definition}
As before, let $\top:=\lnot \bot$. The support clauses for $p$, $\bot$, and $\lnot \phi$ are as before; the support clauses for $\wedgedot$ and $\veedot$ are the same as those for $\wedge$ and $\vee$, respectively. For $\diamondslash$, we have:
$$s\vDash \diamondslash \phi :\iff s \nDashv \phi. $$
The anti-support clause for $\lnot \phi$ is as before. For all formulas $\phi$ whose main connective is not $\lnot$:
$$s\Dashv \phi :\iff \forall t \subseteq s: [t\vDash \phi \implies t=\emptyset].$$
Observe that this anti-support clause is equivalent to our usual one if $\phi$ is $p$ or $\bot$.

We have made two alterations to the system as presented in \cite{hawke}: first, we have formulated the semantics in terms of support and anti-support conditions (Hawke and Steinert-Therlkeld's  semantics are essentially bilateral, but not formulated in terms of support/anti-support conditions). Second, we have added the constant $\bot$. This addition does not substantially change the logic---$\bot$ is strongly uniformly definable in the $\bot$-free fragment using any given $p\in \mathsf{Prop}$ in that $\bot\equiv^\pm p\wedgedot \lnot p$---but it allows us to secure a simple expressive completeness theorem for the logic (by giving us the means to express all team properties over the empty set of propositional variables).

Many of the propositional analogues of the results in Section \ref{neg:section:preliminaries} hold for this logic; we note in particular that the natural analogue of Fact \ref{neg:fact:ML_team_classical_correspondence} holds for the $\diamondslash$-free fragment, meaning that this fragment is expressively equivalent (in the usual positive sense) to $\PL$.
%and that (at least) a weakened version of Fact \ref{neg:fact:replacement} (replacement) holds: if $\chi[p]$ is not in the scope of any negations, then $\phi\equiv \psi$ implies $\chi[\phi/p]\equiv \chi[\psi/p]$. 
The expressive completeness theorem mentioned above follows from this together with the fact that we can uniformly define $\sim$ in $\HS$:
\begin{theorem}[Expressive completeness of $\HSBF$] \label{neg:theorem:HS_epxressive_completeness}
$\HS$ is expressively complete for the class of all team properties.
\end{theorem}
\begin{proof}
    Observe that $s\vDash \diamondslash \lnot \phi\iff s\nDashv \lnot\phi \iff s\nvDash \phi \iff s\vDash \sim \phi$. The result now follows from Theorem \ref{neg:theorem:PL_expressive_completeness} \ref{neg:theorem:PL_expressive_completeness_iii} together with the fact that the $\diamondslash$-free fragment of $\HS$ is expressively equivalent to $\PL$.
\end{proof}

Replacement of equivalents fails in negated contexts; note, for instance, that $\lnot p \equiv \lnot \diamondslash p$, but $\lnot \lnot p \equiv p \notequiv \diamondslash p \equiv \lnot \lnot \diamondslash p$. The incompatibility notion appropriate for $\HS$ is a kind of weaker two-sided version of down-set incompatibility:
\begin{definition}[(Either-sided) Down-set Incompatibility] \label{neg:def:two_side_down_set_incompatibility}
$\phi_0$ and $\phi_1$ are \emph{down-set incompatible (on either side) (E-D-I)} if $\phi_0$ is a down-set incompatibility of $\phi_1$ or $\phi_1$ is a down-set incompatibility of $\phi_0$.
\end{definition}
\begin{proposition} \label{neg:prop:HS_DI}
    For any $\phi\in \HS$, $\phi$ and $\lnot \phi$ are down-set incompatible (on either side).
\end{proposition}
\begin{proof}
    Case 1: $\phi$ is of the form $\lnot^n \psi$, where $\lnot^n$ denotes a string of symbols $\lnot$ of length $n$, and $n$ is even (possibly $0$); and the main connective of $\psi$ is not $\lnot$. Then $\phi\equiv \psi$ and $\lnot \phi\equiv \lnot \psi$. We have $s\vDash \lnot\psi $ iff $\forall t\subseteq s:[t\vDash \psi \implies t=\emptyset]$, and so $\lnot \psi$ is a down-set incompatibility of $\psi$; similarly $\lnot \phi$ is a down-set incompatibility of $\phi$.
    
    Case 2: $\phi$ is of the form $\lnot^n \psi$, where $\lnot^n$ and $\psi$ are as above, but $n$ is odd. Then $\phi\equiv \lnot \psi$ and $\lnot \phi \equiv \psi$. Similarly to the above, we now have that $\phi$ is a down-set incompatibility of $\lnot\phi$.
    \end{proof}

    Note that while the first element of a D-I pair determines the second, it is not the case that the first element of an E-D-I pair always determines the second (recall our example from above: $\lnot p \equiv \lnot \diamondslash p$, but $ \lnot \lnot p\equiv p \nequiv \diamondslash p \equiv \lnot\lnot \diamondslash p$, so an E-D-I pair with first element $\left\Vert\lnot p\right\Vert$ does not determine the second element) or that the second element always determines the first. It is, however, the case that for any given E-D-I pair, either the first element determines the second, or vice versa.

    Clearly D-I implies E-D-I, and it is also easy to see that, like D-I, E-D-I implies W-I. For other implications, see Figure \ref{neg:figure:inconsistency_relations}. Note that in a flat context, E-D-I is equivalent to F-I as well as D-I, so that $\PL$ is also \texttt{bicomplete} for E-D-I pairs.
    
    We now show Burgess and bicompleteness theorems for $\HS$ using E-D-I. This time our proof will not be analogous to Burgess'; rather, the result follows easily from the basic properties of the negation.

    \begin{theorem}[Burgess theorem for $\HSBF$] \label{neg:theorem:negation_HS}
    In $\HS$, the following are equivalent:
    \begin{enumerate}
        \item[(i)] \makeatletter\def\@currentlabel{(i)}\makeatother\label{neg:theorem:negation_HS_i} $\phi$ and $\psi$ are down-set incompatible (on either side).
        \item[(ii)] \makeatletter\def\@currentlabel{(ii)}\makeatother\label{neg:theorem:negation_HS_ii} There is a formula $\theta$ such that $\phi\equiv \theta$ and $\psi \equiv \lnot \theta$ (and $\mathsf{P}(\theta)
        \subseteq \mathsf{P}(\phi)\cup \mathsf{P}(\psi)$).
    \end{enumerate}
\end{theorem}
\begin{proof}
\ref{neg:theorem:negation_HS_ii} $\implies $ \ref{neg:theorem:negation_HS_i} by Proposition \ref{neg:prop:HS_DI}. \ref{neg:theorem:negation_HS_i} $\implies$ \ref{neg:theorem:negation_HS_ii}: If $\psi$ is D-I of $\phi$, let $\theta:=\phi \wedgedot \top$. Then $\theta\equiv \phi$ and $$s\vDash \lnot \theta \iff \forall t\subseteq s: [t\vDash \phi \wedgedot \top\implies t=\emptyset] \iff \forall t\subseteq s: [t\vDash \phi \implies t=\emptyset], $$
so that $\lnot \theta\equiv \psi$. Else if $\phi$ is D-I of $\psi$, let $\theta:=\lnot (\psi\wedgedot \top)$. Similarly to the previous case, $\theta\equiv \phi$ and $\lnot \theta\equiv \psi$.
\end{proof}

\begin{corollary}[Bicompleteness of $\HSBF$] \label{neg:coro:HS_negation_completeness}
$\HS$ is both bicomplete and \texttt{bicomplete} for E-D-I pairs.
\end{corollary}

Note that by Corollaries \ref{neg:coro:PL_negation_completeness} \ref{neg:coro:PL_negation_completeness_iii} and \ref{neg:coro:HS_negation_completeness}, we have $\left\Vert \PL(\sim)\right\Vert = \left\Vert \HS\right\Vert$, whereas $\left\Vert \PL(\sim)\right\Vert^{\pm,\sim} \neq \left\Vert \HS\right\Vert^{\pm,\lnot}$.

Interestingly, Hawke and Steinert-Threlkeld also consider a variant of $\HS$ which uses $\land$ and $\vee$ instead of $\wedgedot$ and $\veedot$, but they dismiss it since the pairs it produces are not $\bot$-incompatible (among other reasons): ``the altered system egregiously allows for the assertibility of contradictions" \cite{hawke}. ($\HS$ allows for a formula and its negation to be simultaneously true in the empty team, but in the variant this also holds for nonempty teams. An example from \cite{hawke}: $\lnot (\diamondslash p \land \diamondslash \lnot p)\land (\diamondslash p \land \diamondslash \lnot p)$, which is true in any team $s$ such that $s\nvDash p$ and $s\nvDash \lnot p$.) The bicompleteness of this variant is an open problem.

%$$s\vDash (\lnot \diamondslash p\land \top)\iff \forall t^+\subseteq s: t\nvDash \lnot \diamondslash p \iff \forall t^+\subseteq s: t\nDashv \diamondslash p\iff \forall t^+\subseteq s:\exists u^+\subseteq t: u\vDash \diamondslash p \iff \forall t^+\subseteq s:\exists u^+\subseteq t: u\nDashv p $$

\subsection{The Propositional Fragments of $\BSMLBF$ and $\BSMLIBF$} \label{neg:section:PLNE} Let $\PL(\NE)$ (the propositional fragment of $\BSML$) and $\PL(\NE,\vvee)$ (the propositional fragment of $\BSMLI$) be the extension of $\PL$ with $\NE$, and with $\NE$ and $\vvee$, respectively. In this section, we show that while $\PL(\NE,\vvee)$ is, like $\BSML$ and $\BSMLI$, \texttt{bicomplete} for ground-incompatible pairs, $\PL(\NE)$ is \texttt{bicomplete} for pairs conforming to a stronger notion of incompatibility. We also note that a variant of $\PL(\NE,\vvee)$ is both bicomplete and \texttt{bicomplete} for all pairs.

%It should be noted that the usual formulations of these logics (as well as those in Sections \ref{neg:section:classical_logics} and \ref{neg:section:propositional_dependence_logic}) in the literature do not include the dual negation; instead, one would typically incorporate a negation that may only occur in front of propositional variables, or one that may only occur in front of classical formulas. The resulting variants are expressively equivalent in the usual sense (of being able to express the same properties) to the present formulations, but since Burgess theorems do not hold for these variants, they are not expressively equivalent in the stronger sense of being able to express the same pairs of properties to the present ones. (This situation is different in $\PL$ or $\PL(\sim)$: since the dual negation may only precede classical formulas in these logics, it does not matter whether the dual negation or this alternative type of restricted negation is used---the different formulations are also equivalent in the stronger sense.)
Define, as before, $\Bot:=\bot \land \NE$ and $\bigvvee\emptyset:=\Bot$. The support and anti-support conditions are again the obvious propositional analogues of the modal ones. The propositional analogues of all results in Section \ref{neg:section:preliminaries} hold; we note in particular that each formula is equivalent to one in negation normal form, that replacement holds in positive contexts, and that we have the following propositional analogue of Theorem \ref{neg:theorem:BSML_expressive_completeness}:

\begin{theorem}[Expressive completeness of $\PLBF\mathbf{(\NE)}$ and $\PLBF\mathbf{(\NE,\vvee)}$] \cite{aknudstorp2024,yang2017}\ \label{neg:theorem:PLNE_expressive_completeness}
\begin{enumerate}
    \item[(i)] \makeatletter\def\@currentlabel{(i)}\makeatother\label{neg:theorem:PLNE_expressive_completeness_i} $\PL(\NE)$ is expressively complete for
    $\{\PPP \mid \PPP $ is convex and union closed$\}$.
    \item[(ii)] \makeatletter\def\@currentlabel{(ii)}\makeatother\label{neg:theorem:PLNE_expressive_completeness_ii} $\PL(\NE,\vvee)$ is expressively complete for the class of all team properties.
\end{enumerate}
\end{theorem}

As already noted in Section \ref{neg:section:BSML_theorems}, the $\BSMLI$ Burgess theorem makes no essential use of the modal operators, and so the propositional analogues of Theorem \ref{neg:theorem:separation} and Theorem \ref{neg:theorem:negation_bsml} (using the second proof given for this theorem) give us:

\begin{theorem}[Burgess theorem for $\PLBF\mathbf{(\NE,\vvee)}$] \label{neg:theorem:negation_PLNEVVEE}
    In $\PL(\NE,\vvee)$, the following are equivalent:
    \begin{enumerate}
        \item[(i)] \makeatletter\def\@currentlabel{(i)}\makeatother\label{neg:theorem:negation_PLNEVVEE_i} $\phi$ and $\psi$ are ground-incompatible.
        \item[(ii)] \makeatletter\def\@currentlabel{(ii)}\makeatother\label{neg:theorem:negation_PLNEVVEE_ii} There is a formula $\theta$ such that $\phi\equiv \theta$ and $\psi \equiv \lnot \theta$ (and $\mathsf{P}(\theta)=\mathsf{P}(\phi)\cup \mathsf{P}(\psi)$).
    \end{enumerate}
\end{theorem}

\begin{corollary}[Bicompleteness of $\PLBF\mathbf{(\NE,\vvee)}$] \label{neg:coro:PLVVEENE_negation_completeness}
    $\PL(\NE,\vvee)$ is both bicomplete and \texttt{bicomplete} for ground-incompatible pairs.
\end{corollary}

However, the propositional analogue of Theorem \ref{neg:theorem:negation_bsml} does not hold for $\PL(\NE)$. To show why, we recursively define the \emph{flattening}\footnote{The name `flattening' is taken from Väänänen, who defines an analogous notion for D in \cite[p. 42]{vaananen2007}. See also endnote \hyperref[neg:footnote:ground_team]{10}.} $\phi^f\in \PL$ of $\phi\in \PL(\NE)$ as follows: $p^f:=p$; $\bot^f:=\bot$; $\NE^f:=\top$; $(\lnot \phi)^f=\lnot \phi^f$; $(\phi\land \psi)^f:=\phi^f\land \psi^f$; and $(\phi\vee \psi)^f:=\phi^f\vee\psi^f$. Note that since $\phi^f$ is classical, we have $|\phi^f|=\llbracket \phi^f\rrbracket$. We require the following lemma for the sequel:

\begin{lemma} \label{neg:lemma:flattening_normal_form}
    Let $\phi'$ be the negation normal form of $\phi\in \PL(\NE)$. Then $\phi^f\equiv \phi'^f$.
\end{lemma}
\begin{proof}
    By induction on the complexity of $\phi$. The cases for the atoms, $\land $, and $\lor$ are all immediate. The cases for negated atoms are likewise immediate.

    If $\phi =\lnot \lnot \psi$, then $\phi'= \psi'$. By the induction hypothesis, $\psi^f\equiv\psi'^f$. We then have $\phi^f=\lnot \lnot \psi^f\equiv \psi^f\equiv \psi'^f =\phi'^f $.

    If $\phi=\lnot (\psi \land \chi)$, then $\phi'=(\lnot \psi)' \vee (\lnot \chi)'$. By the induction hypothesis, $(\lnot \psi)^f\equiv (\lnot \psi)'^f$ and $(\lnot \chi)^f\equiv (\lnot \chi)'^f$. We then have $$\phi^f=\lnot (\psi^f\land \chi^f)\equiv \lnot \psi^f \vee \lnot \chi^f=(\lnot \psi)^f \vee (\lnot \chi)^f\equiv (\lnot \psi)'^f \vee (\lnot \chi)'^f=\phi'^f.$$
    The case for $\phi=\lnot (\psi \vee \chi)$ is similar to that for $\phi=\lnot (\psi \land \chi)$.
\end{proof}

Using which we can show:
\begin{lemma} \label{neg:lemma:PLNE_flattening}
    For all $\phi\in \PL(\NE)$ and all $\mathsf{X}\supseteq \mathsf{P}(\phi)$, either $|\phi|_\mathsf{X}=|\phi^f|_\mathsf{X}$, or $\phi\equiv \Bot$.
\end{lemma}
\begin{proof}
    By induction on the complexity of $\phi$, which we may assume to be in negation normal form by Lemma \ref{neg:lemma:flattening_normal_form}. The base cases---those for atoms and negated atoms---are obvious; note in particular that $|\NE|=|\top|=|\NE^f|$ since for any $w$, $w\in \{w\}\vDash \NE$, and $|\lnot\NE|=|\bot|=|\lnot \top|=|(\lnot \NE)^f|$.

    Consider $\phi=\psi\land \chi$. If either $\psi\equiv \Bot$ or $\chi\equiv \Bot$, then $\psi\land \chi\equiv \Bot$. Otherwise we have both $|\psi|=|\psi^f|$ and $|\chi|=|\chi^f|$. If there is no $s$ such that $s\vDash \psi\land \chi$, then $\psi \land \chi \equiv\Bot$, and we are done; we may therefore assume there is $s$ with $s\vDash \psi \land \chi$. We now show $|\psi\land \chi|=|\psi|\cap |\chi|$. $|\psi\land \chi|\subseteq|\psi|\cap |\chi|$ is immediate; for the converse inclusion, let $w \in |\psi|\cap |\chi|$. Then $w\in t\vDash \psi$ and $w\in u\vDash \chi $. Let $s$ be such that $s\vDash \psi\land \chi$. By union closure, $s\cup t\vDash \psi$ and $s \cup u \vDash \chi$. Since $s\subseteq s\cup \{w\} \subseteq s\cup t$ and $s\subseteq s\cup \{w\} \subseteq s\cup u$, by convexity we have $s\cup \{w\} \vDash \psi$ and $s\cup \{w\} \vDash \chi$, so $w\in |\psi\land \chi|$. We then have
    $$|\psi\land \chi|= |\psi| \cap |\chi|=|\psi^f| \cap |\chi^f|=\llbracket \psi^f\rrbracket \cap \llbracket \chi^f\rrbracket =\llbracket \psi^f\land \chi^f\rrbracket=\llbracket (\psi \land \chi)^f\rrbracket=| (\psi \land \chi)^f|.$$

    Now consider $\phi=\psi\vee \chi$. If either $\psi\equiv \Bot$ or $\chi\equiv \Bot$, then $\psi\vee \chi\equiv \Bot$. Otherwise we have both $|\psi|=|\psi^f|$ and $|\chi|=|\chi^f|$. We show $|\psi\vee \chi|=|\psi|\cup |\chi|$. For $|\psi\vee \chi|\subseteq |\psi|\cup |\chi|$, if $w\in |\psi\vee \chi| $, then $w\in s\vDash \psi \vee \chi$. We have that $s=s_1 \cup s_2$ where $s_1\vDash \psi$ and $s_2\vDash \chi$, so $w\in s_1$ or $w\in s_2$. Either way, $w\in |\psi|\cup |\chi|$. To show the converse inclusion, let $w\in |\psi|\cup |\chi|$. Then $w\in s$ where $s\vDash \psi$ or $s\vDash \chi$; assume without loss of generality that $s\vDash \psi$. Since $\chi\nequiv \Bot$, there is $t$ such that $t\vDash \chi$. Then $s\cup t\vDash \psi\vee \chi$, so $s\in |\psi \vee\chi|$. We then have:
    $$|\psi\vee \chi|= |\psi| \cup |\chi|=|\psi^f| \cup |\chi^f|=\llbracket \psi^f\rrbracket \cup \llbracket \chi^f\rrbracket =\llbracket \psi^f\vee \chi^f\rrbracket=\llbracket (\psi \vee \chi)^f\rrbracket=| (\psi \vee \chi)^f|.$$
\end{proof}
Note that the $\BSML$-analogue of Lemma \ref{neg:lemma:PLNE_flattening} fails, since, for instance, extending the definition of flattening in the obvious way, we would have $|\Box \Bot|=|\bot|$ and $\Box \Bot \nequiv \Bot$, but $|(\Box \Bot)^f|=|\Box (\bot \land \top )|=|\Box \bot|\neq |\bot|$.
%(which, given the intended use of the lemma, would be formulated using bisimulation rather than equality between ground teams) fails, since, for instance, extending the definition of flattening in the obvious way, we would have $|\Box \Bot|\bisim (M,\emptyset)$ and $\Box \Bot \nequiv \Bot$, but $|(\Box \Bot)^f|\bisim|\Box (\bot \land \top )|\bisim|\Box \bot|\not\bisim (M,\emptyset)$ (for any model $M$).

Now assume for contradiction that the propositional analogue of Theorem \ref{neg:theorem:negation_bsml} holds for $\PL(\NE)$, i.e., that for any ground-incompatible $\phi$ and $\psi$, there is a $\theta$ such that $\phi\equiv \theta$ and $\psi\equiv \lnot \theta$. Let $\phi:=p $ and $\psi:=\lnot p  \land q$. Clearly these are ground-incompatible, so there is a $\theta$ such that $\phi\equiv \theta$ and $\psi\equiv \lnot \theta$. By Lemma \ref{neg:lemma:PLNE_flattening},
$$|\phi|=|\theta|=|\theta^f|=\llbracket \theta^f\rrbracket = \llbracket \top \rrbracket \setminus \llbracket \lnot \theta^f\rrbracket=| \top | \setminus | \lnot \theta^f|=| \top | \setminus | (\lnot \theta)^f|=|\top | \setminus | \lnot \theta|=|\top | \setminus | \psi|.$$
But this leads to contradiction, since, for instance, $w_{\lnot p,\lnot q}\notin |\phi|$ and $w_{\lnot p,\lnot q}\notin |\psi|$.

$\PL(\NE)$-formulas and their negations conform, instead, to the following incompatibility notion. 
\begin{definition}[Ground-complementariness (modulo $\Bot$)] \label{neg:def:ground-complementary}\
\begin{enumerate}
   \item[(i)] $\phi_0$ and $\phi_1$ are \emph{ground-complementary} if $|\phi_i|_\mathsf{X}=|\top|_\mathsf{X}\setminus |\phi_{1-i}|_\mathsf{X}$ for all $i\in \{0,1\}$ and $\mathsf{X}\supseteq \mathsf{P}(\phi)\cup \mathsf{P}(\psi)$.
    \item[(ii)] $\phi_0$ and $\phi_1$ are \emph{ground-complementary modulo $\Bot$} if (a) they are ground-complementary, (b) $\phi_0 \equiv \Bot$, or (c) $\phi_1\equiv \Bot$.
\end{enumerate}
\end{definition}
Similarly, let $\PPP$ and $\QQQ$ be \emph{ground-complementary (modulo $\Bot$)} if $\bigcup \PPP= |\top|\setminus \bigcup \QQQ$ (or $\PPP=\emptyset$ or $\QQQ=\emptyset$). %(Note that for flat formulas, this notion is equivalent to flat-incompatibility.) 
We have:

\begin{proposition} \label{neg:prop:PLNE_ground_complements}
    For any $\phi\in \PL(\NE)$, $\phi$ and $\lnot \phi$ are ground-complementary modulo $\Bot$.
\end{proposition}
\begin{proof}
    If $\phi\equiv \Bot$ or $\lnot \phi\equiv \Bot$, we are done. Otherwise by Lemma \ref{neg:lemma:PLNE_flattening}, $$|\phi|=|\phi^f|=\llbracket \phi^f\rrbracket =\llbracket \top \rrbracket \setminus \llbracket \lnot \phi^f \rrbracket =| \top |\setminus | \lnot \phi^f |=| \top |\setminus | (\lnot \phi)^f |=|\top|\setminus|\lnot \phi|.\qedhere$$
\end{proof}

%We can also think ground-complementariness (modulo $\Bot$) as a type of generalization of world-incompatibility to a setting in which downward closure (and the empty team property) fail. To see why, note that we have:
Note that we have the following connection between ground-complementariness and world-incompatibility:

\begin{fact} \label{neg:fact:downward_closed_ground_complements_world_incompatible}
    For downward-closed $\phi,\psi$: $\phi$ and $\psi$ are world-incompatible if and only if they are ground-complementary. Neither implication holds for all formulas.
\end{fact}
\begin{proof}
    Let $\phi$ and $\psi$ be world-incompatible, and let $w\in |\phi|$. Assume for contradiction that $w\in |\psi|$. We then have that $w\in s\vDash \phi$ and $w\in t\vDash \psi$, whence by downward closure, $\{w\}\vDash \phi$ and $\{w\}\vDash \psi$, contradicting world-incompatibility. So $|\phi|\subseteq |\top|\setminus |\psi|$. Now let $w\in |\top|\setminus |\psi|$. By world-incompatibility, either $\{w\}\vDash \phi$ or $\{w\}\vDash \psi$. The latter clearly implies $w\in |\psi|$, so we must have $\{w\}\vDash \phi$. Therefore, $|\top|\setminus |\psi|\subseteq |\phi|$ and so $|\phi|= |\top|\setminus |\psi|$.

    Conversely, let $\phi$ and $\psi$ be ground-complementary. Let $\{w\}\vDash \phi$ and assume for contradiction that $\{w\}\vDash\psi$. Then $w\in |\phi|\cap |\psi|$, contradicting ground-complementariness. So $\{w\}\vDash \phi \implies \{w\}\nvDash\psi$. Now let $\{w\}\nvDash\psi$. Assume for contradiction that $w\in |\psi|$. Then $w\in t\vDash \psi$, so by downward closure, $\{w\}\vDash \psi$, contradicting our assumption. So $w\notin |\psi|$, whence $w\in |\phi|$ by ground-complementariness.

    To see why world-incompatibility does not imply ground-complementariness in general, consider, for instance, the pair $\top$ and $(p\land \NE)\vee (\lnot p \land \NE)$, or the pair $q$ and $\lnot q\vvee ((p\land \NE)\vee (\lnot p \land \NE))$.

    To see why ground-complementariness does not imply world-incompatibility in general, consider, for instance, the pair $\bot$ and $(p\land \NE)\vee (\lnot p \land \NE)$, or the pair $q\land ((p\land \NE)\vee (\lnot p \land \NE))$ and $\lnot q\land ((p\land \NE)\vee (\lnot p \land \NE))$.
\end{proof}

Furthermore, for formulas with the empty team property, G-C modulo $\Bot$ and G-C are clearly equivalent. Therefore, given the fact above, G-C modulo $\Bot$, G-C, and W-I are all equivalent for downward-closed properties with the empty team property. In a flat context, each of these is also equivalent to each of F-I, D-I, and E-D-I, so that $\PL$ is also \texttt{bicomplete} for G-C/G-C mod $\Bot$ pairs. %The way in which G-C (mod $\bot$) generalizes W-I is different from the way in which I-C does so: if a negation conforms to I-C, the behavior of the negation is always classical with respect to worlds.  %The former can therefore be regarded as a type of team-based generalization of the latter. %The generalization is different in kind from that enacted by I-C: I-C always implies W-I so that, with a negation conforming to I-C, the behavior of the negation is always classical with respect to worlds. Under G-C mod $\Bot$, there is a special class 

We now show a Burgess theorem for $\PL(\NE)$ using ground-complementariness modulo $\Bot$. As in the preceding section, our proof will not be analogous to Burgess'; this time we derive our result by making a modification to the proof of Theorem \ref{neg:theorem:PLNE_expressive_completeness} \ref{neg:theorem:PLNE_expressive_completeness_i}. This proof makes use of the propositional analogues $\chi^{\mathsf{X}}_w:=\bigwedge\{p\mid p\in  \mathsf{X}, w\vDash p\} \land \bigwedge\{\lnot p\mid p\in \mathsf{X}, w\nvDash p\}$ and $\chi^{\mathsf{X}}_s:=\bigvee_{w\in s} \chi^{\mathsf{X}}_w$ of Hintikka formulas. As with Hintikka formulas, we have $w'\vDash \chi^{\mathsf{X}}_w \iff w'\restriction \mathsf{X}= w\restriction X$ (so that $w'\vDash \chi^{\mathsf{X}}_w \iff w'= w$ if $w,w'\in 2^\mathsf{X}$) and $s'\vDash \chi^{\mathsf{X}}_s \iff s'\restriction \mathsf{X}= t$ for some $ t\subseteq s\restriction X$, where $s\restriction \mathsf{X}:=\{v\restriction \mathsf{X}\mid v\in s\}$, etc. (so that $s'\vDash \chi^{\mathsf{X}}_s \iff s'\subseteq s$ if $s,s'\subseteq 2^\mathsf{X}$). Given a property $\PPP=\{t_1,\ldots, t_n\}$, define
$$\delta^\mathsf{X}_\PPP:=\bigwedge_{w_1\in t_1,\ldots,w_n\in t_n}(((\chi^\mathsf{X}_{w_1}\vee \ldots \vee \chi^\mathsf{X}_{w_n})\land \NE)\lor \top).$$
Theorem \ref{neg:theorem:PLNE_expressive_completeness} \ref{neg:theorem:PLNE_expressive_completeness_i} essentially follows from:
\begin{proposition} \cite{aknudstorp2024} \label{neg:prop:PLNE_normal_form}
    For any nonempty, convex, and union-closed property $\PPP=\{t_1,\ldots, t_n\}$ over $\mathsf{X}\subseteq \mathsf{Prop}$ (where $\mathsf{X}$ is finite),
    $\PPP=\left\Vert\bigvee_{t\in \PPP}\chi^\mathsf{X}_t \land \delta^\mathsf{X}_\PPP\right\Vert_\mathsf{X}.$
\end{proposition}
Using which we get:
\begin{proposition} \label{neg:prop:PLNE_dual_normal_form}
    For any nonempty, convex, and union-closed properties $\PPP=\{t_1,\ldots, t_n\}$, $\QQQ=\{s_1,\ldots, s_m\}$ over $\mathsf{X}\subseteq \mathsf{Prop}$ (where $\mathsf{X}$ is finite) such that $\PPP,\QQQ$ are ground-complementary modulo $\Bot$,
    $(\PPP,\QQQ)=\left\Vert(\bigvee_{t\in \PPP}\chi^\mathsf{X}_t \land \delta^\mathsf{X}_\PPP)\vee\lnot \delta^\mathsf{X}_\QQQ\right\Vert^\pm_{\mathsf{X}}.$
\end{proposition}
\begin{proof}
    It is easy to check that $\lnot \delta_\QQQ\equiv\bot$, whence $$(\bigvee_{t\in \PPP}\chi_t \land \delta_\PPP)\vee\lnot \delta_\QQQ\equiv (\bigvee_{t\in \PPP}\chi_t \land \delta_\PPP)\vee\bot\equiv \bigvee_{t\in \PPP}\chi_t \land \delta_\PPP,$$
    and therefore, by Proposition \ref{neg:prop:PLNE_normal_form}, $\left\Vert \bigvee_{t\in \PPP}\chi_t \land \delta_\PPP\vee\lnot \delta_\QQQ\right \Vert=\PPP$. It remains to show that $\left\Vert\lnot (( \bigvee_{t\in \PPP}\chi_t \land \delta_\PPP)\vee\lnot \delta_\QQQ)\right \Vert=\QQQ$. We have:
    $$\lnot ( (\bigvee_{t\in \PPP}\chi_t \land \delta_\PPP)\vee\lnot \delta_\QQQ)\equiv (\lnot \bigvee_{t\in \PPP}\chi_t \vee \lnot \delta_\PPP)\land \delta_\QQQ\equiv (\lnot \bigvee_{t\in \PPP}\chi_t \vee \bot)\land \delta_\QQQ\equiv \lnot \bigvee_{t\in \PPP}\chi_t \land \delta_\QQQ.$$
    Observe that $\bigcup \PPP= |\bigvee_{t\in \PPP}\chi_t |=\llbracket\bigvee_{t\in \PPP}\chi_t\rrbracket$, and similarly for $\QQQ$. Since $\PPP$ and $\QQQ$ are ground-complementary modulo $\Bot$, and since they are nonempty, $\bigcup \QQQ = |\top| \setminus \bigcup \PPP$. We therefore have:
    $$\llbracket \lnot \bigvee_{t\in \PPP}\chi_t \rrbracket =\llbracket\top \rrbracket \setminus \llbracket \bigvee_{t\in \PPP}\chi_t \rrbracket= \llbracket\top \rrbracket \setminus \bigcup \PPP= \bigcup \QQQ =\llbracket \bigvee_{s\in \QQQ}\chi_s\rrbracket,$$
    from which it follows by flatness that $\left\Vert \lnot \bigvee_{t\in \PPP}\chi_t \right\Vert=\left\Vert \bigvee_{s\in \QQQ}\chi_s\right\Vert$. We therefore have $$\lnot ( (\bigvee_{t\in \PPP}\chi_t \land \delta_\PPP)\vee\lnot \delta_\QQQ)\equiv \bigvee_{s\in \QQQ}\chi_s\land \delta_\QQQ,$$ so that by Proposition \ref{neg:prop:PLNE_normal_form}, $\left\Vert\lnot (( \bigvee_{t\in \PPP}\chi_t \land \delta_\PPP)\vee\lnot \delta_\QQQ)\right \Vert=\QQQ$.
\end{proof}
The Burgess and bicompleteness theorems now follow easily:
\begin{theorem}[Burgess theorem for $\PLBF(\NE)$] \label{neg:theorem:negation_PLNE}
    In $\PL(\NE)$, the following are equivalent:
    \begin{enumerate}
        \item[(i)] \makeatletter\def\@currentlabel{(i)}\makeatother\label{neg:theorem:negation_PLNE_i} $\phi$ and $\psi$ are ground-complementary modulo $\Bot$.
        \item[(ii)] \makeatletter\def\@currentlabel{(ii)}\makeatother\label{neg:theorem:negation_PLNE_ii} There is a formula $\theta$ such that $\phi\equiv \theta$ and $\psi \equiv \lnot \theta$ (and $\mathsf{P}(\theta)= \mathsf{P}(\phi)\cup \mathsf{P}(\psi)$).
    \end{enumerate}
\end{theorem}
\begin{proof}
    \ref{neg:theorem:negation_PLNE_ii} $\implies $ \ref{neg:theorem:negation_PLNE_i} by Proposition \ref{neg:prop:PLNE_ground_complements}. For \ref{neg:theorem:negation_PLNE_i} $\implies$ \ref{neg:theorem:negation_PLNE_ii}, if $\psi \equiv \Bot$, let $\theta:=\phi\land \lnot \Bot$. Then $\theta\equiv \phi \land (\lnot \bot \lor \lnot \NE)\equiv \phi \land (\top \lor \bot)\equiv \phi$ and $\lnot \theta\equiv \lnot \phi \vee \Bot\equiv \Bot\equiv \psi$. The case for $\phi \equiv \Bot$ is similar. If $\phi\nequiv \Bot$ and $\psi\nequiv \Bot$, let $\mathsf{X}:=\mathsf{P}(\phi)\cup \mathsf{P}(\psi)$, and let $\theta:=(\bigvee_{t\in \PPP}\chi^\mathsf{X}_t \land \delta^\mathsf{X}_\PPP)\vee\lnot \delta^\mathsf{X}_\QQQ$. By Proposition \ref{neg:prop:PLNE_dual_normal_form}, $\theta$ is as desired.
\end{proof}

\begin{corollary}[Bicompleteness of $\PLBF\mathbf{(\NE)}$]\ \label{neg:coro:PLNE_negation_completeness}
    $\PL(\NE)$ is bicomplete for $\{(\PPP,\QQQ) \mid$  $\PPP,\QQQ $ are convex and union closed; $\PPP,\QQQ$ are ground-complementary modulo $\Bot\}$ and hence \texttt{bicomplete} for pairs which are ground-complementary modulo $\Bot$.
\end{corollary}

We conclude this section by showing that a variant of $\PL(\NE,\vvee)$ is both bicomplete and \texttt{bicomplete} for all pairs. Let $\PL(\NE^*,\vvee)$ be $\PL(\NE,\vvee)$ with $\NE$ swapped out for an atom $\NE^*$ with the following support/anti-support clauses: $s\vDash \NE^*:\iff s\neq \emptyset$ and $s\Dashv \NE^*:\iff s\neq \emptyset $.\footnote{$\NE^*$ was introduced by Tomasz Klochowicz in unpublished work.} Clearly many of the properties of $\PL(\NE,\vvee)$ are preserved in this variant: we still have a negation normal form, replacement in non-negated contexts, and expressive completeness for all properties.

 We define $\top:=\lnot\top$, $\Bot^*:=\lnot((\NE^*\vvee \bot)\vee\top)$, and $\NE^{*-}:=\lnot ((\NE^* \land \Bot^*)\vvee \bot)$. One can then check that $\Bot^*\equiv^\pm\Bot $ and $\NE^{*-} \equiv^{\pm} \NE$.

 \begin{theorem}[Burgess theorem for $\PLBF\mathbf{(\NE^*,\vvee)}$]\ \label{neg:theorem:negation_NEstar} For any $\phi,\psi\in \PL(\NE^*,\vvee)$, there is a $\theta\in \PL(\NE^*,\vvee)$ such that $\phi\equiv \theta$ and $\psi\equiv \lnot \theta$.
 \end{theorem}
 \begin{proof}
     Let $\phi',\psi'$ be the negation normal forms of $\phi,\psi$ respectively. Let $\phi^*$ be the result of replacing each occurrence of $\NE^*$ in $\phi'$ with $\NE^{*-}$, and each occurrence of $\lnot \NE^*$ with $\bot$, and similarly for $\psi^*$. We then have, as in Lemma \ref{neg:lemma:negation_empty_team}, that $\phi^*\equiv \phi$, and that $\lnot \phi^*$ has the empty team property, and similarly for $\psi^*$. Now let $\phi_\top:=\lnot (\lnot \phi^*\vvee \lnot (\NE^* \vvee \bot))$, and similarly for $\psi_\top$. Then:
     \begin{align*}
         &\phi_\top&&=&&\lnot (\lnot \phi^*\vvee \lnot (\NE^* \vvee \bot)) &&\equiv && \phi^*\land (\NE^* \vvee \bot) && \equiv && \phi^*&&\equiv && \phi \\
         &\lnot \phi_\top&&\equiv&&\lnot \phi^*\vvee \lnot (\NE^* \vvee \bot) &&\equiv && \lnot \phi^*\vvee (\NE^* \land \top)&& \equiv && \top,&& &&
     \end{align*}
     where the final equivalence holds because $\lnot\phi^*$ has the empty team property. Similarly, $\psi_\top\equiv \psi$ and $\lnot \psi_\top\equiv \top$. Finally, let $\theta:=\phi_\top\vvee (\Bot^* \vee \lnot \psi_\top)$. Then:
          \begin{align*}
         &\theta&&=&&\phi_\top\vvee (\Bot^* \vee \lnot \psi_\top) &&\equiv && \phi_\top && \equiv && \phi\\
         &\lnot \theta&&=&&\lnot(\phi_\top\vvee (\Bot^* \vee \lnot \psi_\top)) &&\equiv && \lnot\phi_\top \land \lnot (\Bot^* \vee \lnot \psi_\top)  && \equiv && \top \land (\lnot \Bot^* \land  \psi_\top) \\
         &&&\equiv&&\top \land (\top \land  \psi_\top) &&\equiv && \psi_\top  && \equiv && \psi \tag*{\qedhere}
     \end{align*}
 \end{proof}
     In particular, for any $\phi$ there is a $\theta$ such that $\theta\equiv \phi\equiv \lnot \theta$.

     \begin{corollary}[Bicompleteness of $\PLBF\mathbf{(\NE^*,\vvee)}$]\ \label{neg:coro:NEstar_negation_completeness}
    $\PL(\NE^*,\vvee)$ is both bicomplete and \texttt{bicomplete} for all pairs.
\end{corollary}

\subsection{Propositional Dependence Logic} \label{neg:section:propositional_dependence_logic}

 \emph{Propositional dependence logic} $\PL(\con{\cdot)}$ is $\PL$ extended with the $n+1$-ary ($n\geq -1$) connectives $\dep{p_1\ldots, p_n}{q}$, where $p_1,\ldots, p_n,q\in \mathsf{Prop}$. The formulas $\dep{p_1\ldots, p_n}{q}$ are called \emph{dependence atoms}. In this section, we show that $\PL(\con{\cdot)}$ is \texttt{bicomplete} for ground-complementary pairs (unlike first-order dependence logic D, which is \texttt{bicomplete} for $\bot$-incompatible pairs).
 
 The support/anti-support conditions for dependence atoms are as follows:\\

    \begin{tabular}{p{2.4cm} p{0.6cm} p{8.5cm}}
        $s \vDash   \dep{p_1\ldots, p_n}{q}$  & $:\iff$& $\forall v,w\in s:[v\vDash p_i \iff w\vDash p_i$ for all $\forall 1\leq i\leq n]\implies$\\
        &&$[v\vDash q \iff w\vDash q]$\\
        $s \Dashv  \dep{p_1\ldots, p_n}{q}$ & $:\iff$ & $s=\emptyset$
       % $s \vDash   \dep{\phi_1\ldots, \phi_n}{\psi}^+$  & $:\iff$& $\forall t,u\subseteq s:[[t\vDash \phi_i \iff u\vDash \phi_i$ for all $\forall 1\leq i\leq n]$\\
       % &&and $[t\vDash \lnot \phi_i \iff u\vDash \lnot \phi_i$ for all $\forall 1\leq i\leq n]]\implies$\\
       % &&$[t\vDash \psi \iff u\vDash \psi]$\\
       % $s \Dashv  \dep{\phi_1\ldots, \phi_n}{\psi}^+$ & $:\iff$ & $s=\emptyset$\\
    \end{tabular}\\

In other words, a dependence atom $\dep{p_1\ldots, p_n}{q}$ is true/supported in a team $s$ if, whenever two valuations $v$ and $w$ in $s$ agree on the truth values of all the $p_i$, they also agree on the truth value of $q$---the values of $p_i,\ldots, p_n$ jointly determine the value of $q$ in any valuation in the team. $\dep{p_1\ldots, p_n}{q}$ is anti-supported in a team just in case the team is empty.\footnote{The anti-support clause for $\dep{p_1\ldots, p_n}{q}$ might seem puzzling at first glance. What is the reason, one might ask, for adopting the clause we use rather than, say, $s\Dashv \dep{p_1\ldots, p_n}{q}$ iff $s\nvDash \dep{p_1\ldots, p_n}{q}$? (Pietro Galliani, in recent unpublished work \cite{galliani2024}, considers a variant of D with a first-order version of this anti-support clause.) The clause we use has been adapted from the analogous first-order clause for D in \cite{vaananen2007}. This first-order clause was, in turn, chosen (it seems to me) for the following the reasons: it (unlike the variant clause above) preserves the equivalence between IF and D, as well as the empty team property and downward closure. Väänänen \cite[p. 24]{vaananen2007} also offers the following explanation (notation amended; adapted for the propositional setting): \\

\indent\begin{minipage}{.8\textwidth}
Why not allow $s\Dashv \dep{p_1\ldots, p_n}{q}$ for non-empty $s$? The reason is that if we negate ``for all $w,w'\in s$ such that $w\vDash p_1\iff w'\vDash p_1,$ $\ldots,$ $w\vDash p_n\iff w'\vDash p_n$, we have $w\vDash q\iff w'\vDash q$'', maintaining the analogy with [$s\Dashv p$ iff for all $w\in s: w\nvDash p$], we get ``for all $w,w'\in s$ we have $w\vDash p_1\iff w'\vDash p_1,$ $\ldots,$ $w\vDash p_n\iff w'\vDash p_n$ and $w\vDash q\not\iff w'\vDash q$'', which is only possible if $s=\emptyset$. 
\end{minipage} \\

Note also that Väänänen's clause preserves the ground-incompatibility and $\NE$-incompatibility of $\phi$ and $\lnot \phi$, whereas the variant above does not.}
%An atom $\dep{\phi_1\ldots, \phi_n}{\psi}^+$ is true in $s$ if, whenever two \emph{subteams} $t$ and $u$ of $s$ agree on the truth values of all the $\phi_i$, they also agree on the truth value of $\psi$---the values of $\phi_1,\ldots, \phi_n$ jointly determine the value of $\psi$ in any subteam.
We call a unary dependence atom $\con{p}$ a \emph{constancy atom}. $\con{p}$ expresses that the value of $p$ is constant in all valuations in a team.
%; $\con{\phi}^+$, that the value of $\phi$ is constant in all subteams of a team.

The propositional analogues of many results in Section \ref{neg:section:preliminaries} also hold in the context of $\PL(\con{\cdot)}$; we note in particular that each formula is equivalent to one in negation normal form, and that replacement holds in positive contexts. As with first-order dependence logic D, formulas of $\PL(\con{\cdot)}$ are downward closed and have the empty team property, which, by Fact \ref{neg:fact:downward-closed_incompatibility} and the definition of $\NE$-incompatibility, means that ground-incompatibility, $\bot$-incompatibility, and $\emptyset$-incompatibility are equivalent for these formulas. Formulas of $\PL(\con{\cdot)}$ need not be union closed: consider $\con{p}$. We have the following expressive completeness theorem:

\begin{theorem}[Expressive completeness of $\PLBF\mathbf{(\con{\cdot)}}$] \cite{yangvaananen2016}\ \label{neg:theorem:dep_expressive_completeness}
$\PL(\con{\cdot)}$ is expressively complete for
    $\{\PPP \mid \PPP $ is downward closed and has the empty team property$\}$.
\end{theorem}

Unlike in first-order dependence logic D, a Burgess theorem employing $\bot$-incompatibility does not hold for $\PL(\con{\cdot)}$. We can show this, as we did the analogous fact for $\PL(\NE)$, using flattenings of formulas. Define the flattening $\phi^f\in \PL$ of $\phi\in \PL(\con{\cdot)}$ as in Section \ref{neg:section:PLNE}, together with $\dep{p_1\ldots, p_n}{q}^f:=\top$. As in Section \ref{neg:section:PLNE}, it can be shown:
\begin{lemma} \label{neg:lemma:flattening_normal_form_dep}
    Let $\phi'$ be the negation normal form of $\phi\in \PL(\con{\cdot)}$. Then $\phi^f\equiv \phi'^f$.
\end{lemma}
And we have:
\begin{lemma} \label{neg:lemma:dep_flattening}
    For all $\phi\in \PL(\con{\cdot)}$ and all $\mathsf{X}\supseteq \mathsf{P}(\phi)$, $|\phi|_\mathsf{X}=|\phi^f|_\mathsf{X}$.
\end{lemma}
\begin{proof}
 By induction on the complexity of $\phi$, which we may assume to be in negation normal form by Lemma \ref{neg:lemma:flattening_normal_form_dep}. The base cases---those for atoms and negated atoms---are obvious; note that $|\dep{p_1\ldots, p_n}{q}|=|\top|=|\dep{p_1\ldots, p_n}{q}^f|$ since $w\in \{w\}\vDash \dep{p_1\ldots, p_n}{q}$ for any $w$, and $|\lnot\dep{p_1\ldots, p_n}{q}|=|\bot|=|(\lnot \dep{p_1\ldots, p_n}{q})^f|$.

 Consider $\phi=\psi\land \chi$. We show $|\psi\land \chi|=|\psi|\cap |\chi|$; the result will then follow as in the proof of Lemma \ref{neg:lemma:PLNE_flattening}. $|\psi\land \chi|\subseteq|\psi|\cap |\chi|$ is immediate; for the converse inclusion, let $w\in |\psi|\cap |\chi|$. Then $w\in t\vDash \psi$ and $w\in u\vDash \chi$. By downward closure, $\{w\}\vDash \psi$ and $\{w\}\vDash \chi$, so $w\in |\psi\land \chi|$.

 Now consider $\phi=\psi\vee \chi$. We show $|\psi\vee \chi|=|\psi|\cup |\chi|$; the result will then follow as in the proof of Lemma \ref{neg:lemma:PLNE_flattening}. $|\psi\vee \chi|\subseteq|\psi|\cup |\chi|$ follows as in the proof of Lemma \ref{neg:lemma:PLNE_flattening}. For the converse inclusion, let $w\in |\psi|\cup |\chi|$. Then $w\in s$ where $s\vDash \psi$ or $s\vDash \chi$; assume without loss of generality that $s\vDash \psi$. By the empty team property, $s=s\cup \emptyset\vDash \psi\vee \chi$, so $w\in |\psi\vee \chi|$.
\end{proof}

One can then show that a Burgess theorem employing $\bot$-incompatibility would lead to contradiction using essentially the same argument as that used for the analogous fact in Section \ref{neg:section:PLNE}. We could instead show a Burgess theorem for $\PL(\con{\cdot)}$ using ground-complementariness modulo $\Bot$, as we did for $\PL(\NE)$, but, as can be gleaned from the lemma above, or, alternatively, from the fact that formulas of $\PL(\con{\cdot)}$ have the empty team property, `modulo $\Bot$' is redundant here---we may simply use ground-complementariness. Note that since $\PL(\con{\cdot})$ is downward closed, by Fact \ref{neg:fact:downward_closed_ground_complements_world_incompatible} we could also equivalently use world-incompatibility.  %Say that $\phi$ and $\psi$ are \emph{ground-complementary} if for all $\mathsf{X}\supseteq \mathsf{P}(\phi)\cup \mathsf{P}(\psi)$, $|\phi|_\mathsf{X}=|\top|_\mathsf{X}\setminus |\psi|_\mathsf{X}$. This is equivalent to ground-complementariness modulo $\Bot$ for formulas with the empty team property. 
We have:

%We have, instead, the following strengthening of ground-complementary modulo $\Bot$. Say that $\phi$ and $\psi$ are \emph{ground-complementary} if for all $\mathsf{X}\supseteq \mathsf{P}(\phi)\cup \mathsf{P}(\psi)$, $|\phi|_\mathsf{X}=|\top|_\mathsf{X}\setminus |\psi|_\mathsf{X}$ and that $(\PPP,\QQQ)$ are \emph{ground-complementary} if $\bigcup \PPP= |\top|\setminus \bigcup \QQQ$.

\begin{proposition} \label{neg:prop:dep_ground_complements}
    For any $\phi\in \PL(\con{\cdot})$, $\phi$ and $\lnot \phi$ are ground-complementary.
\end{proposition}
\begin{proof}
    By Lemma \ref{neg:lemma:dep_flattening}; almost the same as the proof of Lemma \ref{neg:prop:PLNE_ground_complements}.
\end{proof}

We now show a Burgess theorem for $ \PL(\con{\cdot})$ using ground-complementariness. As with $\PL(\NE)$, we prove this by modifying the proof of the relevant expressive completeness theorem, Theorem \ref{neg:theorem:dep_expressive_completeness}. For further details on what follows, see \cite{yangvaananen2016}.

For a finite $\mathsf{X}\subseteq \mathsf{Prop}$, define $\gamma^{\mathsf{X}}_0:=\bot$, $\gamma^{\mathsf{X}}_1:= \bigwedge_{p \in \mathsf{X}}\con{p}$, and for $n\geq 2$, $\gamma^{\mathsf{X}}_n:=\bigvee_n \gamma_1$. Then for $s\subseteq 2^{\mathsf{X}}$, we have  $s\vDash \gamma^{\mathsf{X}}_n\iff |s|\leq n$, where $|s|$ is the size of $s$. Note that $\lnot \gamma^{\mathsf{X}}_0\equiv \top$ and $\lnot\gamma^{\mathsf{X}}_n\equiv \bot$ for $n\geq 1$.

    For a nonempty $s\subseteq 2^\mathsf{X}$, define $\xi^{\mathsf{X}}_s:=\gamma^{\mathsf{X}}_{|s|-1}\vee \chi^{\mathsf{X}}_{|\top|_{\mathsf{X}}\setminus s}$. It can be shown that for $t\subseteq 2^{\mathsf{X}}$, $t\vDash \xi^{\mathsf{X}}_s \iff s\not\subseteq t$. Note that if $|s|=1$, then $\xi^\mathsf{X}_s\equiv \chi^\mathsf{X}_{|\top|_\mathsf{X}\setminus s}$ so $|\xi^\mathsf{X}_s|_\mathsf{X}= |\chi^\mathsf{X}_{|\top|_\mathsf{X}\setminus s}|_\mathsf{X}$, and if $|s|> 2$, then $|\xi^\mathsf{X}_s|_\mathsf{X}=|\top|_\mathsf{X}$. Note also that if $|s|=1$ (say $s=\{w\}$), then $\lnot \xi^{\mathsf{X}}_s\equiv \top \land \lnot \chi^{\mathsf{X}}_{|\top|_{\mathsf{X}}\setminus s}\equiv \chi^{\mathsf{X}}_{s}=\chi^{\mathsf{X}}_w$ and if $|s|>1$, then $\lnot \xi^{\mathsf{X}}_s\equiv \bot$. 

    Theorem \ref{neg:theorem:dep_expressive_completeness} essentially follows from:

    \begin{proposition} \cite{yangvaananen2016} \label{neg:prop:dep_normal_form}
    For any downward-closed property $\PPP$ with the empty team property over $\mathsf{X}\subseteq \mathsf{Prop}$ (where $\mathsf{X}$ is finite),
    $\PPP=\left\Vert\bigwedge_{s\in\left\Vert\top \right\Vert_{\mathsf{X}}\setminus \PPP } \xi^{\mathsf{X}}_s \right\Vert_\mathsf{X}.$
\end{proposition}
For a property $\PPP$ over a finite $\mathsf{X}$, define $\PPP^1:=\{s\in \PPP\mid |s|=1\}$ and $\PPP^{>1}:=\{s\in \PPP\mid |s|>1\}$. We have the following corollary of Proposition \ref{neg:prop:dep_normal_form}:
    \begin{proposition} \label{neg:prop:dep_dual_normal_form}
    For any downward-closed properties $\PPP$ and $\QQQ$ with the empty team property over $\mathsf{X}\subseteq \mathsf{Prop}$ (where $\mathsf{X}$ is finite) such that $\PPP,\QQQ$ are ground-complementary,
    $$(\PPP,\QQQ)=\left\Vert\bigwedge_{s\in\left\Vert\top \right\Vert_{\mathsf{X}}\setminus \PPP } \xi^{\mathsf{X}}_s \vee \lnot \bigwedge_{s\in(\left\Vert\top \right\Vert_{\mathsf{X}}\setminus \QQQ)^{> 1} } \xi^{\mathsf{X}}_s \right\Vert^\pm_\mathsf{X}.$$
\end{proposition}
\begin{proof}
    It is easy to check that $\lnot \bigwedge_{s\in(\left\Vert\top \right\Vert\setminus \QQQ)^{> 1} } \xi_s\equiv \bot$, whence $$\bigwedge_{s\in\left\Vert\top \right\Vert\setminus \PPP } \xi_s \vee \lnot \bigwedge_{s\in(\left\Vert\top \right\Vert\setminus \QQQ)^{> 1} } \xi_s\equiv \bigwedge_{s\in\left\Vert\top \right\Vert\setminus \PPP } \xi_s \vee \bot\equiv \bigwedge_{s\in\left\Vert\top \right\Vert\setminus \PPP } \xi_s ,$$
    and therefore by Proposition \ref{neg:prop:dep_normal_form}, $\left\Vert\bigwedge_{s\in\left\Vert\top \right\Vert\setminus \PPP } \xi_s \vee \lnot \bigwedge_{s\in(\left\Vert\top \right\Vert\setminus \QQQ)^{> 1} } \xi_s\right\Vert =\PPP$. It remains to show $\left\Vert\lnot (\bigwedge_{s\in\left\Vert\top \right\Vert\setminus \PPP } \xi_s \vee \lnot \bigwedge_{s\in(\left\Vert\top \right\Vert\setminus \QQQ)^{> 1} } \xi_s)\right\Vert =\QQQ$. We first prove:
    \begin{align} \label{neg:claim}
        \text{For any team $s$, }\chi^{\mathsf{X}}_s\equiv \bigwedge_{w\in |\top|_\mathsf{X}\setminus s} \chi^{\mathsf{X}}_{|\top|_\mathsf{X}\setminus \{w\}}.
    \end{align}

    Proof of (\ref{neg:claim}): %We prove the equivalence for worlds; the equivalence for teams then follows by flatness.

    %$\chi^{\mathsf{X}}_s\vDash \bigwedge_{w\in |\top|_{\mathsf{X}}\setminus s} \chi^{\mathsf{X}}_{|\top|_{\mathsf{X}}\setminus \{w\}}$: Let $v\vDash \chi^{\mathsf{X}}_s$ (where for simplicity assume $v\in 2^{\mathsf{X}}$). Then $v\in s$. Let $w\in |\top|\setminus s$. Then $v\vDash \lnot \chi_w$. Observe that $\chi_w\equiv \bigwedge_{w'\in |\top|\setminus \{w\}} \lnot \chi_{w'} $. Therefore, $\lnot \chi_w\equiv \bigvee_{w'\in |\top|\setminus \{w\}}  \chi_{w'}=\chi_{|\top|\setminus \{w\}}$, and so $v\vDash \chi_{|\top|\setminus \{w\}}$. Since $w$ was arbitrary, $v\vDash \bigwedge_{w\in |\top|\setminus s} \chi_{|\top|\setminus \{w\}}$. 

    \begin{align*}
        \chi_s=\bigvee_{v\in s}\chi_v\equiv \bigwedge_{w\in |\top|\setminus s}\lnot \chi_w&\equiv\bigwedge_{w\in |\top|\setminus s} \bigvee_{w'\in |\top|\setminus \{w\}}  \chi_{w'} \equiv \bigwedge_{w\in |\top|\setminus s}  \chi_{|\top|\setminus \{w\}}  \tag*{$\dashv$}
        %\chi_s=\bigvee_{v\in s}\chi_v\equiv \bigwedge_{w\in |\top|\setminus s}\lnot \chi_w&\equiv\bigwedge_{w\in |\top|\setminus s}\lnot \bigwedge_{w'\in |\top|\setminus \{w\}} \lnot \chi_{w'} \\
        %&\equiv\bigwedge_{w\in |\top|\setminus s} \bigvee_{w'\in |\top|\setminus \{w\}}  \chi_{w'} \equiv \bigwedge_{w\in |\top|\setminus s}  \chi_{|\top|\setminus \{w\}}  \tag*{$\dashv$}
    \end{align*}
    
    We can now prove $\left\Vert\lnot (\bigwedge_{s\in\left\Vert\top \right\Vert\setminus \PPP } \xi_s \vee \lnot \bigwedge_{s\in(\left\Vert\top \right\Vert\setminus \QQQ)^{> 1} } \xi_s)\right\Vert =\QQQ$. We have:
    \begin{align*}
        &&&\lnot (\bigwedge_{s\in\left\Vert\top \right\Vert\setminus \PPP } \xi_s \vee \lnot \bigwedge_{s\in(\left\Vert\top \right\Vert\setminus \QQQ)^{> 1} } \xi_s)\\
        &\equiv&&\lnot \bigwedge_{s\in\left\Vert\top \right\Vert\setminus \PPP } \xi_s \land \bigwedge_{s\in(\left\Vert\top \right\Vert\setminus \QQQ)^{> 1} } \xi_s\\
        &\equiv &&\lnot (\bigwedge_{s\in(\left\Vert\top \right\Vert\setminus \PPP)^1 } \xi_s \land \bigwedge_{s\in(\left\Vert\top \right\Vert\setminus \PPP)^{>1} } \xi_s  )\land \bigwedge_{s\in(\left\Vert\top \right\Vert\setminus \QQQ)^{> 1} } \xi_s\\
        &\equiv &&(\lnot \bigwedge_{s\in(\left\Vert\top \right\Vert\setminus \PPP)^1 } \xi_s \vee \lnot \bigwedge_{s\in(\left\Vert\top \right\Vert\setminus \PPP)^{>1} } \xi_s  )\land \bigwedge_{s\in(\left\Vert\top \right\Vert\setminus \QQQ)^{> 1} } \xi_s\\
        &\equiv &&(\lnot \bigwedge_{s\in(\left\Vert\top \right\Vert\setminus \PPP)^1 } \xi_s \vee \bot)\land \bigwedge_{s\in(\left\Vert\top \right\Vert\setminus \QQQ)^{> 1} } \xi_s\\
        &\equiv &&\lnot \bigwedge_{s\in(\left\Vert\top \right\Vert\setminus \PPP)^1 } \xi_s \land \bigwedge_{s\in(\left\Vert\top \right\Vert\setminus \QQQ)^{> 1} } \xi_s\\
        &\equiv  &&\bigvee_{s\in(\left\Vert\top \right\Vert\setminus \PPP)^1 } \lnot\xi_s \land \bigwedge_{s\in(\left\Vert\top \right\Vert\setminus \QQQ)^{> 1} } \xi_s\\
        &\equiv  &&\bigvee_{s\in(\left\Vert\top \right\Vert\setminus \PPP)^1 } \chi_s \land \bigwedge_{s\in(\left\Vert\top \right\Vert\setminus \QQQ)^{> 1} } \xi_s\\
        &\equiv  &&\bigvee_{\{w\}\in(\left\Vert\top \right\Vert\setminus \PPP)^1 } \chi_w \land \bigwedge_{s\in(\left\Vert\top \right\Vert\setminus \QQQ)^{> 1} } \xi_s\\
        &\equiv  && \chi_{|\top| \setminus \bigcup\PPP} \land \bigwedge_{s\in(\left\Vert\top \right\Vert\setminus \QQQ)^{> 1} } \xi_s \tag{left-to-right by downward closure of $\PPP$}\\
        &\equiv  && \bigwedge_{w\in \bigcup \PPP}\chi_{|\top|\setminus\{w\}} \land \bigwedge_{s\in(\left\Vert\top \right\Vert\setminus \QQQ)^{> 1} } \xi_s \tag{\ref{neg:claim}}\\
        &\equiv  && \bigwedge_{w\in \bigcup \PPP}\xi_{\{w\}} \land \bigwedge_{s\in(\left\Vert\top \right\Vert\setminus \QQQ)^{> 1} } \xi_s \\
        &\equiv  && \bigwedge_{w\in |\top|\setminus \bigcup \QQQ}\xi_{\{w\}} \land \bigwedge_{s\in(\left\Vert\top \right\Vert\setminus \QQQ)^{> 1} } \xi_s \tag{$\PPP,\QQQ$ ground-complementary}\\
        &\equiv  && \bigwedge_{\{w\}\in (\left\Vert\top\right\Vert\setminus \QQQ)^1}\xi_{\{w\}} \land \bigwedge_{s\in(\left\Vert\top \right\Vert\setminus \QQQ)^{> 1} } \xi_s \tag{right-to-left by downward closure of $\QQQ$}\\
        &\equiv  && \bigwedge_{s\in\left\Vert\top \right\Vert\setminus \QQQ } \xi_s 
    \end{align*}
    The result now follows by Proposition \ref{neg:prop:dep_normal_form}.
\end{proof}

And so we have our Burgess and bicompleteness theorems:
\begin{theorem}[Burgess theorem for $\PLBF\mathbf{(\con{\cdot})}$] \label{neg:theorem:negation_dep}
    In $\PL(\con{\cdot})$, the following are equivalent:
    \begin{enumerate}
        \item[(i)] \makeatletter\def\@currentlabel{(i)}\makeatother\label{neg:theorem:negation_dep_i} $\phi$ and $\psi$ are ground-complementary.
        \item[(ii)] \makeatletter\def\@currentlabel{(ii)}\makeatother\label{neg:theorem:negation_dep_ii} There is a formula $\theta$ such that $\phi\equiv \theta$ and $\psi \equiv \lnot \theta$ (and $\mathsf{P}(\theta)= \mathsf{P}(\phi)\cup \mathsf{P}(\psi)$).
    \end{enumerate}
\end{theorem}
\begin{proof}
    \ref{neg:theorem:negation_dep_ii} $\implies $ \ref{neg:theorem:negation_dep_i} by Proposition \ref{neg:prop:dep_ground_complements}. For \ref{neg:theorem:negation_dep_i} $\implies$ \ref{neg:theorem:negation_dep_ii}, let $\mathsf{X}:=\mathsf{P}(\phi)\cup \mathsf{P}(\psi)$, and let $\theta:=\bigwedge_{s\in\left\Vert\top \right\Vert_{\mathsf{X}}\setminus \PPP } \xi^{\mathsf{X}}_s \vee \lnot \bigwedge_{s\in(\left\Vert\top \right\Vert_{\mathsf{X}}\setminus \QQQ)^{> 1} } \xi^{\mathsf{X}}_s$. By Proposition \ref{neg:prop:dep_dual_normal_form}, $\theta$ is as desired.
\end{proof}

\begin{corollary}[Bicompleteness of $\PLBF\mathbf{(\con{\cdot})}$] \label{neg:coro:dep_negation_completeness}
    $\PL(\con{\cdot})$ is bicomplete for $\{(\PPP,\QQQ) \mid $ $\PPP,\QQQ $ are downward closed and have the empty team property; $\PPP$ and $\QQQ$ are G-C/G-C mod $\Bot$/W-I$\}$ and hence \texttt{bicomplete} for pairs which are G-C/G-C mod $\Bot$/W-I.
\end{corollary}

Note that by Corollaries \ref{neg:coro:PL_negation_completeness} \ref{neg:coro:PL_negation_completeness_ii} and \ref{neg:coro:dep_negation_completeness}, we have $\left\Vert \INQB\right\Vert = \left\Vert \PL(\con{\cdot})\right\Vert$, but $\left\Vert \INQB\right\Vert^{\pm,\lnot_i} \neq \left\Vert \PL(\con{\cdot})\right\Vert^{\pm,\lnot}$.

To conclude, let us note that it is trivial to find a propositional logic which, like D in the first-order setting, is \texttt{bicomplete} for $\bot$-incompatible pairs, $\emptyset$-incompatible pairs, and ground-incompatible pairs: for instance, one could extend $\PL(\con{\cdot})$ with a constant $\theta_0$ such that $s\vDash \theta_0 \iff s\Dashv \theta_0\iff s=\emptyset$. It is not clear whether there is some interesting and non-\emph{ad hoc} propositional logic \texttt{bicomplete} for these pairs.

%% file: sections/conclusion.tex
\section{Concluding Remarks} \label{neg:section:conclusion}

In this paper, we proved analogues of Burgess' theorem for $\BSML$ and $\BSMLI$, for the propositional fragments of these logics, for Hawke and Steinert-Threlkeld's semantic expressivist logic for epistemic modals, as well as for propositional dependence logic with the dual negation. We saw that the notion of incompatibility employed to secure a Burgess theorem had to be adjusted according to the logic in question, formulated notions suitable for each of the logics we considered, and examined the relationships between these notions. We defined the notion of bicompleteness to succinctly describe our results: a logic is \texttt{bicomplete} with respect to a class of pairs of the relevant sort and conforming to a particular incompatibility notion just in case both a Burgess theorem employing that incompatibility notion as well as its converse hold for the logic. We also applied the notion of bicompleteness to logics which do not exhibit failure of determination of negated meanings by positive meanings, and gave an example of a logic which is \texttt{bicomplete} for all pairs.

The incompatibility notions and bicompleteness are interesting in their own right, and might warrant further study. We conclude with some final remarks on these notions, and on the interpretation and potential uses of our results.%and on some possible directions for further research.

We noted in Section \ref{neg:section:introduction} that Burgess characterized his theorem as concerning the degree of failure of contrariness/the dual negation (in H$_p$, and equivalently in IF and in D) to correspond to any operation on classes of models. One can think of our \texttt{bicompleteness} results and incompatibility notions as giving us a way of measuring this degree in different logics, with \texttt{bicompleteness} with respect to a weaker notion of incompatibility corresponding to a stronger degree of failure. 
%restricting our attention to downward-closed logics with the empty team property, 
On one end of the scale we have the classical $\PL$, which exhibits no such failure---the class of teams $\left\Vert \phi\right\Vert $ on which $\phi$ is true determines that $\left\Vert \lnot\phi\right\Vert $ on which $\lnot \phi$ is true  since, as per flat-incompatibility, $\left \Vert \lnot \phi\right \Vert=\wp(|\top|\setminus |\phi|)$---moreover, $\left\Vert \lnot \phi\right\Vert$ also determines $\left\Vert  \phi\right\Vert$: $\left \Vert \phi\right \Vert=\wp(|\top|\setminus |\lnot \phi|)$. To get \texttt{bicompleteness} for $\INQB$, we must go from the stronger flat-incompatibility to the weaker D-I $\left\Vert\lnot \phi\right\Vert=\{s\mid [t\subseteq s $ and $t\vDash \phi]\implies t=\emptyset\}$---$\left\Vert \phi\right\Vert $ still determines $\left\Vert \lnot\phi \right\Vert$, but the converse is no longer true. Further along the scale, we have mild failure of determination: for $\HS$ and E-D-I, it is not the case that $\left\Vert \phi\right\Vert $ always determines $\left\Vert \lnot \phi\right\Vert $ or vice versa, but it is always the case that one of the two determines the other. Propositional dependence logic violates determination in a different way with ground-complements $|\lnot\phi|=|\top|\setminus |\phi|$. Neither $\left\Vert \phi\right\Vert$ nor $\left\Vert \lnot\phi\right\Vert$ determines the other, but $\left\Vert \phi\right\Vert$ does determine the ground team $|\lnot \phi|$ of the negation, and similarly for $\left\Vert \lnot\phi\right\Vert$ and $| \phi|$. Even further along we have first-order dependence logic D, and we must now weaken ground-complements to $\bot$-incompatibility $\phi,\lnot \phi\vDash \bot$. Properties and ground teams are no longer determined, but we do at least know that $\left\Vert\phi\right\Vert$ and $\left\Vert\lnot \phi\right\Vert$ are disjoint modulo the empty team. At the very far end of the scale we have $\PL(\NE^*,\vvee)$, \texttt{bicomplete} for all pairs: $\left\Vert\phi\right\Vert$ does not constrain $\left\Vert\lnot \phi\right\Vert$ in any way.

We saw that dual negations in seemingly similar logics can conform to very different incompatibility notions.  %changes remarkably depending on the logic it is deployed in, with pairs $\phi,\lnot \phi$ in apparently very similar logics conforming to very different incompatibility notions.
For instance, whereas $\BSML$ is \texttt{bicomplete} for ground-incompatible pairs---pairs $(\PPP,\QQQ)$ in which $\PPP$ places only a very weak constraint on $\QQQ$ ($\bigcupdot\PPP\cap \bigcupdot\QQQ=\emptyset$)---simply removing the modalities gives us a logic ($\PL(\NE)$) which is \texttt{bicomplete} for pairs which are ground-complements modulo $\Bot$---pairs $(\PPP,\QQQ)$ such that, as long as $\PPP\neq \emptyset\neq \QQQ$, all the valuation-level information concerning $\PPP$ can be recovered from $\QQQ$ ($\bigcupdot\PPP=|\top|\setminus \bigcupdot \QQQ$). Similarly, whereas first-order dependence logic D is \texttt{bicomplete} for $\bot$-incompatible pairs, the natural propositional analogue $\PL(\con{\cdot})$ is \texttt{bicomplete} for ground-complementary pairs. 

On the one hand, these types of results are to be expected whenever the semantics of the negation $\lnot \phi$ of $\phi$ are defined in the bilateral way we have been observing, with the support conditions of $\lnot \phi$ not depending in a uniform manner on the positive semantic value $\left\Vert\phi\right \Vert$ of the argument $\phi$, but depending, rather, on the anti-support conditions specific to the main connective or atom of $\phi$. The pair properties that a given logic conforms to are a function of each of the anti-support clauses; if new atoms or connectives are added, their anti-support clauses must be configured to guarantee or bring about whatever properties are desired (for instance, switching out $\NE$, in $\PL(\NE,\vvee)$ with the atom $\NE^*$, whose semantics clearly do not conform to any of our incompatibility notions, allowed us to violate all of these notions).\footnote{The dependence of the semantics of the dual negation on the main logical symbol of the argument is obvious in the bilateral support/anti-support semantics we have been considering. It is also obvious in the team semantics for IF \cite{hodges1997} and D \cite{vaananen2007}, and Väänänen's team-based game semantics for D \cite[Section 5.2]{vaananen2007}---each of these systems features two clauses/rules for each connective or atom, one applied when the symbol in question appears in a positive context, the other when it appears in a negative context. However, it should be noted that this dependence is far less clear in the original assignment-based game semantics for IF \cite{hintikka1996} and D \cite[Section 5.3]{vaananen2007}. This type of semantics has a single game rule for each connective or atom, and the behavior of the negation is brought about by the interaction of these rules. The team semantics for these logics (whether bilateral or game-theoretical), accordingly, also make it easier to see why Burgess theorems hold, and enable one to easily adjust the properties of the dual negation as described in the main text. It seems that making similar adjustments would be more challenging in an assignment-, world-, or valuation-based game semantics, and it is not clear to me whether one can formulate (reasonable) world- or valuation-based game semantics for the logics we consider in this paper.} On the other hand, however, the results mentioned above demonstrate the surprising ways in which the anti-support conditions of different connectives can interact to violate aspects of determination which the conditions of any given connective do not violate on their own. If we add the $\BSML$-modalities to $\PL$, we have $\ML$, which, it is easy to see, is \texttt{bicomplete} for flat-incompatible pairs. If, instead of the modalities, we add $\NE$, we will thereby violate determination of $\left\Vert \lnot \phi \right\Vert$, but will still have determination of $| \lnot \phi |$. But if we add both the modalities and $\NE$, we have $\BSML$, which violates determination in a more radical way.

The standard positive expressive completeness theorems with respect to team-semantic closure properties (such as \ref{neg:theorem:BSML_expressive_completeness}, \ref{neg:theorem:PL_expressive_completeness}, \ref{neg:theorem:PLNE_expressive_completeness}, and \ref{neg:theorem:dep_expressive_completeness}) play an important role in the study of team logics, allowing for concise and tractable characterizations of these logics, and providing useful resources for axiomatization completeness proofs (see, e.g., \cite{yangvaananen2016,yang2017}) and the proofs of other properties such as uniform interpolation \cite{dagostino}. It would be interesting to investigate whether any similar technical application could be found for bicompleteness/Burgess theorems. At any rate, theorems of this type have some conceptual application---we have observed for instance, that we may think of them as providing a measure of the degree of failure of determination of negative meanings by positive ones, and that which incompatibility notions a given logic conforms to can provide insight into which types of information are clashing in a contradiction $\phi\land \lnot \phi$ of the logic. To that end, it may prove interesting and fruitful to establish bicompleteness results for more logics in the philosophical logic and formal semantics literature featuring bilateral negations, including logics which do not employ team semantics.

%% file: sections/bibliography.tex
\bibliography{bibb}